\documentclass{article}
\usepackage[a4paper, total={6in, 8in}]{geometry}
\usepackage[english]{babel}
\usepackage[T1]{fontenc}
\usepackage[utf8]{inputenc}
\usepackage{amsmath}
\usepackage{amsthm}
\usepackage{amssymb}
\usepackage{natbib}
\usepackage{mwe}
\usepackage[title]{appendix}
\usepackage{accents}
\usepackage{mathrsfs}
\usepackage{graphicx,mathtools,tikz,hyperref}
\usepackage{subcaption}
\usepackage{algorithm}
\usepackage{algpseudocode}
\usepackage{stackengine}

\usepackage{amsmath,amsfonts,bm}

\def\eqref#1{(\ref{#1})}

\def\1{\bm{1}}

\def\rmP{{\mathbf{P}}}

\def\vc{{\bm{c}}}

\def\ve{{\bm{e}}}

\def\vp{{\bm{p}}}

\def\vx{{\bm{x}}}
\def\vy{{\bm{y}}}
\def\vz{{\bm{z}}}
\def\vX{{\bm{X}}}

\def\mA{{\bm{A}}}
\def\mB{{\bm{B}}}
\def\mC{{\bm{C}}}

\def\mG{{\bm{G}}}

\def\mJ{{\bm{J}}}

\def\mL{{\bm{L}}}

\def\mP{{\bm{P}}}

\def\mS{{\bm{S}}}

\def\mY{{\bm{Y}}}

\DeclareMathAlphabet{\mathsfit}{\encodingdefault}{\sfdefault}{m}{sl}
\SetMathAlphabet{\mathsfit}{bold}{\encodingdefault}{\sfdefault}{bx}{n}

\DeclareMathOperator*{\argmin}{arg\,min}

\newcommand{\bI}{\textbf{I}}

\newcommand{\bA}{\textbf{A}}

\newcommand{\bX}{\textbf{X}}
\newcommand{\bY}{\textbf{Y}}

\newcommand{\la}{\langle}
\newcommand{\ra}{\rangle}

\newcommand{\norm}[1]{\left\lVert#1\right\rVert}
\newcommand{\dd}{\mathrm{d}}
\newcommand{\bQ}{\mathbf{Q}}
\numberwithin{equation}{section}
\usetikzlibrary{positioning}
\PassOptionsToPackage{hyphens}{url}\usepackage{hyperref}
\usepackage{cleveref}

\theoremstyle{plain}
\newtheorem{theorem}{Theorem}[section]
\newtheorem{proposition}[theorem]{Proposition}
\newtheorem{lemma}[theorem]{Lemma}

\newtheorem{corollary}[theorem]{Corollary}
\theoremstyle{definition}

\newtheorem{definition}[theorem]{Definition}

\theoremstyle{remark}
\newtheorem{remark}[theorem]{Remark}

\usepackage{enumitem}
\setlist[itemize]{itemsep=0pt, topsep=0pt}
\setlist[enumerate]{itemsep=-2pt, topsep=0pt}

\title{Gradient-adjusted underdamped Langevin dynamics for sampling}
\author{Xinzhe Zuo\thanks{Department of Mathematics, University of California, Los Angeles, CA, 90095.}, \and  Stanley Osher \footnotemark[1],\and  Wuchen Li \thanks{Department of Mathematics, University of South Carolina, Columbia, SC, 29208.}}
\date{}
\begin{document}
\maketitle
\begin{abstract}
Sampling from a target distribution is a fundamental problem with wide-ranging applications in scientific computing and machine learning. Traditional Markov chain Monte Carlo (MCMC) algorithms, such as the unadjusted Langevin algorithm (ULA), derived from the overdamped Langevin dynamics, have been extensively studied. From an optimization perspective, the Kolmogorov forward equation of the overdamped Langevin dynamics can be treated as the gradient flow of the relative entropy in the space of probability densities embedded with Wasserstein-2 metrics. Several efforts have also been devoted to including momentum-based methods, such as underdamped Langevin dynamics for faster convergence of sampling algorithms. Recent advances in optimizations have demonstrated the effectiveness of primal-dual damping and Hessian-driven damping dynamics for achieving faster convergence in solving optimization problems. Motivated by these developments, we introduce a class of stochastic differential equations (SDEs) called gradient-adjusted underdamped Langevin dynamics (GAUL), which add stochastic perturbations in primal-dual damping dynamics and Hessian-driven damping dynamics from optimization. We prove that GAUL admits the correct stationary distribution, whose marginal is the target distribution. The proposed method outperforms overdamped and underdamped Langevin dynamics regarding convergence speed in the total variation distance for Gaussian target distributions. Moreover, using the Euler-Maruyama discretization, we show that the mixing time towards a biased target distribution only depends on the square root of the condition number of the target covariance matrix. Numerical experiments for non-Gaussian target distributions, such as Bayesian regression problems and Bayesian neural networks, further illustrate the advantages of our approach over classical methods based on overdamped or underdamped Langevin dynamics. 
\end{abstract}


\section{Introduction}
Sampling from a target distribution is a long-standing quest and has numerous applications in scientific computing, including Bayesian statistical inference \cite{mackay2003information,robert1999monte,liu2001monte,gelman1995bayesian}, Bayesian inverse problems \cite{stuart2010inverse,idier2013bayesian,dashti2013bayesian,garbuno2020interacting}, as well as Bayesian neural networks \cite{welling2011bayesian,andrieu2003introduction,teh2016consistency,izmailov2021bayesian,mackay1995bayesian,neal2012bayesian}. In this direction, various algorithms have been developed to sample a target distribution $\pi\propto \exp(-f)$ for a given function $f:\mathbb R^d \to \mathbb R$, where $\pi$ is only known up to a normalization constant. In this area, a simple and popular algorithm is the unadjusted Langevin algorithm (ULA): 
\begin{equation}\label{eq:ula1}
  \vx_{k+1} = \vx_k - h \nabla f(\vx_k) + \sqrt{2h} \vz_k  \,,
\end{equation}
where $\vx_k\in \mathbb{R}^d$, $k$ is the iteration number, $f$ is assumed to be a differentiable function, $h>0$ is a step size, and $\vz_k$ is a $d$-dimensional random variable with independently and identically distributed (i.i.d) entries following standard Gaussian distributions. The ULA algorithm \cref{eq:ula1} comes from the forward Euler discretization of a stochastic differential equation (SDE) known as overdamped Langevin dynamics: 
\begin{equation}\label{eq:langevin}
    d\vx_t = - \nabla f(\vx_t) dt + \sqrt{2}d\mB_t \,,
\end{equation}
where $\vx_t\in\mathbb{R}^d$ and $\mB_t$ is a standard $d$-dimensional Brownian motion. Under some mild conditions on $f$, it has been shown that the SDE \cref{eq:SDE} has a unique strong solution $\{\vx_t,t\geq0\}$ that is a Markov process \cite{roberts1996exponential,meyn2012markov}. Moreover, the distribution of $\vx_t$ converges to the invariant distribution $\pi \propto \exp(-f)$ as $t\to \infty$. The asymptotic convergence guarantees of \cref{eq:ula1} have been established decades ago \cite{talay1990expansion,gelfand1991simulated,mattingly2002ergodicity}. In more recent years, non-asymptotic behaviors of \cref{eq:ula1} have also been explored by several works \cite{dalalyan2017further,dalalyan2017theoretical,durmus2017nonasymptotic,dalalyan2019user,cheng2018convergence,vempala2019rapid}.

An important result by \cite{jordan1998variational} states that the Kolmogorov forward equation of Langevin dynamics corresponds to the gradient flow of the relative entropy functional in the space of probability density functions with the Wasserstein-2 metric. This observation serves as a bridge between the sampling community and the optimization community by studying optimization problems in Wasserstein-2 space. {In the field of optimization, Nesterov's accelerated gradient \cite{nesterov1983method} is a first order algorithm for finding the minimum of a convex/strongly convex objective function $f$. The intuition is that Nesterov's method incorporates momentum into the updates. It is much faster than the traditional gradient descent method, in the sense that the convergence speed for convex functions is $\mathcal{O}(\frac{1}{k^2})$ where $k$ is the number of iterations compared to $\mathcal{O}(\frac{1}{k})$ for gradient descent. The convergence speed of Nesterov's method for $L$-smooth, $m$-strongly convex functions is $\mathcal{O}\big(\exp(-k/\sqrt{\kappa})\big)$ where $\kappa=L/m$ is the condition number of $f$ compared to $\mathcal{O}\big(\exp(-k/\kappa)\big)$ for gradient descent. By taking the step size to 0, one obtains a second-order ODE for Nesterov's method called the Nesterov's accelerated gradient flow or Nesterov's ODE \cite{su2016differential,attouch2019fast}. In recent years, one extends the gradient flow of the relative entropy into Nesterov's accelerated gradient flow \cite{su2016differential}, which is explored in \cite{wang2022accelerated,taghvaei2019accelerated,ma2021there} from different perspectives. } For the optimization in Wasserstein-2 space perspective, \cite{wang2022accelerated,taghvaei2019accelerated,chen2023accelerating} study a class of accelerated dynamics with depending on the score function, i.e., the gradient of logarithm of density function. This results in the approximation of a non-linear partial differential equation, known as the damped Euler equation \cite{CarrilloChoiTse2019_convergence}. In this case, the optimal choices of parameters for sampling a target distribution share similarities with the classical Nesterov's accelerated gradient flow. On the other hand, from a stochastic dynamics perspective, a line of research has been devoted to study the accelerated version of Langevin dynamics, known as the underdamped Langevin dynamics \cite{cao2023explicit,cheng2018underdamped,ma2021there,zhang2023improved}. As explained later in \Cref{sec:Nesterov_ode_UL}, the underdamped Langevin dynamics consists of a deterministic component and a stochastic component. The deterministic component exactly corresponds to the Nesterov's accelerated gradient flow. The marginal of invariant distribution in $x$-axis satisfies the target distribution.  However, the optimal choice of parameters in underdamped Langevin dynamics might not directly follow the classical Nesterov's method \cite{cheng2018underdamped}. 

Recently, \cite{zuo2023primal} proposed to use the primal-dual hybrid gradient (PDHG) method \cite{chambolle2011first,valkonen2014primal} to solve unconstrained optimization problems. The original PDHG method is designed for optimization problem with linear constraints. \cite{zuo2023primal} formulated the optimality condition $\nabla f(\vx) = 0$ of a strongly convex function $f$ into the solution of a saddle point problem 
$$\inf_{\vx\in\mathbb{R}^d} \sup_{\vp\in \mathbb{R}^d}\quad\la \nabla f(\vx), \vp \ra - \frac{\gamma}{2}\|\vp\|^2\,, $$ 
where $\gamma>0$ is a selected regularization parameter. They proceed by using the PDHG algorithm with appropriate preconditioners to solve the above saddle point problem. By taking the limit as the step size goes to zero, their algorithm yields a continuous-time flow, which is a second-order ordinary differential equation (ODE) called the primal-dual damping (PDD) dynamics. In particular, the PDD dynamic contains Nesterov's ODE \cite{su2016differential}. In other words, Nesterov's ODE is a special case of PDD dynamics. The PDD dynamics also shares similarities with the Hessian--driven damping dynamics that has been studied in recent years \cite{attouch2019fast,attouch2020first,attouch2021convergence}. 
The main difference between the PDD dynamics and the Nesterov's ODE is a second-order term $ \nabla^2 f(\vx) \dot \vx$ that appears in the former. This term is also presented in the Hessian driven damping dynamics. It has been observed that the PDD dynamics and the Hessian driven damping dynamics yield faster convergence towards the global minimum than the traditional gradient flow and Nesterov's ODE. Therefore, it is natural to extend the PDD dynamics and Hessian driven damping dynamics to SDEs for sampling a target distribution.

In this paper, we take inspirations from \cite{zuo2023primal,attouch2020first} to design a system of SDE called gradient-adjusted underdamped Langevin dynamics (GAUL) that resembles the primal-dual damping dynamics and the Hessian driven damping dynamics. Consider 
\begin{align}\label{eq:gaul}
\begin{pmatrix}
    d\vx_t \\ d\vp_t 
\end{pmatrix} &= \begin{pmatrix}
    -a\mC \nabla f(\vx_t)dt + \mC \vp_t dt \\
    -\nabla f(\vx_t) dt - \gamma \vp_t dt
\end{pmatrix} + \sqrt{\begin{pmatrix}
    2a\mC & \bI -\mC \\
    \bI -\mC & 2\gamma \bI 
\end{pmatrix}} \begin{pmatrix}
    d\mB^{(1)}_t \\ d\mB^{(2)}_t
\end{pmatrix}\,, 
\end{align}
for some constants $a,\gamma>0$, whose detailed choices will be explained later. $\mC$ is a preconditioner such that the diffusion matrix in front of the Brownian motion term is well-defined and positive semidefinite. And $\mB^{(i)}_t$ is a standard Brownian motion in $\mathbb{R}^d$ for $i=1,2$. The supercript on $\mB_t$ indicates that $\mB^{(1)}_t$ and $\mB^{(2)}_t$ are independent. We show that the stationary distribution GAUL \cref{eq:gaul} is the desired target distribution of the form $\frac{1}{Z}\exp(-f(\vx)-\|\vp\|^2/2)$. Noticeably, the $\vx$-marginal distribution is the target distribution $\pi$.  Additionally, we demonstrate that for a quadratic function $f$, GAUL achieves the exponential convergence and outperforms both overdamped and underdamped Langevin dynamics. A series of numerical examples are provided to demonstrate 
the advantage of the proposed method. 

To illustrate the main idea, we summarize main theoretical results into the following informal theorem.
\begin{theorem}[Informal]\label{thm:informal}
    Suppose that $f:\mathbb R^d \to \mathbb R^d$ is given by $f(\vx)=\frac{1}{2}\vx^T \Lambda \vx$ with a symmetric positive definite matrix $\Lambda\in\mathbb{R}^{d\times d}$ with eigenvalues $s_1 \geq s_2 \geq \ldots \geq s_d >0$. Let $\kappa=s_1/s_d$ be the condition number of matrix $\Lambda$.  And let $\mC = \bI$.
    \begin{itemize}
        \item[(1)]  Denote by $\rho_x(\vx,t)$ the law of $\vx_t$ driven by \cref{eq:gaul}, and $\pi(\vx) \propto \exp(-f(\vx))$ the target distribution. Let $a>0$, $\gamma=as_d + 2\sqrt{s_d}$. Then it takes at most $t = \mathcal{O}(\log(d/\delta))/(as_d + 2\sqrt{s_d})$ for the total variation distance between $\rho_x(\vx,t)$ and $\pi(\vx)$ to decrease to $\delta$. 
        \item [(2)] Denote by $\tilde\rho_x(\vx,k)$ the law of $\vx$ after $k$ iterations of the Euler-Maruyama discretization of \cref{eq:gaul}. Suppose $\sqrt{s_1}-\sqrt{s_d}\geq 2$, $a=1$, $\gamma=s_d + 2\sqrt{s_d}$ and consider the Euler-Maruyama discretization of \cref{eq:gaul} with step size $h=1/5s_1$. Then it takes at most $N=\mathcal{O}(\log(d/\delta)/(\kappa^{-1} + (\kappa s_1)^{-1/2})$ iterations for the total variation distance between $\tilde \rho_x(\vx,k)$ and $\tilde\pi(\vx)$ to decrease to $\delta$, where $\tilde\pi(\vx)$ is a biased target distribution given by \Cref{eq:x_element}. 
        \item[(3)] When taking $a=\frac{2}{\sqrt{s_1}-\sqrt{s_d}}$, $\gamma=as_d + 2\sqrt{s_d}$ and $h = \frac{1}{2(as_1 + \gamma)}$, we can improve the number of iterations in (2) to $N =\mathcal{O}(\sqrt{\kappa}\log(d/\delta)) $.
    \end{itemize}
\end{theorem}
The detailed version of \Cref{thm:informal} is given in \Cref{thm:cont_mix_time}, \Cref{thm:dis_mix_time} and \Cref{thm:choice_of_a}. It is worth noting that GAUL \cref{eq:gaul} reduces to underdamped Langevin dynamics when $a=0$ and $\mC=\bI$. Our theorem implies that in the Gaussian case, GAUL converges to the target measure faster than underdamped Langvein dynamics. In particular, we demonstrate that the Euler-Maruyama discretization admits a mixing time proportional to the square root of the condition number of covariance matrix. While this work primarily focuses on Gaussian distributions, our numerical experiments also explore non-log-concave target distributions in Bayesian linear regressions and Bayesian neural networks, which demonstrate potential advantages of GAUL over overdamped and underdamped Langevin dynamics. Extending these results to more general distributions and discretization schemes is an important future research direction. The choice of preconditoner $\mC$ is tricky as one needs to guarantee that the diffusion matrix in \cref{eq:gaul} is positive semidefinite. Therefore, we mainly focus on the case when $\mC =\bI$. We address on our results for $\mC\neq \bI$ in \cref{rmk:cont_mix_time_C} and \cref{rmk:dis_mix_C}. For $\mC=\bI$, \cite{li2022hessian} also explored dynamics \cref{eq:gaul}, which they called Hessian-Free High-Resolution (HFHR) dynamics. For this closely related work, we provide some comparisons later in \cref{rmk:comparison}. 

This paper is organized as follows. In \Cref{sec:prelim}, we review the connection between optimization methods and sampling dynamics, which leads to the construction of our proposed SDE called gradient-adjusted underdamped Langevin dynamics (GAUL). Our main results are presented in \Cref{sec:analysis_GAUL}, where we prove the exponential convergence of GAUL to the target distribution when the target measure follows a Gaussian distribution. We also study the Euler-Maruyama discretization of GAUL and prove its linear convergence to a biased target distribution. Lastly, in \Cref{sec:numerics}, we present several numerical examples to compare GAUL with both overdamped and underdamped Langevin dynamics.

\section{Preliminaries}\label{sec:prelim}

In this section, we briefly review the relation among Euclidean gradient flows, overdamped Langevin dynamics and Wasserstein gradient flows. We then draw the connection between the underdamped Langevin dynamics and Nesterov's ODEs. We next review primal-dual damping (PDD) flows \cite{zuo2023primal} and Hessian driven damping dynamics. Finally, we introduce a new SDE called gradient-adjusted underdamped Langevin dynamics (GAUL) for sampling, which resembles the PDD flow and the Hessian-driven damping dynamics with designed stochastic perturbations in terms of Brownian motions. 

\subsection{Gradient descent, unadjusted Langevin algorithms, and optimal transport gradient flows}
Let $f:\mathbb R^d\to \mathbb R$ be a differentiable convex function with $L$-Lipschitz gradient. The classical gradient descent algorithm for finding the global minimum of $f(\vx)$ is an iterative algorithm that reads:
\begin{equation}\label{eq:gd}
    \vx_{k+1} = \vx_k - h \nabla f(\vx_k)\,,
\end{equation}
where $h>0$ is the step size. When $f$ is convex and the step size is not too large, this algorithm converges at a rate of $\mathcal{O}(k^{-1})$. When $f$ is $m$-strongly convex, the same algorithm can be shown to converge at a rate of $\mathcal{O}\big((1-m/L)^k\big)$, if the step size is chosen appropriately. The gradient descent algorithm \cref{eq:gd} can be understood as the forward Euler time discretization of the gradient flow  
\begin{equation}\label{eq:gf}
    \dot \vx(t) = -\nabla f(\vx(t))\,, 
\end{equation}
where $\vx(t)$ describes a trajectory in $\mathbb R^d$ that travels in the direction of the steepest descent. Similar convergence results can be obtained for the gradient flow \cref{eq:gf}. When $f$ is convex, the gradient flow \cref{eq:gf} converges at a rate of $\mathcal{O}(t^{-1})$. When $f$ is assumed to be $m$-strongly convex, the gradient flow \cref{eq:gf} converges at a rate of $\mathcal{O}\big(\exp(-mt)\big)$. 

While the goal of optimization is to find the global minimum of $f$, the goal of sampling algorithm is to sample from a distribution of the form $\frac{1}{Z_1}\exp(-f(\vx))$, where the normalization constant $Z_1>0$ is assumed to be finite, i.e.,  $$Z_1=\int_{\mathbb{R}^d}e^{-f(x)}dx<+\infty.$$ The classical unadjusted Langevin algorithm (ULA) given in \cref{eq:ula1} is a simple modification to the gradient descent method. Recall that ULA is given by 
\begin{equation}\label{eq:ula_2}
     \vx_{k+1} = \vx_k - h \nabla f(\vx_k) + \sqrt{2h} \vz_k  \,,
\end{equation}
where $\vz_k$ is a $d$-dimensional standard Gaussian random variable and $h$ is the step size. We obtain \cref{eq:ula_2} from \cref{eq:gd} by adding a Gaussian noise term $\vz_k$ scaled by $\sqrt{2h}$. Similar to how \cref{eq:gd} can be viewed as the Euler discretization of \cref{eq:gf}, ULA \cref{eq:ula_2} represents the forward Euler discretization of the overdamped Langevin dynamics: 
\begin{equation}\label{eq:langevin_2}
     d\vx_t = - \nabla f(\vx_t) dt + \sqrt{2}d\mB_t \,,
\end{equation}
where $\mB_t$ is a standard $d$-dimensional Brownian motion. Denote by $\rho(\vx,t)$ the probability density function for $\vx_t$. Then the Kolmogorov forward equation (also known as the Fokker-Planck equation) of the overdamped Langevin dynamics \cref{eq:langevin_2} is given as  
\begin{equation}\label{eq:KFE}
    \frac{\partial \rho}{\partial t} = \nabla \cdot (\rho \nabla f) + \Delta \rho\,.
\end{equation}
Clearly,  $\pi(\vx) = \frac{1}{Z_1}\exp(-f(\vx))$ is a stationary solution of the Fokker-Planck equation \cref{eq:KFE}. In other words, note that $\nabla \pi=-\pi\nabla f$, then 
\begin{equation*}
  0=\partial_t\pi=\nabla \cdot (\pi \nabla f) + \Delta \pi=\nabla\cdot((\pi\nabla f+\nabla \pi))\,.  
\end{equation*}

In the literature, one can also study the gradient drift Fokker-Planck equation \cref{eq:KFE} from a gradient flow point of view. This means that equation \cref{eq:KFE} is a gradient flow in the probability space embedded with a Wasserstein-2 metric. We review some facts on a formal manner; see rigorous treatment in \cite{ambrosio2008gradient}. 

Define the probability space on $\mathbb{R}^d$ with finite second-order moment:
\[ \mathcal{P}(\mathbb{R}^d) = \left\{\rho(\cdot)\in C^{\infty}:~ \int_{\mathbb{R}^d} \rho(\vx)d\vx = 1, \; \int_{\mathbb{R}^d} |\vx|^2\rho(\vx)\ d\vx < \infty, \quad \rho(\cdot)\geq 0\right\}.\]
We note that $\mathcal{P}(\mathbb{R}^d)$ can be equipped with the $L_2$--Wasserstein metric $g_W$ at each $\rho \in \mathcal{P}(\mathbb{R}^d)$ to form a Riemannian manifold $(\mathcal{P}(\mathbb{R}^d),g_W)$. Let $\mathcal{F}: \mathcal{P}(\mathbb{R}^d)\to \mathbb R$ be an energy functional on $\mathcal{P}(\mathbb{R}^d)$. To be more precise, 
denote the Wassertein gradient operator of functional $\mathcal{F}(\rho)$ at the density function $\rho\in\mathcal{P}(\mathbb{R}^d)$, such that 
\begin{equation*}
  \mathrm{grad}_W \mathcal{F}(\rho) := -\nabla\cdot\Big(\rho \nabla \frac{\delta }{\delta \rho} \mathcal{F} (\rho) \Big)\,,  
\end{equation*}
where $\frac{\delta }{\delta \rho}$ is the $L_2$--first variation with respect to $\rho$. This yields that the gradient descent flow in the Wasserstein-2 space satsifies 
\begin{equation*}
     \frac{\partial \rho}{\partial t} = - \mathrm{grad}_W \mathcal{F}(\rho)=\nabla\cdot\Big(\rho \nabla \frac{\delta }{\delta \rho} \mathcal{F} (\rho) \Big) \,.
\end{equation*}
The above PDE is also named the {\em Wasserstein gradient descent flow}, in short Wasserstein gradient flows, which depend on the choices of the energy functionals $\mathcal{F}(\rho)$. 

An important example observed by \cite{jordan1998variational} is as follows. Consider the relative entropy functional, also named Kullback--Leibler(KL) divergence
\begin{equation*}
 \mathcal{F}(\rho):=\mathrm{D}_{\mathrm{KL}}(\rho\|\pi)= \int_{\mathbb{R}^d} \rho (\vx)\log\big(\frac{\rho(\vx)}{\pi(\vx)}\big)d\vx \,.   
\end{equation*}
One can show that the Fokker-Planck equation \cref{eq:KFE} is the gradient flow of the relative entropy in $(\mathcal{P}(\mathbb{R}^d),g_W)$. Upon recognizing $\frac{\delta }{\delta \rho}\mathrm{D}_{\mathrm{KL}}(\rho\|\pi) = \log\big(\frac{\rho}{\pi}\big)+1$, we obtain that \cref{eq:KFE} can be expressed as 
\begin{equation}
\begin{split}
    \frac{\partial \rho}{\partial t}  = &- \mathrm{grad}_W\mathrm{D}_{\mathrm{KL}}(\rho\|\pi) 
    =\nabla \cdot \Big(\rho \nabla \log\big(\frac{\rho}{\pi}\big) \Big)\\
    =&\nabla\cdot(\rho\nabla\log\rho)-\nabla\cdot(\rho \nabla\log\pi) \\
    =&\Delta \rho+\nabla \cdot (\rho \nabla f),
    \end{split}
\end{equation}
where we use facts that $\rho\nabla\log\rho=\nabla\rho$ and  $\nabla\log\pi=-\nabla f$. 

We note that the gradient of the logarithm of the density function, i.e. $\nabla\log\rho$, is often called the score function. The analysis of score functions are essential in understanding the convergence behavior of the Fokker-Planck equation \cref{eq:KFE} toward its invariant distribution; see related analytical studies in \cite{feng2024fisher}. 

\subsection{Nesterov's ODEs and underdamped Langevin dynamics}\label{sec:Nesterov_ode_UL}
Consider the problem of minimizing $f:\mathbb R^d \to \mathbb R$ for some convex function $f$ with $L$-Lipschitz gradient. \cite{nesterov1983method} proposed the following iterations:
\begin{subequations}
    \begin{align}
        \vx_{k+1} &= \vp_k - h \nabla f(\vp_k) \\
        \vp_{k+1} &= \vx_{k+1} + \gamma_k(\vx_{k+1}-\vx_k)\,,
    \end{align}
\end{subequations}
where $\gamma_k = (k-1)/(k-2)$. \cite{nesterov1983method} showed that the above method converges at a rate of $\mathcal{O}(k^{-2})$ instead of $\mathcal{O}(k^{-1})$ which is the convergence rate of the classical gradient descent method. If $f$ is further assumed to be $m$-strongly convex, then taking $h=1/L$ and $\gamma_k = \frac{1-\sqrt{\kappa}}{1+\sqrt{\kappa}}$ where $\kappa = L/m$, yields a convergence rate of $\mathcal{O}\big(\exp(-k/\sqrt{\kappa})\big)$. This is also considerably faster than gradient descent, which is $\mathcal{O}\big((1-\kappa^{-1})^k\big)$.  \cite{su2016differential} showed that the continuous-time limit of Nesterov's accelerated gradient method \cite{nesterov1983method} satisfies a second order ODE: 
\begin{equation}\label{eq:nag_ode}
    \Ddot{\vx} + \gamma_t \dot \vx + \nabla f(\vx) = 0 \,.
\end{equation}
If $f$ is a convex function, then $\gamma_t = 3/t$; if $f$ is a $m$-strongly convex function, then $\gamma_t =\gamma = 2\sqrt{m}$. As observed in \cite{maddison2018hamiltonian}, \cref{eq:nag_ode} can be formulated as a damped Hamiltonian system: 
\begin{equation}\label{eq:nag_sys_ode}
    \begin{pmatrix}
        \dot \vx \\
        \dot \vp 
    \end{pmatrix} = \begin{pmatrix}
        0\\ -\gamma_t \vp 
    \end{pmatrix} + \begin{pmatrix}
        0&\bI \\
        -\bI & 0 
    \end{pmatrix}\begin{pmatrix}
        \nabla_x H(\vx,\vp) \\
        \nabla_p H(\vx,\vp) 
    \end{pmatrix} = \begin{pmatrix}
        0&\bI \\
        -\bI & -\gamma_t \bI 
    \end{pmatrix}\begin{pmatrix}
        \nabla_x H(\vx,\vp) \\
        \nabla_p H(\vx,\vp) 
    \end{pmatrix}  \,, 
\end{equation}
where the Hamiltonian function is defined as $H(\vx,\vp)=f(\vx) + \|\vp\|^2/2$, $\vp\in\mathbb R^d$. On the other hand, the underdamped Langevin dynamics for sampling $\Pi(\vx,\vp) \propto \exp(-f(\vx) - \|\vp\|^2/2) $ is given by the system of SDE:  
\begin{align}
    d\vx_t &= \vp_t dt,  \nonumber \\
    d\vp_t &= -\nabla f(\vx_t) dt - \gamma_t \vp_t dt + \sqrt{2\gamma_t} d\mB_t, \nonumber 
\end{align}
where $\gamma_t$ is some damping parameter, and $\mB_t$ is a $d$-dimensional standard Brownian motion. This can be reformulated as 
\begin{equation}\label{eq:ul_sde}
    \begin{pmatrix}
        d\vx_t \\
        d\vp_t 
    \end{pmatrix} = \begin{pmatrix}
        0 & \bI \\
        -\bI & -\gamma_t \bI 
    \end{pmatrix}\begin{pmatrix}
        \nabla_x H(\vx,\vp) \\
        \nabla_p H(\vx,\vp) 
    \end{pmatrix}dt + \begin{pmatrix}
        0&0\\
        0&\sqrt{2\gamma_t}\bI
    \end{pmatrix} d\mB_t \,,
\end{equation}
where $\mB_t$ is a $2d$-dimensional standard Brownian motion. Observe that by adding a suitable Brownian motion term (the last term on the right hand side of \cref{eq:ul_sde}) to \cref{eq:nag_sys_ode}, Nesterov's accelerated gradient method for convex optimization becomes an algorithm for sampling $\Pi(\vx,\vp) =\frac{1}{Z} \exp(-f(\vx) - \|\vp\|^2/2) $, where $Z:=\int_{\mathbb{R}^{2d}}\exp(-f(\vx) - \|\vp\|^2/2) d\vx d\vp<+\infty$ is a noramlization constant. Moreover, the $\vx$-marginal of $\Pi(\vx,\vp)$ is simply $\pi(\vx)=\frac{1}{Z_1}\exp(-f(\vx))$ up to a normalizing constant $Z_1:=\int_{\mathbb{R}^{2d}}\exp(-f(\vx) - \|\vp\|^2/2) d\vx d\vp<+\infty$. Therefore, \cref{eq:ul_sde} can be used to sample distributions of the form $\exp(-f(\vx))/Z_1$. We postpone the proofs in terms of Fokker-Planck equations and there invariant distributions in Proposition \ref{prop:decomposition} and \ref{prop:pi}. 

\subsection{Primal-dual damping dynamics and Hessian driven damping dynamics}
Recently, \cite{zuo2023primal} proposed to solve an unconstrained strongly convex optimization problem using the PDHG method by considering the saddle point problem 
$$
\inf_{\vx\in\mathbb R^d} \sup_{\vp \in \mathbb R^d} \quad\la \nabla f(\vx), \vp \ra - \frac{\gamma}{2}\|\vp\|^2\,,
$$
where $\gamma$ is a damping parameter, and $f:\mathbb R^d \to \mathbb R$ is $m$-strongly convex. Note that the saddle point $(\vx^*,\vp^*)$ for the above inf-sup problem satisfies $\nabla f(\vx^*) = \vp^* = 0$. Then the primal-dual damping (PDD) algorithm \cite{zuo2023primal} admits the following iterations
\begin{align*}
    \vp_{k+1} & = \frac{1}{1+\tau_1 \gamma}\vp_k + \frac{\tau_1}{1+\tau_1 \gamma}\nabla f(\vx_k)\,, \\
    \tilde \vp_{k+1} &= \vp_{k+1} + \omega (\vp_{k+1} - \vp_k)\,,\\
    \vx_{k+1} &= \vx_k - \tau_2 \mC(\vx_k) \tilde \vp_{k+1} \,,
\end{align*}
where $\tau_1,\tau_2>0$ are dual and primal step sizes, $\omega>0$ is an extrapolation parameter, and $\mC\in\mathbb{R}^{d\times d}$ is a preconditioning positive definite matrix that could depend on $\vx_k$ and $t$. The continuous-time limit of the PDD algorithm can be obtained by letting $\tau_1, \tau_2 \to 0$ while keeping $\tau_1 \omega \to a$ for some $a>0$.  This yields a second-order ODE called the PDD flow: 
\begin{equation}\label{eq:pdd_C}
    \Ddot{\vx} + \Big(\gamma  + a \mC  \nabla^2 f(\vx) - \dot \mC \mC^{-1} \Big) \dot \vx + \mC  \nabla f(\vx) = 0 \,.
\end{equation}
In the case when $\mC$ is constant, \cref{eq:pdd_C} reads
\begin{equation}\label{eq:pdd_C_no_t}
    \Ddot{\vx} + \Big(\gamma  + a \mC  \nabla^2 f(\vx) \Big) \dot \vx + \mC  \nabla f(\vx) = 0 \,.
\end{equation}
And when $\mC = \bI$, the PDD flow simplifies to 
\begin{equation}\label{eq:hessian_driven_damping}
    \Ddot{\vx} + \gamma \dot \vx +  a \nabla^2 f(\vx) \dot \vx  + \nabla f(\vx) = 0\,.
\end{equation}
This corresponds to the Hessian driven damping dynamic \cite{attouch2020first} when $\gamma = 2\sqrt{m}$. The terminology `Hessian driven damping' comes from the Hessian term $\nabla^2 f(\vx) \dot \vx$ in \cref{eq:hessian_driven_damping}, which is controlled by a constant $a\geq 0$. When $a=0$, equation \cref{eq:hessian_driven_damping} reduces to Nesterov's ODE \cref{eq:nag_ode}. As in dynamics \cref{eq:nag_sys_ode}, we can express equation \cref{eq:pdd_C} as  
\begin{equation}\label{eq:pdd_system}
    \begin{pmatrix}
        \dot \vx \\
        \dot \vp 
    \end{pmatrix} = \begin{pmatrix}
        -a \mC & \mC \\
        (\gamma a-1)\bI & -\gamma \bI 
    \end{pmatrix}\begin{pmatrix}
        \nabla_x H(\vx,\vp) \\
        \nabla_p H(\vx,\vp) 
    \end{pmatrix} \,,
\end{equation}
where as before the Hamiltonian function is $H(\vx,\vp) = f(\vx) + \|\vp\|^2/2$. Note that one of the key differences between \cref{eq:nag_sys_ode} and \cref{eq:pdd_system} is that the top left block of the preconditioner matrix is nonzero in \cref{eq:pdd_system}, which gives rise to the Hessian damping term $\nabla^2 f(\vx) \dot \vx$. Throughout this paper, we focus on the dynamical system \cref{eq:pdd_system}.

\subsection{Gradient-adjusted underdamped Langevin dynamics}
We design a sampling dynamics that resembles the PDD flow and the Hessian driven damping with stochastic perturbations by Brownian motions. Our goal is still to sample a distribution proportional to $\exp(-f(\vx))$ for some $f:\mathbb R^d \to \mathbb R$. Let $H(\vx,\vp) = f(\vx) + \norm{\vp}^2/2 $. And denote by $\vX=(\vx,\vp)\in \mathbb R^{2d} $. We consider the following SDE. 
\begin{equation}\label{eq:SDE}
    d\vX_t = -\mathbf{Q}\nabla H(\vX_t)dt + \sqrt{2\,\mathrm{sym}(\mathbf{Q})}d\mB_t\,,
\end{equation}
where $\bQ \in \mathbb R^{2d\times 2d}$ is of the form 
\begin{equation} \label{eq:general_Q}
    \bQ = \begin{pmatrix}
        a\mC & -\mC \\
        \bI & \gamma\bI
    \end{pmatrix}\,,
\end{equation}
for some constant $a,\gamma\in \mathbb{R}$, and symmetric positive definite $\mC\in \mathbb R^{d\times d}$. $ \nabla H(\vX_t)=(\nabla_x H(\vX_t),\nabla_p H(\vX_t))^T$. And $\mathrm{sym}(\bQ)=\frac{1}{2}(\bQ + \bQ^T)$ is the symmetrization of $\bQ$. We assume that $\mathrm{sym}(\mathbf{Q})$ is positive semidefinite.

Throughout this paper, we will limit our discussion to $a,\gamma \geq 0$. $\mB_t$ is a $2d$-dimensional standard Brownian motion. Observe that when $a=0$, \cref{eq:SDE} reduces to underdamped Langevin dynamics \cref{eq:ul_sde}. When $a>0$, \cref{eq:SDE} has an additional gradient term $a \mC \nabla f(\vx_t)$ in the $d\vx_t$ equation. Thus, we call \cref{eq:SDE} gradient-adjusted underdamped Langevin dynamics. Let us examine the probability density function $\rho(\vX,t)$ of the diffusion governed by \cref{eq:SDE}. This is described by the following Fokker-Planck equation:
\begin{equation}\label{eq:FPK_general}
    \frac{\partial \rho}{\partial t} = \nabla \cdot (\bQ \nabla H \rho) + \sum_{i,j=1}^{2d} \frac{\partial^2}{\partial X_i \partial X_j} (Q_{ij} \rho)\,. 
\end{equation}
We assume that $f$ is differentiable and $\nabla f$ is a smooth Lipschitz vector field.  This ensures that the Fokker-Planck equation \cref{eq:FPK_general} has a smooth solution when $t>0$ for a given initial condition, such that $\rho(\vX, 0)\geq 0$ and $\int_{\mathbb{R}^{2d}}\rho(\vX, 0) d\vX=1$. 

Denote by $\Pi(\vX) = \frac{1}{Z}e^{-H(\vX)}$, where $Z$ is a normalization constant such that $\Pi(\vX)$ integrates to one on $\mathbb R^{2d}$. We show that $\Pi(\vX)$ is the stationary distribution of \cref{eq:FPK_general}. First, we have the following decomposition for \cref{eq:FPK_general}. 
\begin{proposition}[\cite{feng2024fisher} Proposition 1]\label{prop:decomposition}
    The Fokker-Planck equation \cref{eq:FPK_general} can be decomposed as 
    \begin{equation}
     \frac{\partial \rho}{\partial t} = \nabla \cdot \Big(\rho \,\mathrm{sym}(\bQ) \nabla \log \frac{\rho}{\Pi} \Big) + \nabla \cdot (\rho \Gamma) \label{eq:FPK_decomp} \,, 
\end{equation}
where 
\begin{equation}\label{eq:FPK_Gamma}
\begin{split}
         \Gamma(\vX):=& \mathrm{sym}(\bQ) \nabla \log(\Pi(\vX)) + \bQ\nabla H(\vX)\\
         =& \frac{1}{2}(\bQ-\bQ^T)\nabla H(\vX)\,.
         \end{split}
\end{equation}
In particular, the following equality holds:
\begin{equation*}
   \nabla\cdot(\Pi(\vX) \Gamma(\vX))=0. 
\end{equation*}

\end{proposition}
The proof is presented in \Cref{sec:post_proof}. Observe that the first term on the right-hand side of \cref{eq:FPK_decomp} is a Kullback–Leibler (KL) divergence functional that appears in a Fokker-Planck equation associated with the overdamped Langevin dynamics  \cref{eq:KFE}. The second term is due to the fact that the drift term $-\bQ \nabla H$ in \cref{eq:SDE} is a non-gradient vector field.
\begin{proposition}\label{prop:pi}
    $\Pi(\vX)$ is a stationary distribution for \cref{eq:FPK_general}. 
\end{proposition}
The proof is based on a straightforward calculation: When $\rho=\Pi$, we have $\nabla \cdot (\rho\Gamma)=0$, and therefore $\frac{\partial \rho}{\partial t}=0$. For completeness, we have included this calculation in \Cref{sec:post_proof}. This shows that $\Pi(\vX)$ is indeed the stationary distribution of $\cref{eq:FPK_general}$. Like the underdamped Langevin dynamics, the $\vx$-marginal of the stationary distribution is $\exp(-f(\vx))$ up to some normalization constant. Therefore, \cref{eq:SDE} can be used for sampling $\frac{1}{Z_1}\exp(-f(\vx))$ by first jointly sampling $\vX=(\vx,\vp)$ and then taking out the $\vx$-marginal. 

\begin{remark}
GAUL can also be viewed as a preconditioned overdamped Langevin dynamics on the space of $(\vx,\vp)\in \mathbb R^{2d}$. Designing optimal preconditioning matrix and optimal diffusion matrix have been studied in literature; see \cite{casas2022split,bennett1975mass,girolami2011riemann,leimkuhler2018ensemble,goodman2010ensemble,chen2023gradient, lelievre2024optimizing,lelievre2013optimal}. In particular, \cite{lelievre2024optimizing} considered the necessary condition on the optimal diffusion coefficient by studying the spectral gap of the generator assosiated with the SDE, which requires the solution to an optimization subproblem. While the problem considered by \cite{lelievre2024optimizing} is more general, our diffusion matrix \cref{eq:general_Q} is much simpler and does not require solving an optimization problem. Another closely related work is \cite{lelievre2013optimal}, which considered preconditioning of the form $\bQ=\bI + \mJ$. Here $\bI$ is the identity matrix and $\mJ$ is skew-symmetric, i.e. $\mJ = -\mJ^T$. \cite{lelievre2013optimal} studied the optimal $\mJ$ when the potential $f$ is a quadratic function, which is also the focus of this work. 
\end{remark}

\begin{remark}\label{rmk:comparison}
In \cite{li2022hessian}, the authors also studied \cref{eq:gaul} with $\mC = \bI$ which they called Hessian-Free High-Resolution (HFHR) dynamics. They considered potential functions $f$ that are $L$-smooth and $m$-strongly convex. They proved a convergence rate of $\frac{\sqrt{m}}{2\sqrt{\kappa}}$ in continuous time in terms of Wasserstein-2 distance between the target and sample measure. \cite{li2022hessian} used the randomized midpoint method \cite{shen2019randomized} combined with Strang splitting as their discretization and showed an interation complexity of $\widetilde{\mathcal{O}}(\sqrt{d}/\varepsilon)$. Specifically, \cite{li2022hessian} showed that for a two-dimensional Gaussian target measure, under the optimal choice of parameter (damping parameter $\gamma$ and step size $h$) for underdamped Langevin dynamics with Euler-Maruyama discretization, the convergence rate is $\mathcal{O}\big((1-\kappa^{-1})^{k}\big)$. This rate is recovered in \Cref{cor:ul_mixing_time}. On the other hand, \cite{li2022hessian} showed that under their choice of parameter for HRHF, the convergence rate is $\mathcal{O}\big((1-2\kappa^{-1})^{k}\big)$, which is a slight improvement compared with underdamped Langevin dynamics. In this work, we performed a detailed eigenvalue analysis of GAUL on Gaussian target measure. We showed that under our choice of parameters $(\gamma,a,h)$, the convergence rate towards the biased target measure is $\mathcal{O}\big((1-c\sqrt{\kappa})^{k}\big)$ for some constant $c$.
\end{remark}

\section{Analysis of GAUL on quadratic potential functions}\label{sec:analysis_GAUL}
In this section, we establish the convergence rate for the proposed SDE \cref{eq:FPK_general} towards the target distribution following a Gaussian distribution. 
\subsection{Problem set-up}\label{sec:setup}
In this subsection, we present the main problem addressed in this paper. We are interested in sampling from a distribution whose probability density function is proportional to $\exp(-f(\vx))$ for $f:\mathbb R^d \to \mathbb R$. In this paper, we focus on a concrete example in which the potential function $f$ is quadratic, and thus the target distribution is a Gaussian distribution. Let 
\begin{equation}\label{eq:quadratic_f}
    f(\vx)=\frac{1}{2}\vx^{T}\Sigma_*^{-1} \vx, 
\end{equation}
where $\vx\in\mathbb{R}^d$ and $\Sigma_* \succ 0$ is a symmetric positive definite matrix in $\mathbb{R}^{d\times d}$. Define 
\begin{equation}\label{eq:widetilde_Sigma}
    \widetilde \Sigma = \begin{pmatrix}
        \Sigma_* & 0 \\
        0 & \bI 
    \end{pmatrix}\,.
\end{equation}
As in the previous section, denote by $\vX=(\vx,\vp)\in \mathbb R^{2d}$. And $H(\vX) = f(\vx) + \|\vp\|^2/2$. Then, we can write 
\begin{equation}\label{eq:quadratic H} 
H(\vX) = \frac{1}{2} \vX^T \begin{pmatrix}
    \Sigma_*^{-1} & 0 \\
    0 & \bI 
\end{pmatrix} \bX := \frac{1}{2} \bX^T \widetilde\Sigma^{-1} \bX\,.
\end{equation}
Define the target density $\pi:\mathbb R^{2d} \to \mathbb R$ to be 
\begin{equation}\label{eq:target_measure}
    \Pi(\bX) = \frac{1}{Z} \exp(-H(\bX))\,,
\end{equation}
where $H(\bX)$ is given by \cref{eq:quadratic H} and $Z=\int_{\mathbb{R}^{2d}} \exp(-H(\bX))d\bX$ is a normalization constant such that $\Pi(\bX)$ integrates to one on $\mathbb R^{2d}$. We also define the $\vx$-marginal target density to be 
\begin{equation}\label{eq:target_measure_x}
    \pi(\vx) = \frac{1}{Z_1}\exp(-f(\vx))\,,
\end{equation}
where $f(\vx)$ is given by \cref{eq:quadratic_f} and $Z_1=\int_{\mathbb{R}^d}\exp(-f(\vx)) d\vx$ is a normalization constant. 
\begin{remark}\label{rmk:diagonalization}
    Note that for any symmetric positive definite $\Sigma_*$, we have that $\Sigma_*^{-1} = \rmP \Lambda \rmP^T$ for some orthogonal matrix $\rmP$ and diagonal matrix $\Lambda = \mathrm{diag}(s_1, \ldots,s_d)$ with $s_1 \geq \cdots s_d > 0$. 
    By a change of variable $\vy = \rmP^T \vx$, one can rewrite $f(\vx)$ in terms of $\vy$, such that 
    $$ f(\vx)=\frac{1}{2}\vx^T\Sigma_*^{-1}\vx=\frac{1}{2}\vx^T \rmP \Lambda \rmP^T\vx=\frac{1}{2}\vy^T\Lambda \vy.$$
       For simplicity of notation, we assume that $\rmP=\bI$ and $\Sigma_*^{-1} = \Lambda$ is a diagonal matrix. We denote by $\kappa=s_1/s_d$ the condition number of $f$. We will also assume that $s_1>1>s_d$ throughout this paper. Furthermore, to simplify our analysis, we consider $\mC = \mathrm{diag}(c_1,\ldots,c_d)$.  
\end{remark}
\subsection{Continuous time analysis}\label{sec:continuous}
In this subsection, we study the convergence of GAUL. In particular, we analyze the convergence of the Fokker-Planck equation \cref{eq:FPK_general} to the target density \cref{eq:target_measure}, \cref{eq:target_measure_x} by directly studying an ODE system of the covariance of the distribution. 
\begin{proposition}\label{prop:Sigma_ODE}
Let $\bX_t$ be the solution of \cref{eq:SDE} where $H(\bX)$ is given by \cref{eq:quadratic H}, and $\bX_0 \sim \mathcal{N}(0,\bI_{2d\times 2d})$. Then $\bX_t \sim \mathcal{N}(0,\Sigma(t))$ where the covariance $\Sigma(t)$ satisfies the following matrix ODE: 
\begin{equation}\label{eq:Sigma}
    \dot \Sigma(t) =2\,\mathrm{sym}(\bQ(\bI-\widetilde\Sigma^{-1}\Sigma(t) ))\,. 
\end{equation}
{Moreover, equation \cref{eq:Sigma} is well-defined, and has a  solution for all $t\geq 0$, such that $\Sigma(t)$ is symmetric semi-positive definite.} 
\end{proposition}
The proof is postponed in \Cref{sec:post_proof}. We denote by $\Sigma_{ij}(t) \in \mathbb R^{d\times d}$ the block components of $\Sigma(t)\in \mathbb R^{2d\times 2d}$: 
$$
\Sigma(t) = \begin{pmatrix}
    \Sigma_{11}(t)&\Sigma_{12}(t) \\
    \Sigma_{12}^T(t) & \Sigma_{22}(t)
\end{pmatrix}\,.
$$
Then we can write \cref{eq:Sigma} in terms of the block components. 
\begin{corollary}\label{cor:Sigma_ij}
The componentwise covariance matrix $\Sigma_{ij}(t)$ satisfies the following ODE system 
\begin{subequations}\label{eqs:Sigma_ij}
\begin{align}
    \dot \Sigma_{11} &= -2a(\mathrm{sym}(\mC \Sigma_*^{-1}\Sigma_{11})-\mC) + 2\,\mathrm{sym}(\mC \Sigma_{12})\,, \\
    \dot \Sigma_{22} &= -2 \,\mathrm{sym}(\Sigma_*^{-1}\Sigma_{12}) - 2\gamma(\Sigma_{22}-\bI)\,, \label{eq:Sigma_ij_b}\\
    \dot\Sigma_{12}&=-a \mC \Sigma_*^{-1}\Sigma_{12} - (\mC  - \mC \Sigma_{22}) + (\bI - \Sigma_{11}\Sigma_*^{-1}) - \gamma\Sigma_{12}\,,
\end{align}
Moreover, with initial conditions $\Sigma_{11}(0) = \Sigma_{22}(0) = \bI$ and $\Sigma_{12}(0) = 0$, the stationary states of $\Sigma_{11}(t)$, $\Sigma_{22}(t)$ and $\Sigma_{12}(t)$ are given by $\Sigma_*$, $\bI$ and 0 respectively.
\end{subequations}
\end{corollary}

From now on, we consider $\mC=\bI$ in our analysis. We address our results for $\mC \neq \bI$ in \Cref{rmk:cont_mix_time_C} and \Cref{rmk:dis_mix_C}. Note that when $\mC =\bI$, we have $\bQ=\mathrm{sym}(\bQ)$ is always positive semidefinite for $a,\gamma \geq 0$. 
Our next theorem makes sure that the stationary state of equation \cref{eq:Sigma} is actually unique and characterizes the convergence speed of the convariance matrix towards its stationary state.
\begin{theorem} \label{thm:gaussian_convergence}
Let $\bX_t$ be the solution of \cref{eq:SDE} where $H(\bX)$ is given by \cref{eq:quadratic H}, and $\bX_0 \sim \mathcal{N}(0,\bI_{2d\times 2d})$. Then $\Sigma(t)$ converges to the unique stationary state $\widetilde \Sigma$ given in \cref{eq:widetilde_Sigma}. The optimal choice of $\gamma$ is given by $\gamma^*= as_d + 2\sqrt{s_d} $ under which we have $\| \Sigma_{11}(t)-\Sigma_*\|_\mathrm{F} = \mathcal{O}(te^{-( 2as_d+ 2\sqrt{s_d})t})$ and $\| \Sigma_{22}(t)-\bI\|_\mathrm{F} = \mathcal{O}(te^{-( 2as_d+ 2\sqrt{s_d})t})$ for $t\geq 1$. 
\end{theorem}

\begin{proof}
As mentioned in \Cref{rmk:diagonalization}, we  consider $\Sigma_*^{-1} = \Lambda $. By our assumption on $\bX_0 $, \cref{eqs:Sigma_ij} implies that $\Sigma_{11}(t)$, $\Sigma_{22}(t)$ and $\Sigma_{12}(t)$ will be diagonal matrices for all $t >0$. This simplifies the ODE system \cref{eqs:Sigma_ij}. After some manipulation, we obtain 
\begin{equation}\label{eq:Sigma_ODE}
     \begin{pmatrix}
        \dot \Sigma_{11}\\
        \dot \Sigma_{22}\\
        \Ddot \Sigma_{22}
    \end{pmatrix} = \underbrace{\begin{pmatrix}
        -2a\mC \Sigma_*^{-1} & -2\gamma \mC\Sigma_* & -\mC\Sigma_* \\
        0&0& \bI \\
        2 \Sigma_*^{-2} & 2(-1-a\gamma)\mC\Sigma_*^{-1}-2\gamma^2 \bI & -3\gamma\bI -a\mC\Sigma_*^{-1}
    \end{pmatrix}}_{\mathcal{D}} \begin{pmatrix}
        \Sigma_{11}\\
        \Sigma_{22}\\
        \dot \Sigma_{22}
    \end{pmatrix} + \mathbf{T}\,,
\end{equation}
where 
$$
\mathbf{T} = \begin{pmatrix}
    2a\mC +2\gamma\mC\Sigma_* \\
    0 \\
    2a\gamma\Sigma_*^{-1}\mC+2\gamma^2 \bI  + 2 \Sigma_*^{-1}\mC-2\Sigma_*^{-1}
\end{pmatrix}\,,
$$
And $\mC = \bI$. We have already seen in \Cref{cor:Sigma_ij} that the stationary state of $\Sigma(t)$ is $\widetilde \Sigma$ given in \cref{eq:widetilde_Sigma}. To show uniqueness, we compute the eigenvalues of $\mathcal{D}$: 
\begin{align*}
    \lambda^{(i)}_0 &= -a s_i - \gamma \,,\nonumber \\
    \lambda^{(i)}_1 &= -as_i - \gamma  - \sqrt{\gamma^2-2a\gamma s_i + s_i(-4+a^2 s_i)}\,, \\
    \lambda^{(i)}_2 &= -as_i - \gamma  + \sqrt{\gamma^2-2a\gamma s_i + s_i(-4+a^2 s_i)}\,,
\end{align*}
where $s_i$'s are the diagonal elements of $\Lambda$ for $i=1,\ldots,d$. It is clear that 0 is not an eigenvalue of $\mathcal{D}$. Therefore, $\widetilde \Sigma$ is the unique stationary state for $\Sigma(t)$. The convergence speed of \cref{eq:Sigma_ODE} is essentially controlled by the largest real part of the eigenvalues of $\mathcal{ D}$. Note that for all $i$, $$\Re(\lambda^{(i)}_2)\geq\Re(\lambda^{(i)}_0)\geq\Re(\lambda^{(i)}_1)\,, $$
where $\Re(z)$ denotes the real part of $z\in \mathbb C$. 
Therefore, to characterize the convergence speed of \cref{eq:Sigma_ODE}, it suffices to control $\max_i \Re(\lambda^{(i)}_2)$. By \Cref{lemma:matrix_gamma}, we know that for any given $a\geq 0$, the optimal choice of $\gamma$ is 
$$
\gamma^* = \argmin_{\gamma>0} \max_i \Re(\lambda^{(i)}_2) = a s_d  + 2\sqrt{s_d}\,. 
$$
With $\gamma = \gamma^*$, we get that 
$$
\max_{i,j} \Re(\lambda^{(i)}_j) \leq \max_{i} \Re(\lambda^{(i)}_2) \leq  -2as_d - 2\sqrt{s_d}\,. 
$$
This leads to 
\begin{equation} \label{eq:frobenius_bound }
\left\|  \begin{pmatrix}
        \Sigma_{11}(t) - \Sigma_*\\
        \Sigma_{22}(t)- \bI \\
        \dot \Sigma_{22}(t)
    \end{pmatrix}  \right\|_\mathrm{F} \leq  C_1 te^{-( 2as_d + 2\sqrt{s_d})t}\,,
\end{equation}
which is valid for $t\geq 1$. The constant $C_1$ depends on $d$, $s_1$, $s_d^{-1}$ at most polynomially according to \Cref{lemma:const_C1_bound}. Note that the extra $t$ dependence comes from the repeated eigenvalue $\lambda^{(d)}_0 =\lambda^{(d)}_1=\lambda^{(d)}_2 $ when $\gamma=\gamma^*$. By a triangle inequality, we get 
$$
\| \Sigma_{11}-\Sigma_*\|_\mathrm{F} \leq \left\|  \begin{pmatrix}
        \Sigma_{11}(t) - \Sigma_*\\
        \Sigma_{22}(t)- \bI \\
        \dot \Sigma_{22}(t)
    \end{pmatrix}  \right\|_\mathrm{F} \leq  C_1 te^{-( 2as_d + 2\sqrt{s_d})t}\,. 
$$
And similarly,  
$$
\| \Sigma_{22}-\bI \|_\mathrm{F} \leq C_1 te^{-( 2as_d + 2\sqrt{s_d})t}\,. 
$$
\end{proof}
\begin{remark}
The choice $a=0$ corresponds to underdamped Langevin dynamics (UL). Taking $a>0$ gives an extra factor of $e^{-2as_d t}$ in terms of convergence. 
\end{remark}
\begin{definition}[Mixing time]
    The total variation between two probability measures $\mathcal{P}$ and $\mathcal Q$ over a measurable space $(\mathbb{R}^d, \mathcal{F})$ is 
    $$
    \mathrm{TV}(\mathcal{P},\mathcal{Q}) = \sup_{A\in \mathcal{F}} | \mathcal{P}(A) - \mathcal{Q}(A) | \,. 
    $$
    Let $\mathcal{T}_p$ be an operator on the space of probability distributions. Assume that $\mathcal{T}_p^k (\nu_0) \to \nu$ as $k\to \infty$ for some initial distribution $\nu_0$ and stationary distribution $\nu$. The discrete $\delta$-mixing time $(\delta \in (0,1))$ is given by 
    $$
    t_{\mathrm{mix}}^{\mathrm{dis}}(\delta ; \nu_0,\nu) = \min \{ k \,|\, \mathrm{TV} (\mathcal{T}_p^k (\nu_0),\nu ) \leq \delta \} \,. 
    $$
    Similarly, if $\mathcal{T}_p(t;\cdot)$ is an operator for each $t\geq 0$ with $\mathcal{T}_p(0;\cdot) = \mathrm{id}(\cdot)$ and assume that $\mathcal{T}_p(t;\nu_0) \to \nu$ as $t\to \infty$. The continuous $\delta$-mixing time $(\delta \in (0,1))$ is given by
    $$
    t_{\mathrm{mix}}^{\mathrm{cont}}(\delta ; \nu_0, \nu) = \min \{ t \,|\, \mathrm{TV} (\mathcal{T}_p (t;\nu_0),\nu ) \leq \delta \} \,. 
    $$
\end{definition}
\begin{theorem}[\cite{devroye2018total}]
    Let $\mu\in \mathbb R^d$, $\Sigma_1$, $\Sigma_2$ be two positive definite covariance matrices, and $\lambda_1, \ldots,\lambda_d$ denote the eigenvalues of $\Sigma_1^{-1}\Sigma_2 - \bI$. Then the total variation satisfies 
    $$
    \mathrm{TV}(\mathcal{N}(\mu,\Sigma_1), \mathcal{N}(\mu,\Sigma_2)) \leq 
    \frac{3}{2} \min \left \{1, \sqrt{\sum_{i=1}^d \lambda_i^2} \right\} \,. 
    $$
\end{theorem}
A straightforward corollary follows from Schur decomposition theorem. 
\begin{corollary} \label{cor:TV_frobenius}
We have 
 $$
    \mathrm{TV}(\mathcal{N}(\mu,\Sigma_1), \mathcal{N}(\mu,\Sigma_2)) \leq 
    \frac{3}{2} \min \left \{1, \|\Sigma_1^{-1}\Sigma_2 - \bI\|_\mathrm{F} \right\} \,.
$$
\end{corollary}
Using \Cref{thm:gaussian_convergence} and \Cref{cor:TV_frobenius}, we obtain the following mixing time theorem when the potential function $f$ is quadratic. 
\begin{theorem}[Continuous mixing time]\label{thm:cont_mix_time}
Consider the same setting as in \Cref{thm:gaussian_convergence}. Consider $0<  \delta \ll 1$. Then 
    $$t_{\mathrm{mix}}^{\mathrm{cont}}(\delta ; \nu_0,\pi) \leq \frac{\mathcal{O}(\log(d) + \log(\kappa))+ \log(1/\delta)}{as_d + 2\sqrt{s_d}}\,.
    $$
  Here $\nu_0$ is the distribution of $\vx$, which is $\mathcal{N}(0,\bI_{d\times d})$. $\pi$ is the target density in the $\vx$ variable given in \cref{eq:target_measure_x}.
\end{theorem}

\begin{proof}
We shall use \Cref{cor:TV_frobenius} with
$$
\Sigma_1 = \Sigma_*\,, \qquad \Sigma_2 = \Sigma_{11}(t) \,. 
$$
We have 
\begin{align*}
    \|\Sigma_1^{-1}\Sigma_2 - \bI\|_\mathrm{F} & = \left\| 
        \Sigma_*^{-1}(\Sigma_{11}(t)-\Sigma_*)
     \right\|_\mathrm{F}  \\
    &\leq C_1 t e^{-( 2as_d + 2\sqrt{s_d})t}s_1 \,.
\end{align*}
By a direct computation, we get 
\begin{align*}
    t_{\mathrm{mix}}^{\mathrm{cont}}(\delta ; \nu_0,\pi) \leq \frac{\log(\tilde C_1 / \delta) }{as_d + 2\sqrt{s_d}}\,,
\end{align*}
where $\tilde C_1 = \frac{3}{2} C_1 s_1$. By \Cref{lemma:const_C1_bound}, we have that 
$$
 t_{\mathrm{mix}}^{\mathrm{cont}}(\delta ; \nu_0,\pi) \leq \frac{\mathcal{O}(\log(d \kappa))+ \log(1/\delta)}{as_d + 2\sqrt{s_d}}\,.
$$
\end{proof}

\begin{remark}\label{rmk:cont_mix_time_C}
When $\mC = \mathrm{diag}(c_1,\ldots,c_d)$ and $\mathrm{sym}(\bQ) \succeq 0$ in \cref{eq:general_Q}, our proof can be easily adapted to show similar results in \Cref{thm:cont_mix_time}: $$t_{\mathrm{mix}}^{\mathrm{cont}}(\delta ; \nu_0,\pi) \leq \frac{\mathcal{O}(\log(d) + \log(\hat{\kappa}))+ \log(1/\delta)}{a\hat{s}_d + 2\sqrt{\hat{s}_d}}\,,
    $$
where $\hat{s}_i$ is the $i$-th largest eigenvalue of matrix $\mC \Sigma_*^{-1}$. And $\hat{\kappa}=\hat{s}_1/\hat{s}_d$.  In other words, the matrix $\mC$ can be viewed as a preconditioner for the target covariance matrix in the sampling problem.  
\end{remark}

\subsection{Discrete time analysis}\label{sec:discrete}
To implement \cref{eq:SDE}, we need to consider its time discretization. As discretization is not the focus of this paper, we will only analyze the simplest discretization using the Euler-Maruyama method in \Cref{sec:euler_maruyama}. 

Let us first make a few observations regarding the discretization in \Cref{sec:euler_maruyama}.  After a straightforward computation, we obtain the following update rule. 
\begin{proposition}\label{prop:discrete_update_cheng}
The Euler-Maruyama discretization of \cref{eq:SDE} given in \Cref{sec:euler_maruyama} with step size $h$ can be written in the following form
        \begin{align}\label{eq:discrete_update}
            \begin{pmatrix}
                \vx_{n+1} \\
                \vp_{n+1} 
            \end{pmatrix} = \mA\begin{pmatrix}
                \vx_{n} \\
                \vp_{n} 
            \end{pmatrix} + \mL \vz \,, 
        \end{align}
where \begin{equation} \label{eq:A}
\mA = \bI_{2d \times 2d}  - \underbrace{h\begin{pmatrix}
    a  \Lambda & - \bI_{d\times d} \\
    \Lambda & \gamma  \bI_{d\times d} 
\end{pmatrix}}_{\mG }\,, \qquad \mL = \begin{pmatrix}
    \sqrt{2ah}\bI & 0 \\
    0& \sqrt{2\gamma h} \bI 
\end{pmatrix}\,. 
\end{equation}
And $\vz$ is a $2d$-dimensional Brownian motion, i.e., $\vz  \sim \mathcal{N}(0,\bI_{2d\times 2d})$. 
\end{proposition}

Using $\cref{eq:discrete_update}$, we can derive the evolution of the mean and covariance at each time step. As before, let us denote by $\vX_n = (\vx_n, \vp_n)$. 
\begin{corollary}
Suppose that $\mathbb E(\vx_0) = \mathbb E(\vp_0) = 0$. Then $$\mathrm{cov}(\vX_{n+1},\vX_{n+1}) = \mA \mathrm{cov}(\vX_{n},\vX_{n})\mA^T + \mL \mL^T\,. $$
\end{corollary}
\begin{proof}
From \cref{eq:discrete_update}, it is clear that $\mathbb E(\vx_n) = \mathbb E(\vp_n) = 0$ for all $n\geq 0$. We calculate 
\begin{align*}
    \mathrm{cov}(\vX_{n+1},\vX_{n+1}) &= \mathbb E \big( \mA \vX_n \vX_n^T \mA^T  + \mA \vX_n \vz^T \mL^T + \mL \vz \vX_n^T \mA^T + \mL \vz \vz^T \mL^T \big) \\
    &= \mA \mathrm{cov}(\vX_{n},\vX_{n})\mA^T + \mL \mL^T\,.
\end{align*}
\end{proof}

\begin{corollary}\label{cor:Yn}
    Denote by $\bY^*$ a solution to the fixed point equation $\mY = \mA \mY \mA^T + \mL \mL^T$. And let $\mY_n = \mathrm{cov}(\vX_n,\vX_n) - \mY^*$. Then 
    $$
    \mY_{n+1} = \mA \mY_n \mA^T\,. 
    $$
\end{corollary}

\begin{theorem}\label{thm:discrete_time_convergence}
    Suppose $a\geq \frac{2}{\sqrt{s_{1}}-\sqrt{s_d}}$ and the step size $h$ satisfies $0<h<1/(as_1 + \gamma)$ and $\gamma = \gamma^* = as_d + 2\sqrt{s_d}$.  Then there exists a unique $\mY^*$ satisfying 
    $$
    \mY^* = \mA \mY^*\mA^T + \mL \mL^T\,. 
    $$
    Moreover, the iteration $\mY_{k+1} = \mA \mY_k \mA^T + \mL \mL^T$ converges to $\mY^*$ linearly: $ \| \mY_{k} - \mY^* \|_\mathrm{F} \leq \widetilde C h^2k^2 (1 - \frac{h}{2}(as_d + \sqrt{s_d}))^{2k-2} $, where the constant $\widetilde C =  d^2 \cdot \mathcal{O}(\mathrm{poly}(\kappa))$. 
\end{theorem}
\begin{proof}
    Existence: we directly compute this stationary point in \Cref{lemma:fixed_point}. Uniqueness: by \Cref{lemma:discrete_eig_bound} and \Cref{cor:uniqueness} we see that $\mY^*$ is unique. The convergence rate is proved in \Cref{lemma:discrete_eig_bound} and \Cref{thm:conv_rate}. 

\end{proof}
\begin{theorem}[Discrete mixing time]\label{thm:dis_mix_time}
    Suppose $\sqrt{s_1}-\sqrt{s_d}\geq 2$. We take $a=1$, $\gamma=\gamma^*= s_d + 2\sqrt{s_d}$, $h = 1/5s_1$. If we use the Euler-Maruyama scheme for \cref{eq:SDE}, then for $0<\delta \ll 1$,
    \begin{equation}\label{eq:dis_mix_time_gaul}
        t_{\mathrm{mix}}^{\mathrm{dis}}(\delta ; \nu_0, \tilde\pi) = \mathcal{O}\left( \frac{\log(\kappa)+ \log(1/\delta) + \log(d)}{\frac{1}{\kappa}+\frac{1}{\sqrt{\kappa s_1}}} \right)\,.
    \end{equation}
    Here $\nu_0$ is the distribution of $\vx$, which is $\mathcal{N}(0,\bI_{d\times d})$. $\tilde\pi$ is the target density in the $\vx$ variable which is a zero mean Gaussian distribution with a variance given by \cref{eq:x_element}. 
\end{theorem}

\begin{proof}
    Note that from our previous notation, we have that 
    $$
    \mathrm{cov}(\vx_k,\vx_k)  = \begin{pmatrix}
        \bI_{d\times d} & 0 
    \end{pmatrix}
    \mathrm{cov}(\vX_k,
    \vX_k)  \begin{pmatrix}
        \bI_{d\times d} \\ 0 
    \end{pmatrix} =: \widetilde \mY_k \,. 
    $$
    Moreover, let us define 
    $$
    \widetilde \mY^* = \begin{pmatrix}
        \bI_{d\times d} & 0 
    \end{pmatrix}
    \mY^*  \begin{pmatrix}
        \bI_{d\times d} \\ 0 
    \end{pmatrix}
    $$
    to be the limiting covariance in the $\vx$ variable for the discretization ($\mY^*$ is defined in \Cref{thm:discrete_time_convergence}). Clearly, we have that 
    \begin{equation}\label{eq:Y_tilde_norm}
    \|\widetilde \mY_k - \widetilde \mY^* \|_\mathrm{F} \leq \| \mY_k - \mY^* \|_\mathrm{F} \leq \widetilde C h^2k^2 (1 - \frac{h}{2}(as_d + \sqrt{s_d}))^{2k-2} \,. 
     \end{equation}
    Using \Cref{cor:TV_frobenius}, we compute 
    \begin{align*}
        \|(\widetilde \mY^*)^{-1} \widetilde \mY_k - \bI \|_\mathrm{F} &= \|(\widetilde \mY^*)^{-1} ( \widetilde \mY_k - \widetilde \mY^*) \|_\mathrm{F} \nonumber \\
        &\leq \| (\widetilde \mY^*)^{-1} \|_\mathrm{F} \|\widetilde \mY_k - \widetilde \mY^* \|_\mathrm{F}\,.
    \end{align*}
    By \Cref{lemma:fixed_point}, $\widetilde \mY^*$ is a diagonal matrix. Therefore $(\widetilde \mY^*)^{-1}$ is also a diagonal matrix. Moreover, from \cref{eq:x_element}, we see that $\|(\widetilde \mY^*)^{-1}\|_\mathrm{F} \leq  \sqrt{d}\mathcal{O}(\mathrm{poly}(\kappa)) $. Therefore, we obtain 
    \begin{align*}
        \|(\widetilde \mY^*)^{-1} \widetilde \mY_k - \bI \|_\mathrm{F}  &\leq d^{5/2} \cdot \mathcal{O}(\mathrm{poly}(\kappa)) h^2k^2 (1 - \frac{h}{2}(s_d + \sqrt{s_d}))^{2k-2} \nonumber \\
        &\leq d^{5/2} \cdot \mathcal{O}(\mathrm{poly}(\kappa)) h^2k^2 e^{-(k-1)h (s_d + \sqrt{s_d}) }\,,
    \end{align*}
    where we used $1-x \leq e^{-x}$ for $x\in \mathbb R$ to get the second inequality. Letting $h = 1/5s_1$ and taking logarithm on both hand sides, we conclude that  
    $$
    t_{\mathrm{mix}}^{\mathrm{dis}}(\delta ; \nu_0, \tilde\pi) \leq \frac{\mathcal{O}(\log(d)) + \mathcal{O}(\log(\kappa)) + \log(1/\delta)}{\frac{1}{10}(\frac{1}{\kappa} + \frac{1}{\sqrt{\kappa s_1}})}\,. 
    $$
\end{proof}
\begin{theorem}[A better choice of $a$]\label{thm:choice_of_a} The denominator of the mixing time given in \Cref{thm:dis_mix_time} can be improved to $\kappa^{-1/2}$ by choosing $a= \frac{2}{\sqrt{s_1}-\sqrt{s_d}}$, $\gamma=as_d + 2\sqrt{s_d}$ and $h = \frac{1}{2(as_1 + \gamma)}$. To be more precise, we have 
\begin{equation}\label{eq:dis_mix_time_optimal_a}
    t_{\mathrm{mix}}^{\mathrm{dis}}(\delta ; \nu_0, \tilde\pi) = \mathcal{O}\left( \frac{\log(\kappa)+ \log(1/\delta) + \log(d)}{\frac{1}{\sqrt{\kappa}}} \right)\,.
\end{equation}
\end{theorem}
\begin{proof}
The proof will be very similar to that of \Cref{thm:dis_mix_time}. We start with \cref{eq:Y_tilde_norm}. And we can explicitly calculate 
\begin{align*}
    1 - \frac{h}{2}(as_d + \sqrt{s_d}) &= 1 - \frac{as_d + \sqrt{s_d}}{4(as_1 + as_d + 2\sqrt{s_d})} \\
    &= 1-\frac{2s_d + \sqrt{s_d}(\sqrt{s_1} - \sqrt{s_d})}{8(s_1 +s_d + \sqrt{s_d}(\sqrt{s_1}-\sqrt{s_d}))} \\
    &= 1- \frac{\sqrt{s_1 s_d} + s_d}{8(s_1 + \sqrt{s_1 s_d})}\\
    &\leq 1- \frac{1}{16 \sqrt{\kappa}}\,. 
\end{align*}
The rest of the proof is the same as the proof of \Cref{thm:dis_mix_time} and we will suppress it for brevity. 
\end{proof}

The following corollary follows from \Cref{lemma:h_size_ul} and the proof of \Cref{thm:dis_mix_time}. 
\begin{corollary}[Underdamped Langevin mixing time]\label{cor:ul_mixing_time}
Suppose $a=0$, $\gamma=2\sqrt{s_d}, h=\sqrt{s_d}/s_1$. If we use the Euler-Maruyama scheme for \cref{eq:SDE}, then for $0<\delta \ll 1$, \begin{equation}\label{eq:dis_mix_time_ul}
        t_{\mathrm{mix}}^{\mathrm{dis}}(\delta ; \nu_0, \tilde\pi) = \mathcal{O}\left( \frac{\log(\kappa)+ \log(1/\delta) + \log(d)}{\frac{1}{\kappa}} \right)\,,
    \end{equation}
    $\nu_0$ is the distribution of $\vx$, which is $\mathcal{N}(0,\bI_{d\times d})$. $\tilde\pi$ is the target density in the $\vx$ variable which is a zero mean Gaussian with variance given by \cref{eq:x_element} with $a=0$.
\end{corollary}
\begin{remark}
    $a=0$ in \cref{eq:SDE} corresponds to the underdamped Langevin dynamics. In this case, we show in \Cref{lemma:h_size_ul} that to guarantee convergence (to a biased target) the step size restriction on $h$ is more strict than when $a=1$. In particular, when $a=0$ it follows from \Cref{lemma:h_size_ul} that the choice $h=1/5s_1$ does not guarantee convergence if $s_d < 10^{-2}$. Comparing \cref{eq:dis_mix_time_optimal_a} and \cref{eq:dis_mix_time_ul}, we see that the mixing time for GAUL beats that of underdamped Langevin dynamics under the Euler-Maruyama discretization. We are aware that this does not imply the same result will hold when comparing the mixing time towards the true target distribution $\pi(\vx)$ given in \cref{eq:target_measure_x}, due to the presence of bias in the Euler-Maruyama scheme. Designing better discretization and reducing the bias in the stationary distribution is left as future works. 
\end{remark}

\begin{remark}\label{rmk:dis_mix_C}
When $\mC=\mathrm{diag}(c_1,\ldots,c_d)$ and $\mathrm{sym}(\bQ) \succeq 0$ in \cref{eq:general_Q}, we also have a similar mixing time described in \Cref{thm:choice_of_a}, which is \newline$\mathcal{O}\big( \sqrt{\hat \kappa}(\log(\hat \kappa)+ \log(1/\delta) + \log(d)) \big)$ when $a= \frac{2}{\sqrt{\hat s_1}-\sqrt{\hat s_d}}$, $\gamma=a \hat s_d + 2\sqrt{\hat s_d}$ and $h = \frac{1}{2(a\hat s_1 + \gamma)}$. The notation $\hat{s}_i$ and $\hat{\kappa}$ are defined in \Cref{rmk:cont_mix_time_C}. 
\end{remark}
\begin{remark}
    When the target potential $f$ is not a quadratic function, it is more technical in proving the convergence speed. A common technique to prove convergence in the Wasserstein-2 distance is by a coupling argument (see \cite{cheng2018underdamped,dalalyan2018sampling}). \cite{cao2023explicit} proved $L_2$ convergence under a Poincar\`e-type inequality using Bochner's formula. In the $L_1$ distance and KL divergence, \cite{feng2024fisher} design convergence analysis towards these problems. We leave the convergence analysis of general $f$ with optimal choices of preconditioned matrices $\bQ$ in future works. 
\end{remark}

\section{Numerical experiment}\label{sec:numerics}
In this section, we implement several numerical examples to compare the proposed SDE with the overdamped (labeled `ol') and underdamped (labeled `ul') Langevin dynamics. We use the same step size for all three algorithms. Recall that `ol' corresponds to the choice $a=1,\gamma=0$ and `ul' corresponds to $a=0$ in \cref{eq:SDE}. We set $\mC = \bI$. 
\subsection{Gaussian examples}
\subsubsection{One dimension}
\begin{figure*}[h!]
        \centering
        \begin{subfigure}[b]{0.24\linewidth}
            \centering
            \includegraphics[scale=0.20]{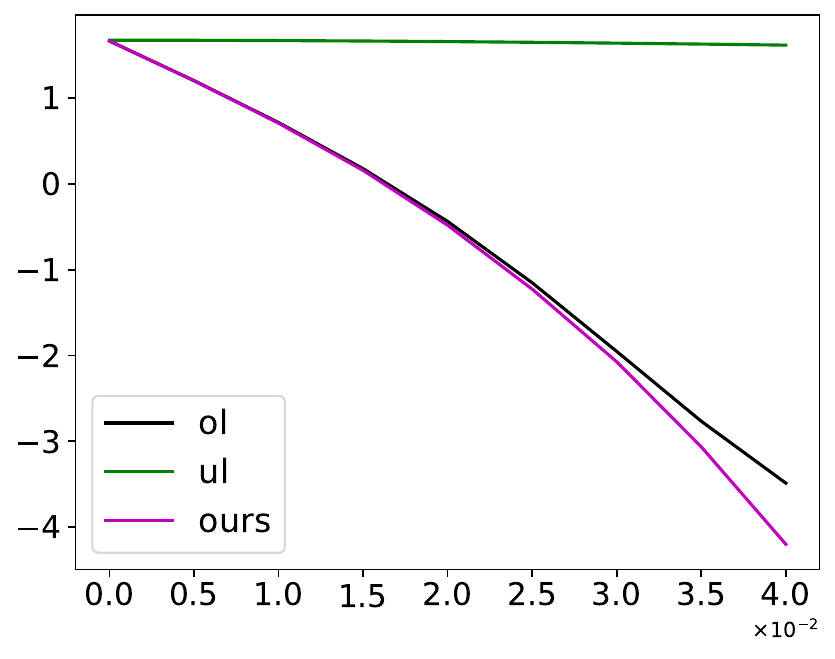}
            \caption{KL decay }\label{fig:1d_gauss_sig0.01}
        \end{subfigure}
        \hfill
        \begin{subfigure}[b]{0.24\linewidth}  
            \centering 
            \includegraphics[scale=0.21]{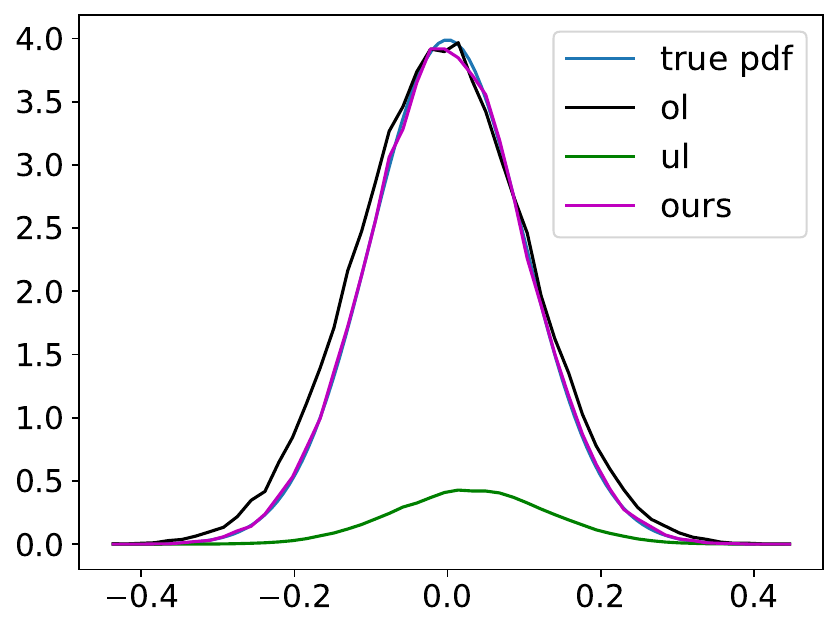}
             \caption{Density}\label{fig:1d_gauss_sig0.01_density}
        \end{subfigure}
        \hfill
        \begin{subfigure}[b]{0.24\linewidth}   
            \centering 
            \includegraphics[scale=0.21]{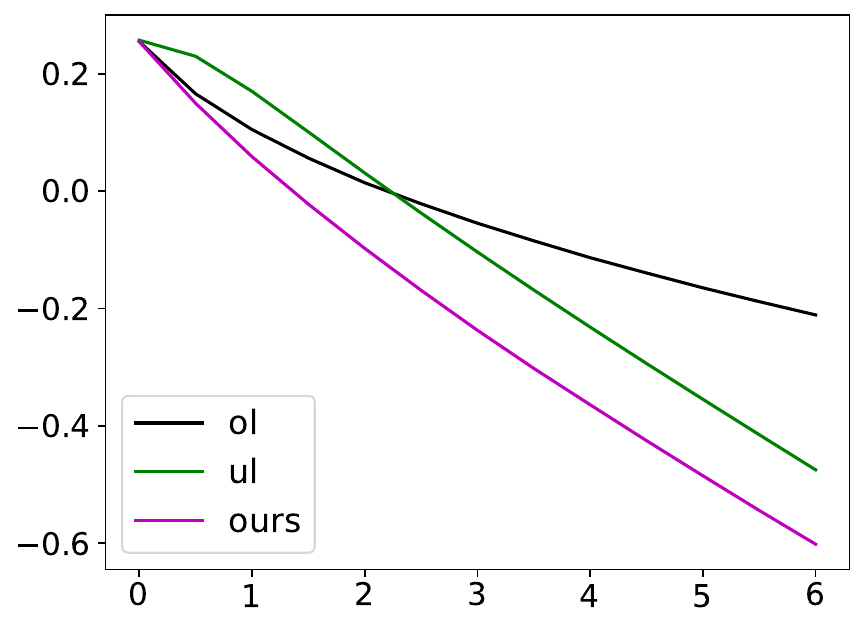}
            \caption{KL decay}\label{fig:1d_gauss_sig100}
        \end{subfigure}
        \hfill
        \begin{subfigure}[b]{0.24\linewidth}   
            \centering 
            \includegraphics[scale=0.21]{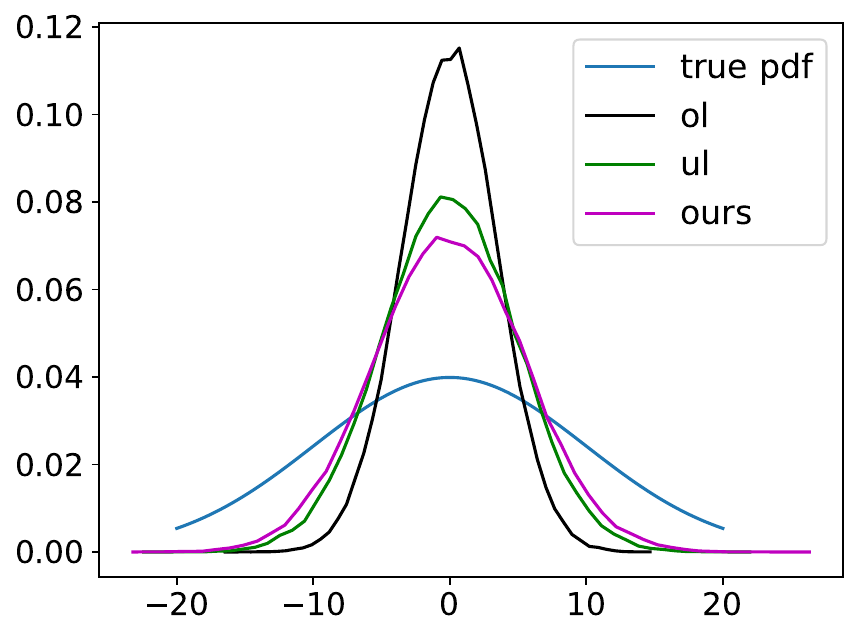}
            \caption{Density}\label{fig:1d_gauss_sig100_density}
        \end{subfigure}
         \caption{Convergence and density comparisons of three methods. (a) and (c): KL divergence between the sample and the target distribution, which is a one-dimensional Gaussian with zero mean and variance 0.01 (a), 100 (c). `ol' represents overdamped Langevin dynamics; `ul' represents underdamped Langevin dynamics. x-axis represents time and y-axis is in $\log_{10}$ scale. (b) and (d): density comparison at the end of the experiment between the three methods and the true density. }
        \label{fig:1d_gaussian}
    \end{figure*}

We begin with a simple example, a one dimensional Gaussian distribution with zero mean. In \Cref{fig:1d_gaussian}, we consider two cases where the variances are given by 0.01 and 100 respectively. We first sample $M=10^5$ particles from $\mathcal{N}(0,\bI_{2\times 2})$ (although our experiment is in one dimension, we need both $\vx$ and $\vp$ variables). When measuring the convergence speed, we use KL divergence in Gaussian distributions to measure the change of covariances. Note that we will only measure the KL divergence in the $x$ variable, since we are primarily interested in sampling distribution of the form $\frac{1}{Z}e^{-f(x)}$. In this experiment, we can make use of the fact that the sample distribution and the target distribution are both Gaussians. And the KL divergence between two centered Gaussians has a closed form expression: 
\begin{equation}\label{eq:KL for gaussian}  
\mathrm{D}_{\mathrm{KL}}(\Sigma(t),\widetilde \Sigma) = \frac{1}{2} \left(\mathrm{tr}(\Sigma(t)\widetilde\Sigma^{-1}) - \log \mathrm{det}(\Sigma(t) \widetilde\Sigma^{-1}) - d  \right)\,.
\end{equation}

In this one dimensional example, we study two cases where $\widetilde \Sigma = 0.01$ or $100$. $\Sigma(t)$ can be approximated by the unbiased sample variance. For $\widetilde \Sigma=0.01$, we choose time step size $h=10^{-4}$, total number of steps $N=400$, $\gamma_{ul}= 2\widetilde\Sigma^{-1/2}=20$, $\gamma_{pdd}=2\widetilde\Sigma^{-1/2}+\widetilde\Sigma^{-1}=120$. For $\widetilde\Sigma=100$, we choose the time step size $h=10^{-2}$, total number of steps $N=600$, $\gamma_{ul}= 2\widetilde\Sigma^{-1/2}=0.2$, $\gamma_{pdd}=2\widetilde\Sigma^{-1/2}+\widetilde\Sigma^{-1}=0.21$. In \Cref{fig:1d_gaussian}, we observe that our proposed method outperforms both overdamped and underdamped Langevin dynamics in both cases.   

\subsubsection{20 dimensions}
\begin{figure*}
        \centering
        \begin{subfigure}[b]{0.24\linewidth}
            \centering
            \includegraphics[scale=0.21]{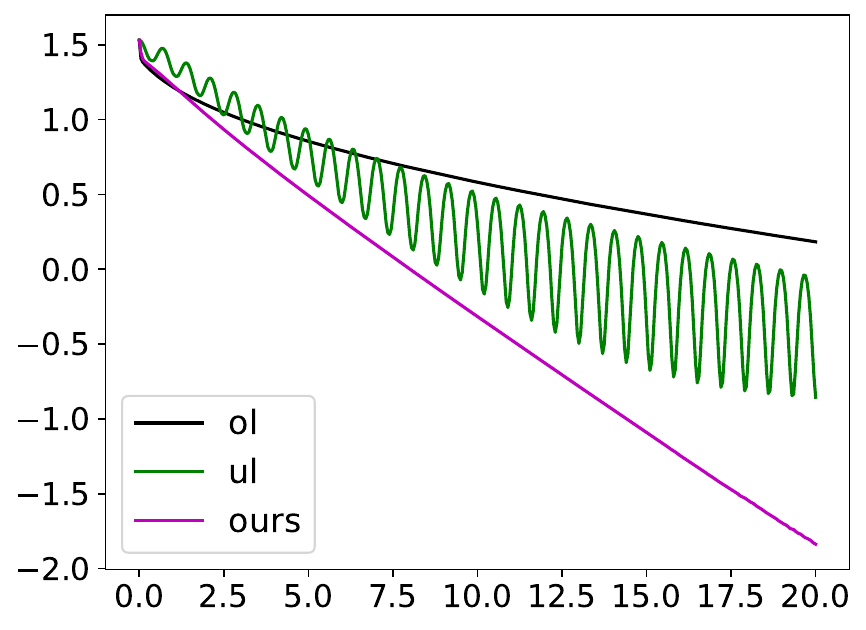}
            \caption{KL decay}\label{fig:2d_gauss_kl}
        \end{subfigure}
        \hfill
        \begin{subfigure}[b]{0.24\linewidth}  
            \centering 
            \includegraphics[scale=0.21]{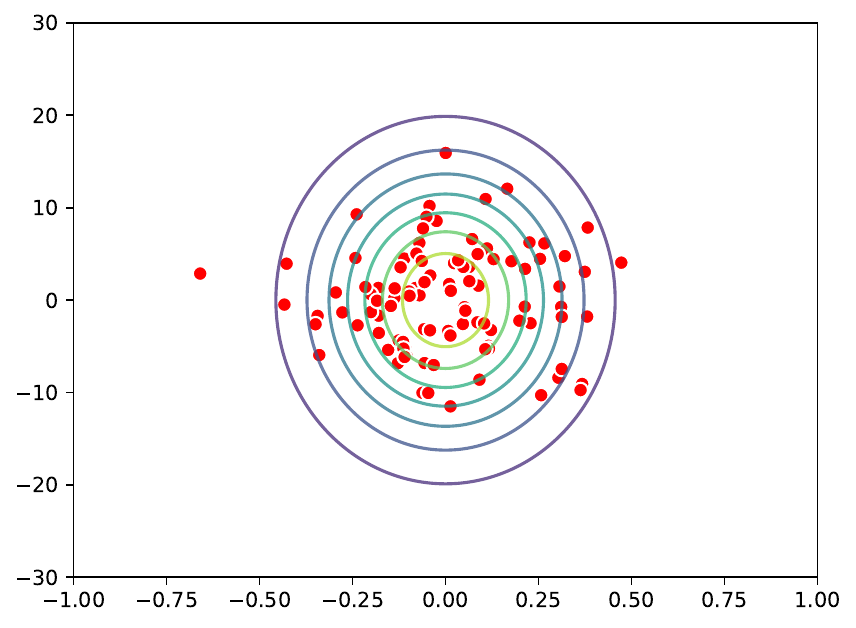}
             \caption{ol}\label{fig:2d_gauss_ol}
        \end{subfigure}
        \begin{subfigure}[b]{0.24\linewidth}   
            \centering 
            \includegraphics[scale=0.21]{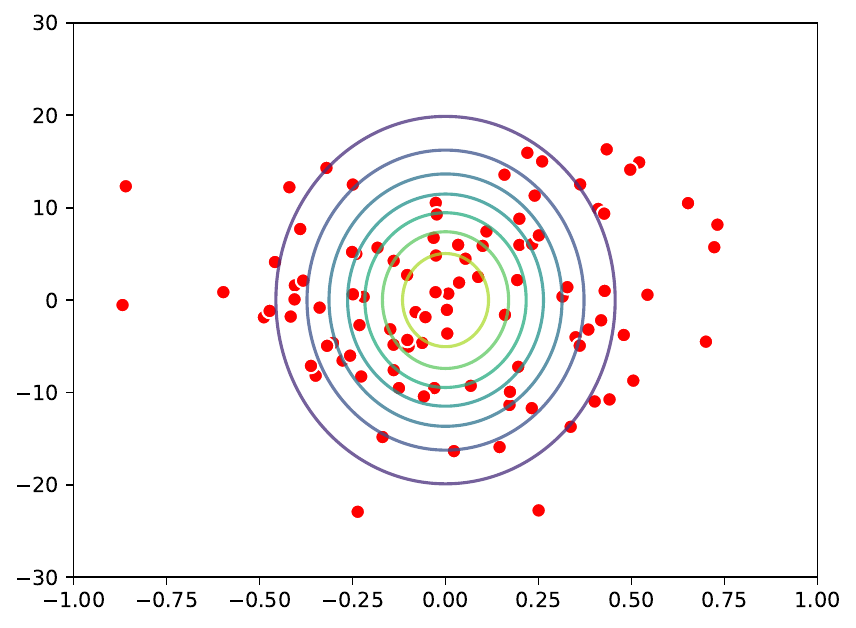}
            \caption{ul}\label{fig:2d_gauss_ul}
        \end{subfigure}
        \hfill
        \begin{subfigure}[b]{0.24\linewidth}   
            \centering 
            \includegraphics[scale=0.21]{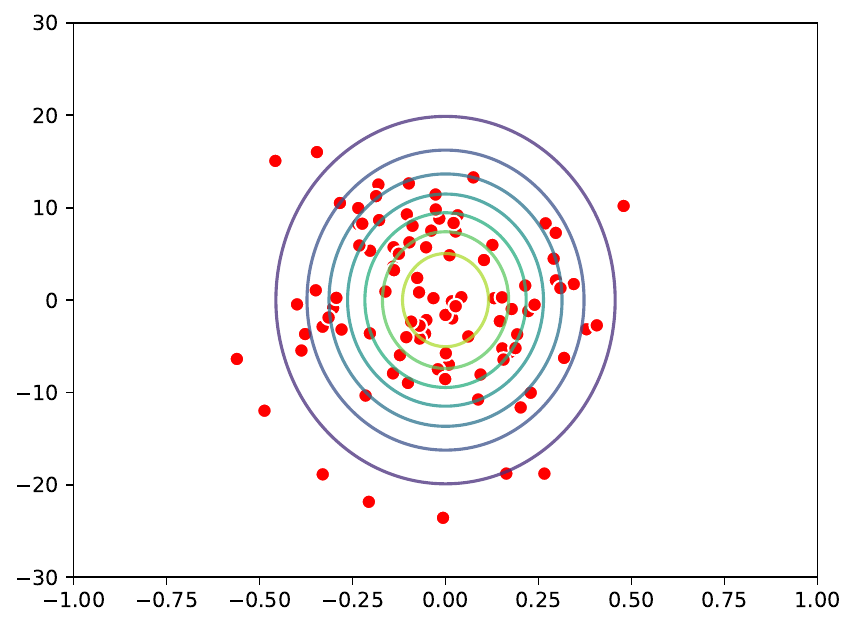}
            \caption{gaul}\label{fig:2d_gauss_pdd}
        \end{subfigure}
        \vskip\baselineskip
        \begin{subfigure}[b]{0.24\linewidth}
            \centering
            \includegraphics[scale=0.21]{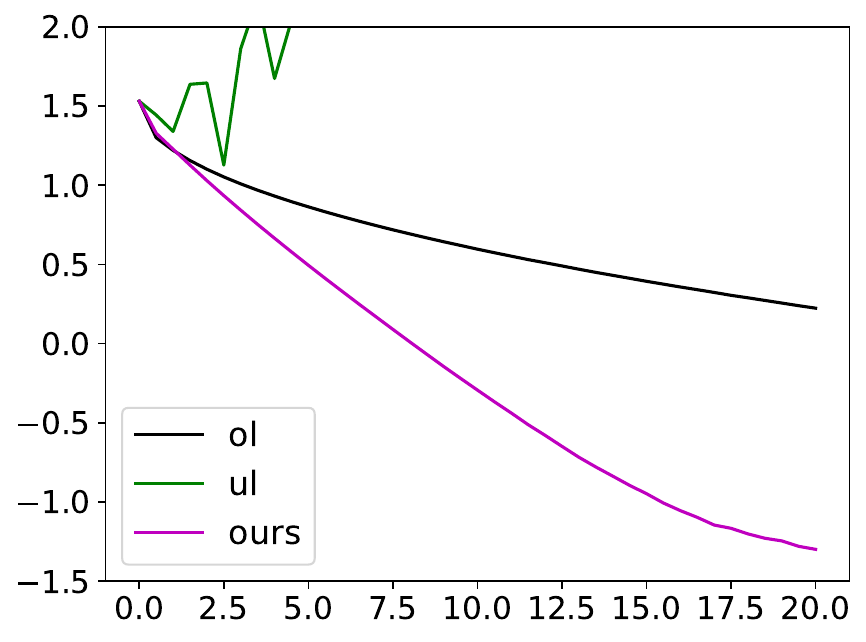}
            \caption{KL decay}\label{fig:2d_gauss_kl_large}
        \end{subfigure}
        \hfill
        \begin{subfigure}[b]{0.24\linewidth}  
            \centering 
            \includegraphics[scale=0.21]{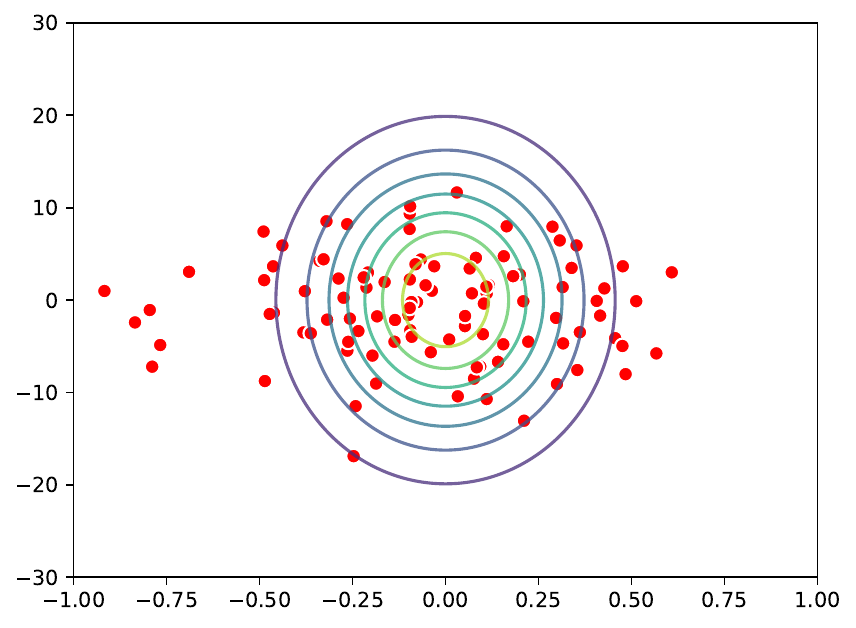}
             \caption{ol}\label{fig:2d_gauss_ol_large}
        \end{subfigure}
        \begin{subfigure}[b]{0.24\linewidth}   
            \centering 
            \includegraphics[scale=0.21]{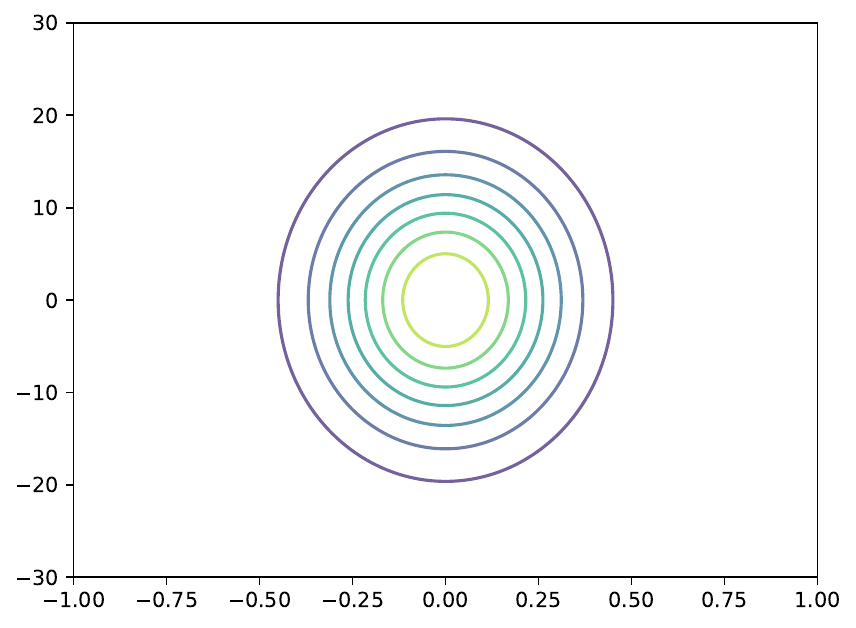}
            \caption{ul}\label{fig:2d_gauss_ul_large}
        \end{subfigure}
        \hfill
        \begin{subfigure}[b]{0.24\linewidth}   
            \centering 
            \includegraphics[scale=0.21]{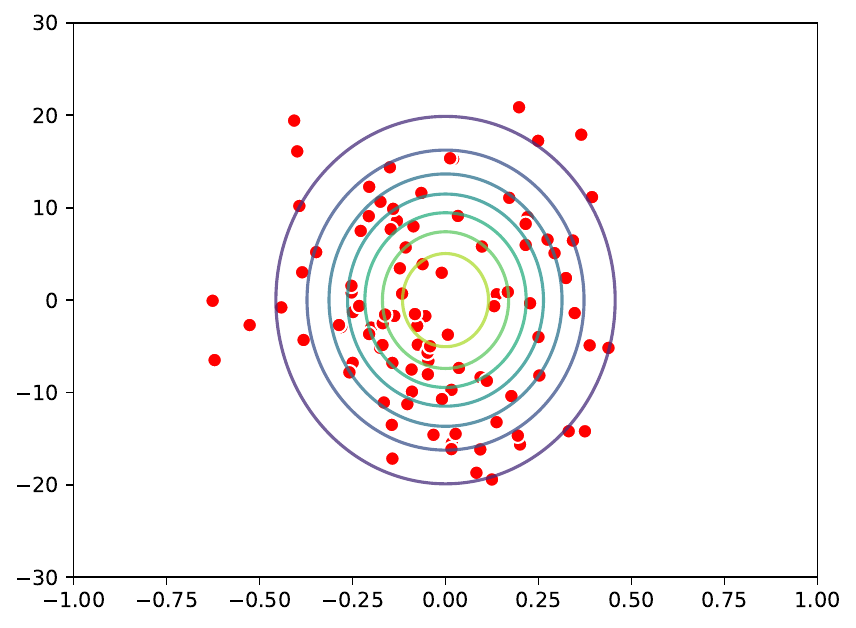}
            \caption{gaul}\label{fig:2d_gauss_pdd_large}
        \end{subfigure}
         \caption{Convergence and scatter plots. (a)--(d): $h=0.005$. (e)--(h): $h=0.05$. (a) and (e): KL divergence between the sample and target distribution. The x-axis represents time and the y-axis is in $\log_{10}$ scale. Rest panels: scatter plot of the three methods at the end of the experiment for different step sizes. Contours of the true density are also provided for comparisons. In (g) there are no scatter points shown as `ul' does not converge for this choice of $h$. }
        \label{fig:2d_gauss}
    \end{figure*}

Let the target distribution be a 20-dimensional Gaussian with zero mean and covariance given by a diagonal matrix with entries $0.05 + 5i$ for $i=0,\ldots,19$. The last dimension has the largest variance, which is $\sigma_{\mathrm{max}}^2 = 95.05$. Therefore, we choose $a =\frac{2}{\sigma_{\mathrm{min}}^{-1/2} -\sigma_{\mathrm{max}}^{-1/2} } $, $\gamma_{ul} = 2\sigma_{\mathrm{max}}^{-1}$ and $\gamma_{pdd} =2\sigma_{\mathrm{max}}^{-1} + a \sigma_{\mathrm{max}}^{-2} $. In this experiment, we use (1) time step size $h=5\times 10^{-3}$ and run for 4000 steps; (2) time step size $h=5\times 10^{-2}$ and run for 400 steps. The KL divergence can still be computed using \cref{eq:KL for gaussian}. To visualize the final distribution in a two-dimensional plane, we plot the scatter plot of the samples in the first and the last dimensions. All results are presented in \Cref{fig:2d_gauss}. 

\subsection{Mixture of Gaussian}
\subsubsection{Strongly log-concave}\label{example:mix_gauss_log_concave}
\begin{figure*}
        \centering
        \begin{subfigure}[b]{0.24\linewidth}
            \centering
            \includegraphics[scale=0.21]{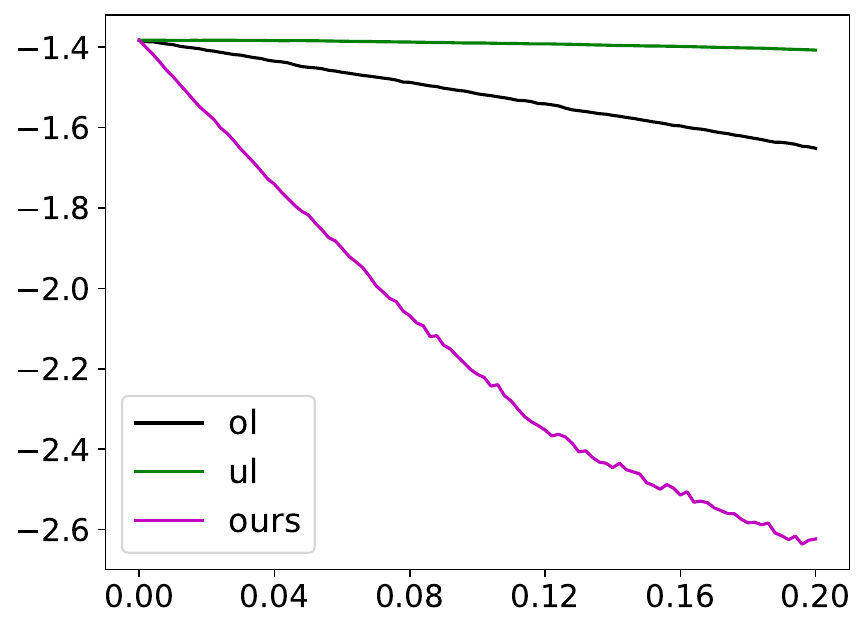}
            \caption{KL decay}\label{fig:mix_kl}
        \end{subfigure}
        \hfill
        \begin{subfigure}[b]{0.24\linewidth}  
            \centering 
            \includegraphics[scale=0.21]{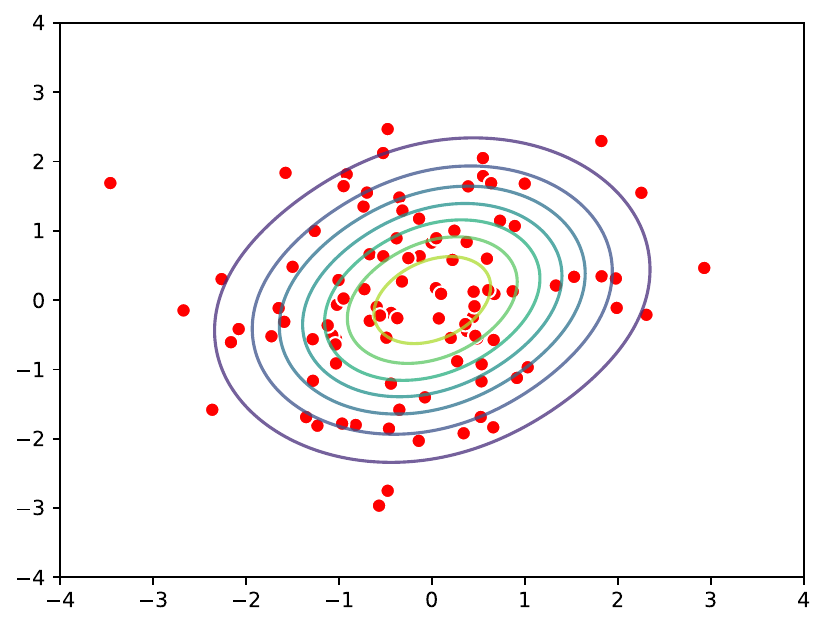}
             \caption{ol}\label{fig:mix_ol}
        \end{subfigure}
        \begin{subfigure}[b]{0.24\linewidth}   
            \centering 
            \includegraphics[scale=0.21]{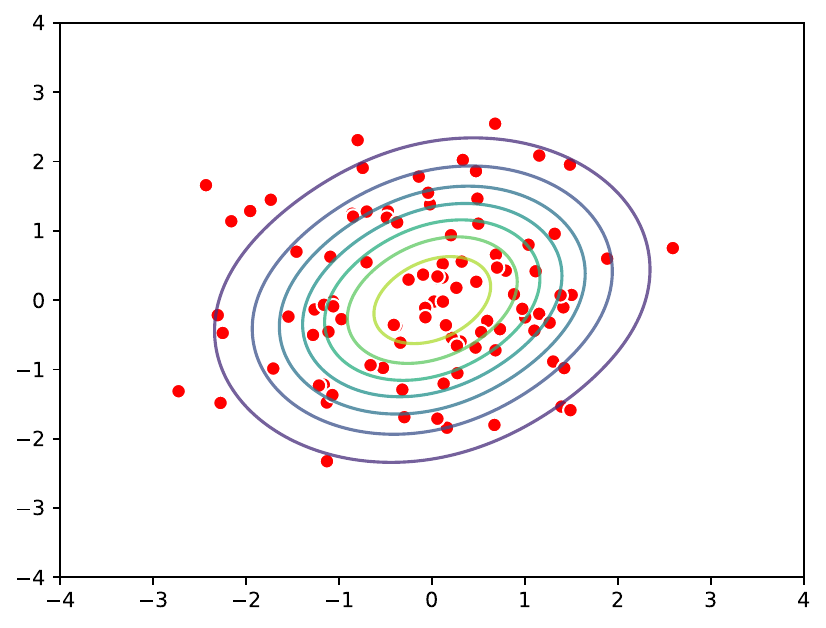}
            \caption{ul}\label{fig:mix_ul}
        \end{subfigure}
        \hfill
        \begin{subfigure}[b]{0.24\linewidth}
            \centering 
            \includegraphics[scale=0.21]{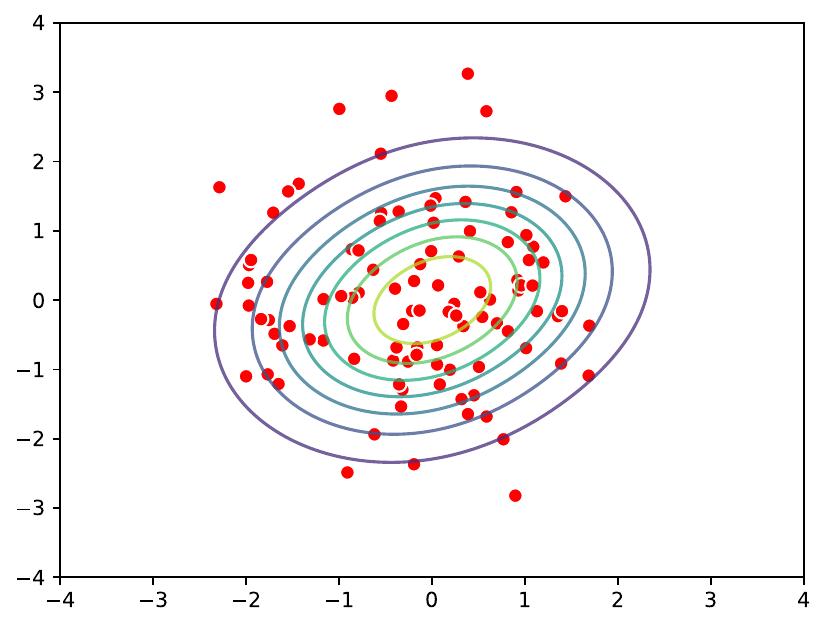}
            \caption{gaul}\label{fig:mix_pdd}
        \end{subfigure}
         \caption{Convergence and scatter plots. (a): KL divergence between the sample and target distribution, which is a mixture of two unit variance Gaussians located at $(1/2,1/2)$ and $(-1/2,-1/2)$. x-axis represents time and y-axis is in $\log_{10}$ scale. (b)--(d): scatter plot of the three methods a the end of the experiment. Contour of the true density is also provided for comparison. }
        \label{fig:mix}
    \end{figure*}

Consider the problem of sampling from a mixture of Gaussian distributions $\mathcal{N}(\alpha, \bI)$ and $\mathcal{N}(-\alpha, \bI)$, whose density satisfies:
\[
p(\vx) = \frac{1}{2(2\pi)^{d/2}} \left( e^{-\|\vx-\alpha\|_2^2/2} + e^{-\|\vx+\alpha\|_2^2/2} \right).
\]
The corresponding potential is given as 
\begin{equation}
f(\vx) = \frac{1}{2} \|\vx - \alpha\|_2^2 - \log \left( 1 + e^{-2\vx^\top \alpha} \right),
\end{equation}
\begin{equation}
\nabla f(\vx) = \vx - \alpha + 2\alpha (1 + e^{2\vx^\top \alpha})^{-1}.
\end{equation}
Following \cite{dwivedi2019log,dalalyan2017theoretical}, we set $\alpha=(1/2,1/2)$ and $d=2$. This choice of parameters yields strong convexity parameter $m=1/2$ and Lipschitz constant $L=1$. We choose $a=\frac{2}{\sqrt{L}-\sqrt{m}}$, $\gamma_{ul} =2m^{1/2}$ and $\gamma_{pdd}=2m^{1/2}+am$. Initially particles are sampled from $\mathcal{N}(0,\bI)$. We use time step $h=2\times 10^{-4}$ and run for 2000 steps. We use $5\times 10^5$ particles and $n^2=2500$ bins to approximate the KL divergence between the sample points and the target distribution (see \Cref{rmk:non-gaussian kl}). The results are shown in \Cref{fig:mix}. 

\begin{remark}\label{rmk:non-gaussian kl}
To compute the KL divergence between sample points and a non-Gaussian target distribution in two dimension, we first get the 2d histogram of the samples points using $n^2$ bins ($n$ in each dimension). We then use this 2d histogram as an approximation of the empirical distribution of the samples. Similarly, we can get a discretized target distribution by evaluating the target distribution at the center of each bins. Finally, we can compute the discrete KL divergence using $n^2$ values from the histogram and the discretized target distribution. 
\end{remark}

\subsubsection{Non log-concave }
\begin{figure*}
        \centering
        \begin{subfigure}[b]{0.24\linewidth}
            \centering
            \includegraphics[scale=0.21]{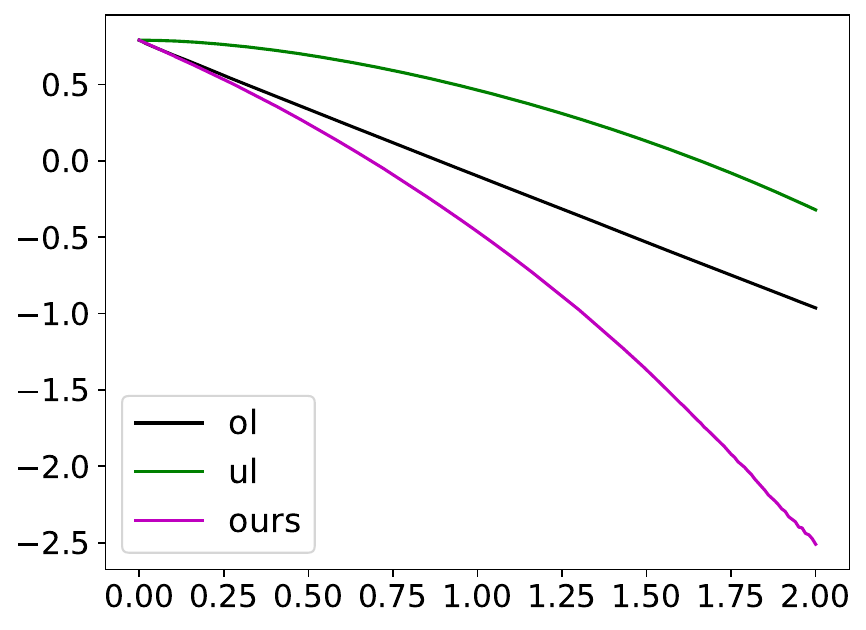}
            \caption{KL decay}\label{fig:mix3_kl}
        \end{subfigure}
        \hfill
        \begin{subfigure}[b]{0.24\linewidth}  
            \centering 
            \includegraphics[scale=0.21]{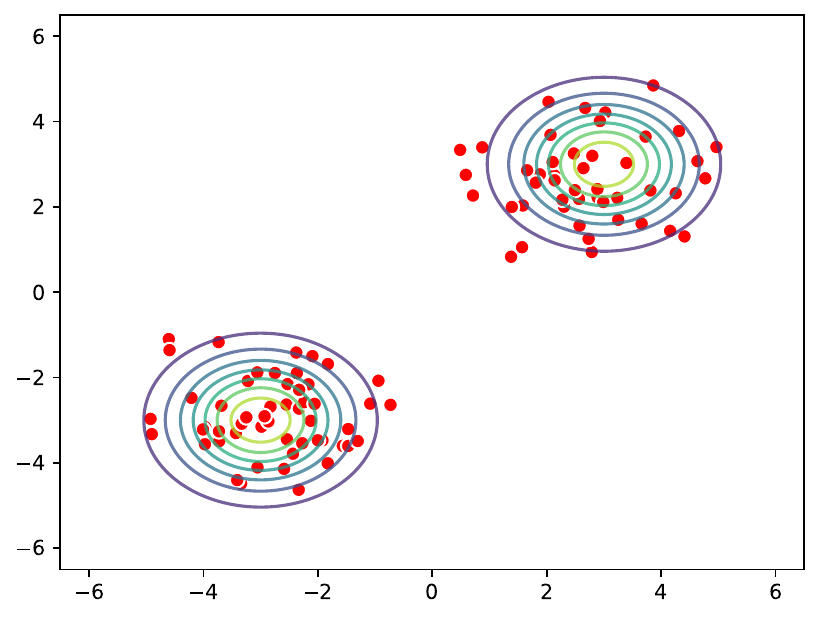}
             \caption{ol}\label{fig:mix3_ol}
        \end{subfigure}
        \begin{subfigure}[b]{0.24\linewidth}   
            \centering 
            \includegraphics[scale=0.21]{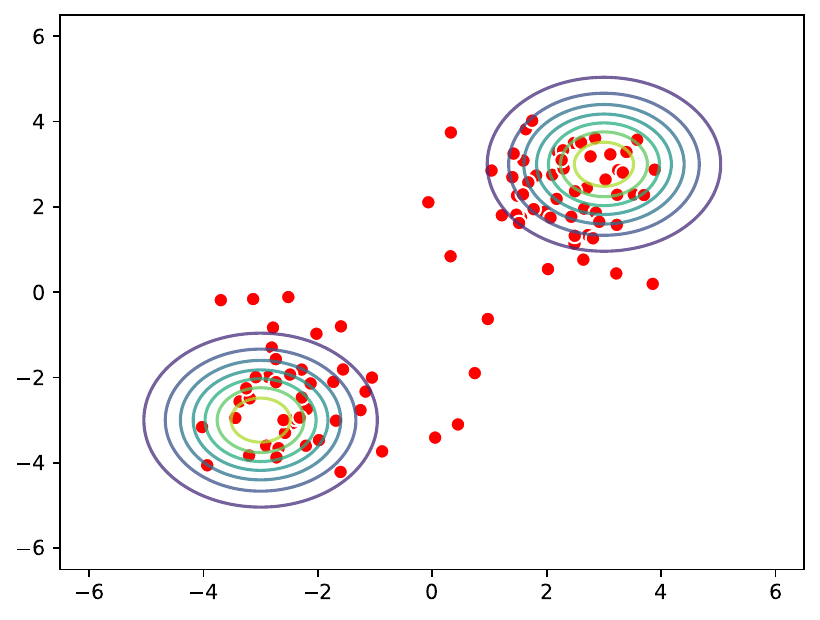}
            \caption{ul}\label{fig:mix3_ul}
        \end{subfigure}
        \hfill
        \begin{subfigure}[b]{0.24\linewidth}   
            \centering 
            \includegraphics[scale=0.21]{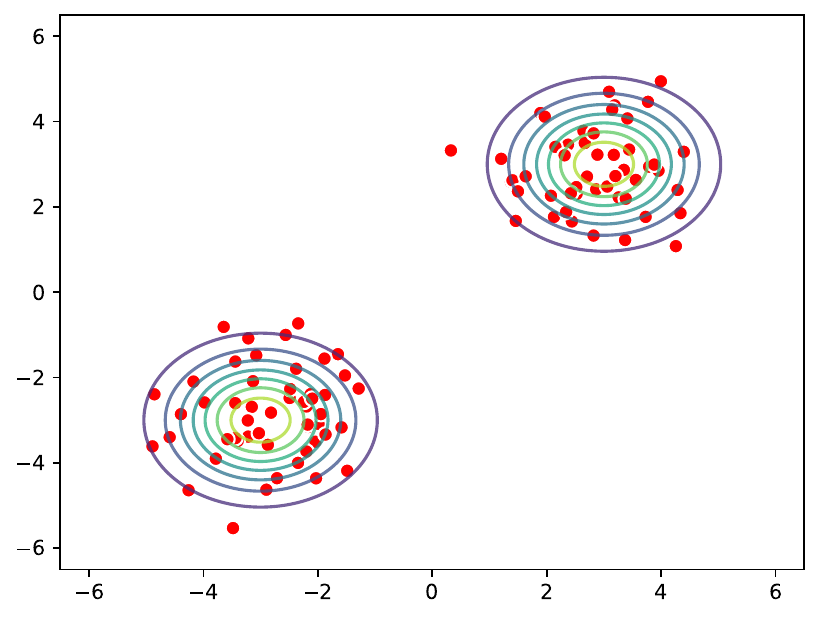}
            \caption{gaul}\label{fig:mix3_pdd}
        \end{subfigure}
         \caption{Convergence and scatter plots for mixture of Guassians centered at $(3,3)$ and $(-3,-3)$.}
        \label{fig:mix3}
    \end{figure*}

We also consider the same example as in \Cref{example:mix_gauss_log_concave} with $\alpha = (3,3)$. As the distance between the two Gaussians increases, the target density is no longer log-concave. We use time step size $h=10^{-3}$ and run for 2000 steps. We use $a=1$, $\gamma_{ul}=\sqrt{2}$, and $\gamma_{pdd}=\sqrt{2}+1/2$. We use $5\times 10^5$ particles and $n^2=2500$ bins to evaluate the KL divergence. The results are demonstrated in \Cref{fig:mix3}. 

\subsection{Quadratic cosine}
\begin{figure*}
        \centering
        \begin{subfigure}[b]{0.24\linewidth}
            \centering
            \includegraphics[scale=0.21]{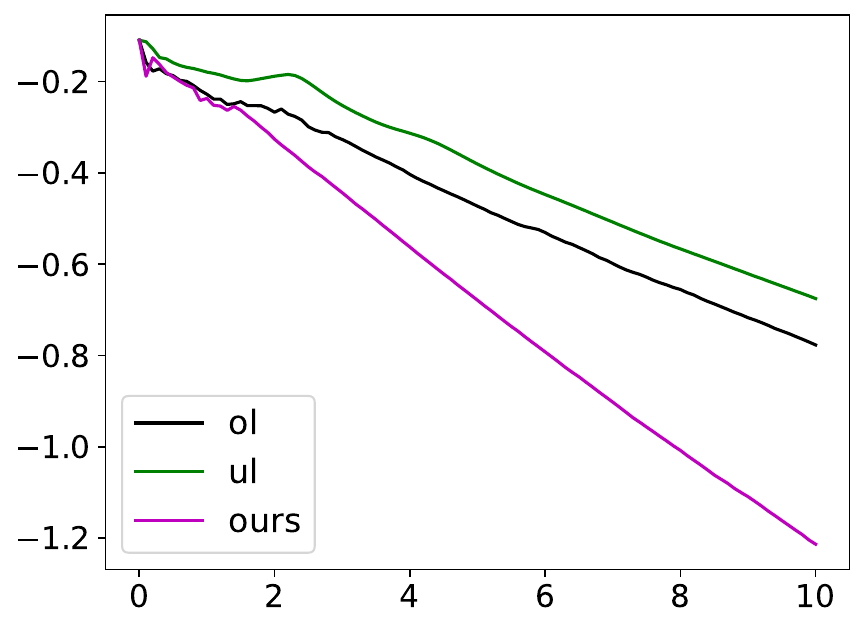}
            \caption{KL decay}\label{fig:qc_kl}
        \end{subfigure}
        \hfill
        \begin{subfigure}[b]{0.24\linewidth}  
            \centering 
            \includegraphics[scale=0.21]{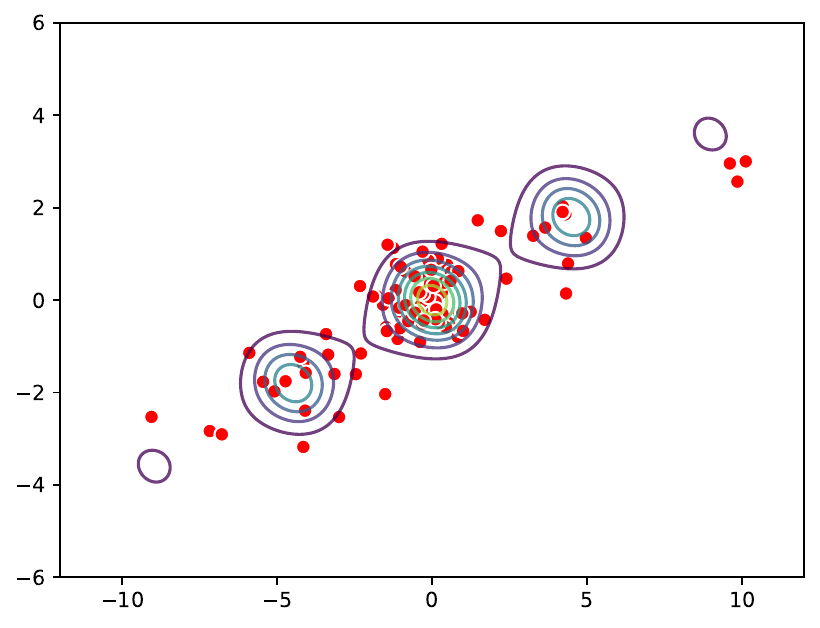}
             \caption{ol}\label{fig:qc_ol}
        \end{subfigure}
        \begin{subfigure}[b]{0.24\linewidth}   
            \centering 
            \includegraphics[scale=0.21]{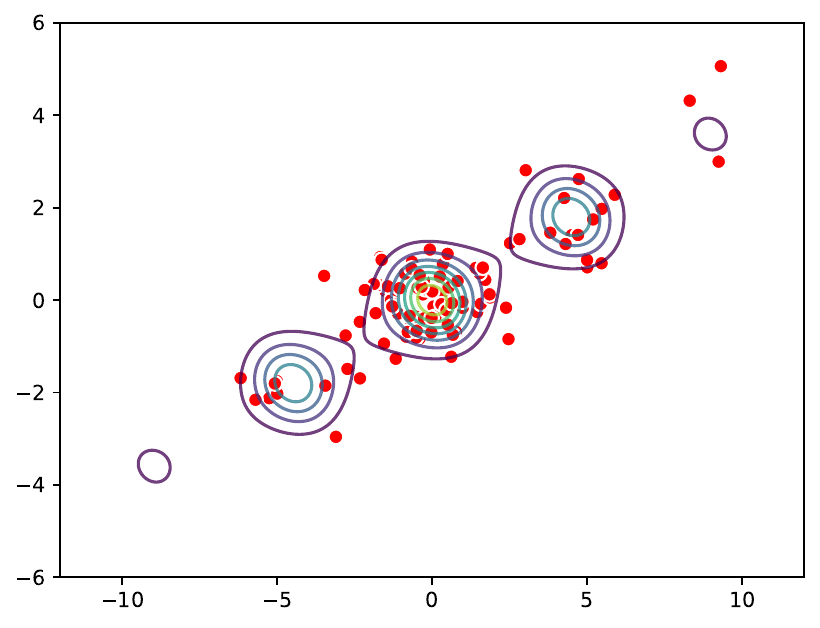}
            \caption{ul}\label{fig:qc_ul}
        \end{subfigure}
        \hfill
        \begin{subfigure}[b]{0.24\linewidth}   
            \centering 
            \includegraphics[scale=0.21]{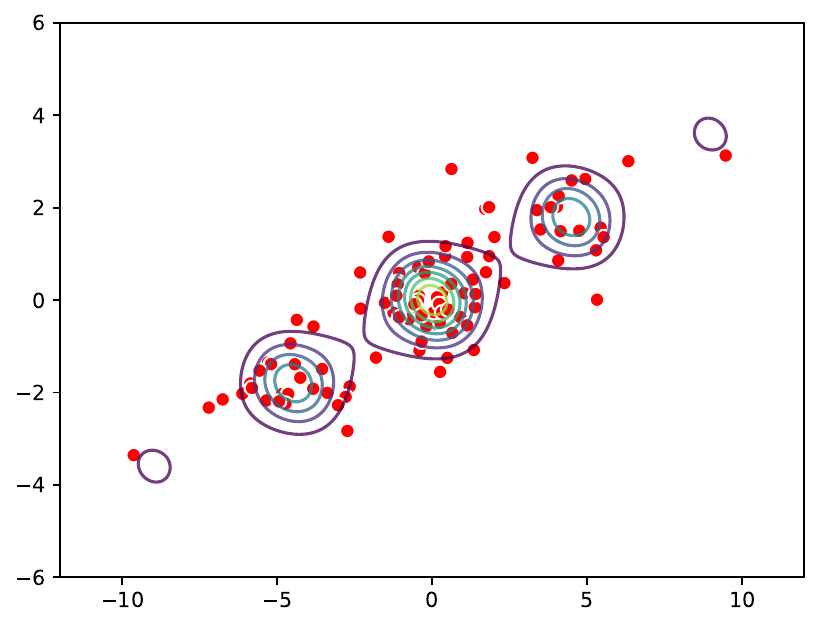}
            \caption{gaul}\label{fig:qc_pdd}
        \end{subfigure}
         \caption{Convergence and scatter plots for the quadratic cosine example.}
        \label{fig:qc}
    \end{figure*}

Consider a potential function given by a quadratic function and a cosine term: 
$$
f(\vx) = \frac{1}{2}\vx^T B^{-1} \vx - \cos (\vc^T \vx) 
$$
where $B = \mP\, \mathrm{diag}(1,25)\, \mP^T$ for an orthogonal matrix $\mP$ and $\vc = \sqrt{0.95}\,(1,1)^T$. Here $\mP$ is generated by using torch.linalg.qr(torch.randn(d)) in Pytorch, where $d=2$ is the dimension.  We set $a=1$, $\gamma_{ul} =2m^{1/2}$ and $\gamma_{pdd}=2m^{1/2}+m$ where we choose $m=1/25$. We use time step size $h=10^{-2}$ and run for 1000 steps. We use $5\times 10^5$ particles and $n^2=2500$ bins to evaluate the KL divergence. The results are demonstrated in \Cref{fig:qc}. 

\subsection{Bimodal}

\begin{figure*}
        \centering
        \begin{subfigure}[b]{0.24\linewidth}
            \centering
            \includegraphics[scale=0.21]{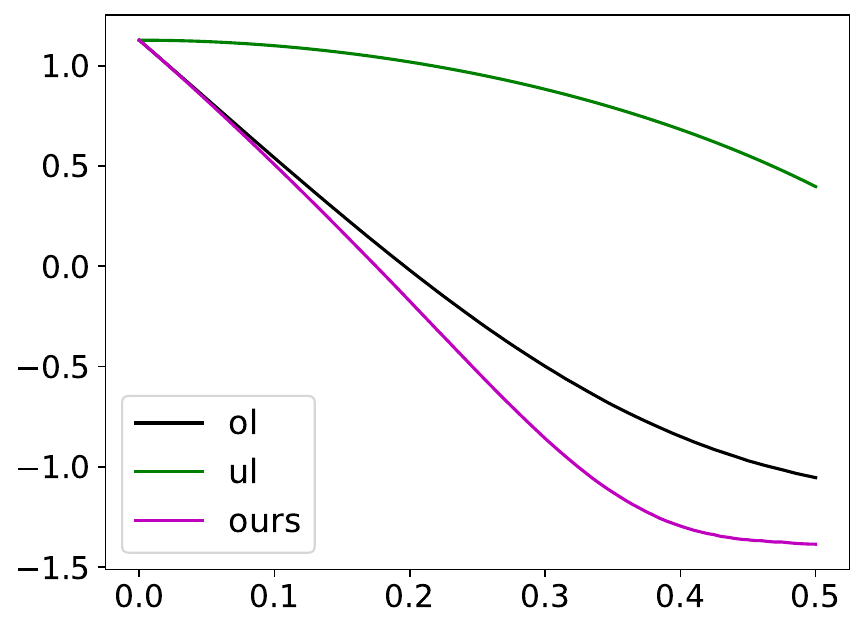}
            \caption{KL decay}\label{fig:bm_kl}
        \end{subfigure}
        \hfill
        \begin{subfigure}[b]{0.24\linewidth}  
            \centering 
            \includegraphics[scale=0.21]{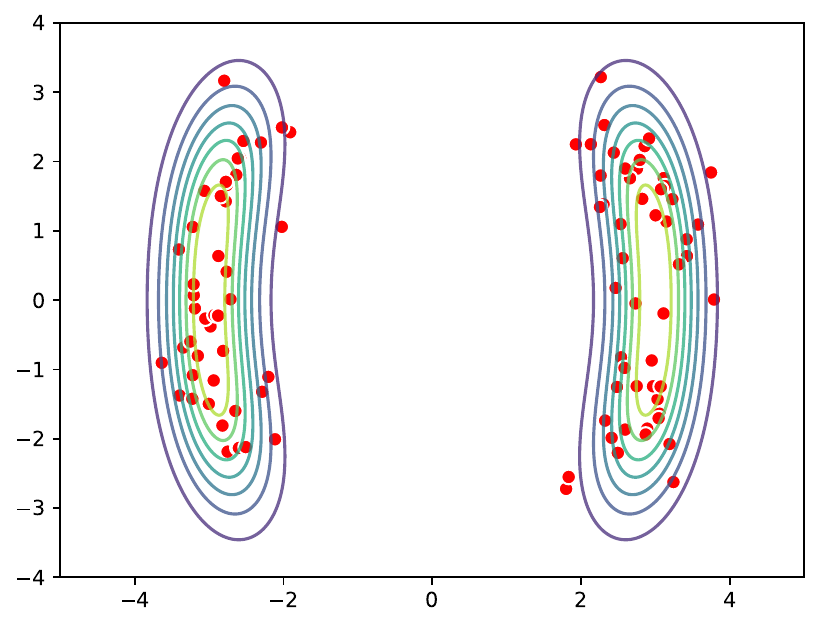}
             \caption{ol}\label{fig:bm_ol}
        \end{subfigure}
        \begin{subfigure}[b]{0.24\linewidth}   
            \centering 
            \includegraphics[scale=0.21]{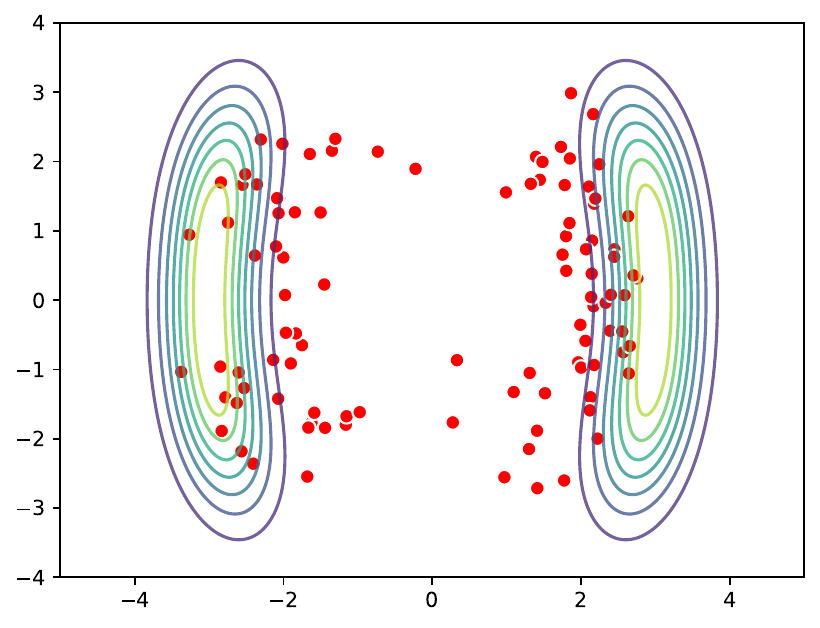}
            \caption{ul}\label{fig:bm_ul}
        \end{subfigure}
        \hfill
        \begin{subfigure}[b]{0.24\linewidth}   
            \centering 
            \includegraphics[scale=0.21]{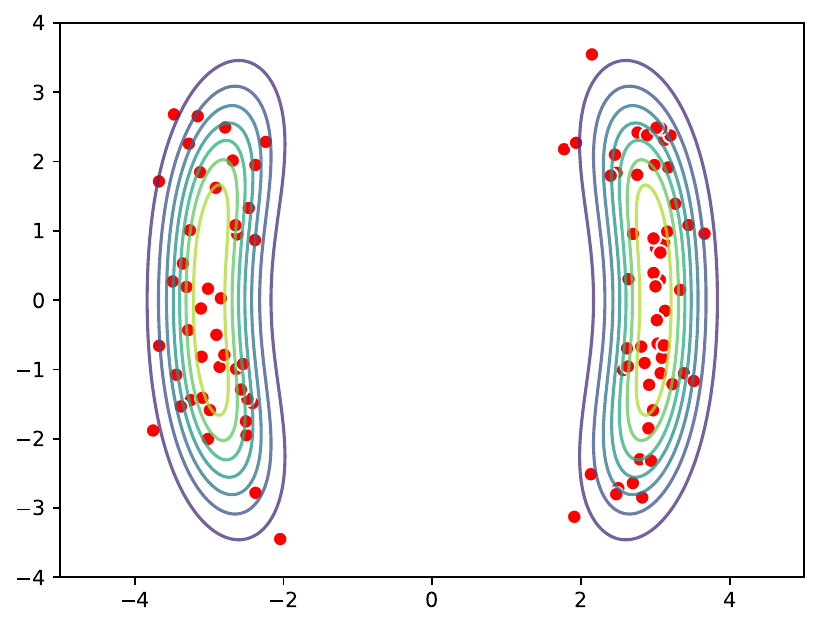}
            \caption{gaul}\label{fig:bm_pdd}
        \end{subfigure}
         \caption{Convergence and scatter plots for the bimodal example.}
        \label{fig:bm}
    \end{figure*}

We consider a two-dimensional bimodal distribution studied in \cite{wang2022accelerated} whose target density has the following form: 
$$
p(\vx) \propto \exp\big(-2(\|\vx\|-3)^2 \big) \Big[\exp\big(-2(x_1-3)^2\big) + \exp\big(-2(x_1+3)^2\big)\Big]\,. 
$$
The corresponding potential function is given by 
$$
f(\vx) = 2(\|\vx\|-3)^2 - 2\log \Big[\exp\big(-2(x_1-3)^2\big) + \exp\big(-2(x_1+3)^2\big)\Big]\,.
$$
The gradient is 
\begin{align*}
    \nabla f(\vx) &=  \frac{4(x_1-3)\exp\big(-2(x_1-3)^2\big) +4(x_1+3)\exp\big(-2(x_1+3)^2\big) }{\exp\big(-2(x_1-3)^2\big) + \exp\big(-2(x_1+3)^2\big)} \ve_1 \\
    &\qquad + 4\frac{(\|\vx\|-3) \vx}{\|\vx\|} \,, 
\end{align*}
where $\ve_1 = (1,0)^T$ is the first standard coordinate vector. We set $\gamma_{ul} =2m^{1/2}$ and $\gamma_{pdd}=2m^{1/2}+m$ where we choose $m=1/2$.
We use time step size $h=10^{-3}$ and run for 500 iterations. We use $10^6$ particles and $n^2=2500$ bins to evaluate the KL divergence. The results are shown in \Cref{fig:bm}.  

\subsection{Bayesian logistic regression}
We consider the Bayesian logistic regression problem studied in \cite{dwivedi2019log,dalalyan2017theoretical,tan2023noise}. We give a brief description of the problem. Suppose we are given a feature matrix $X\in\mathbb R^{n\times d}$ with rows $x_i\in\mathbb R^d$. Correspondingly we are given $Y\in \{0,1\}^n$ the binary response vector for each of the covariates in our feature matrix. The logistic model for the probability of $y_i=1$ given $x_i\in \mathbb R^d$ and a parameter $\theta\in \mathbb R^d$ is 
\begin{equation}\label{eq:generate_y}
\mathbb P(y_i = 1 | x_i,\theta) = \frac{\exp(\theta^T x_i)}{1+\exp(\theta^T x_i)}\,. 
\end{equation}

\begin{figure*}
        \centering
        \begin{subfigure}[b]{0.24\linewidth}
            \centering
            \includegraphics[scale=0.21]{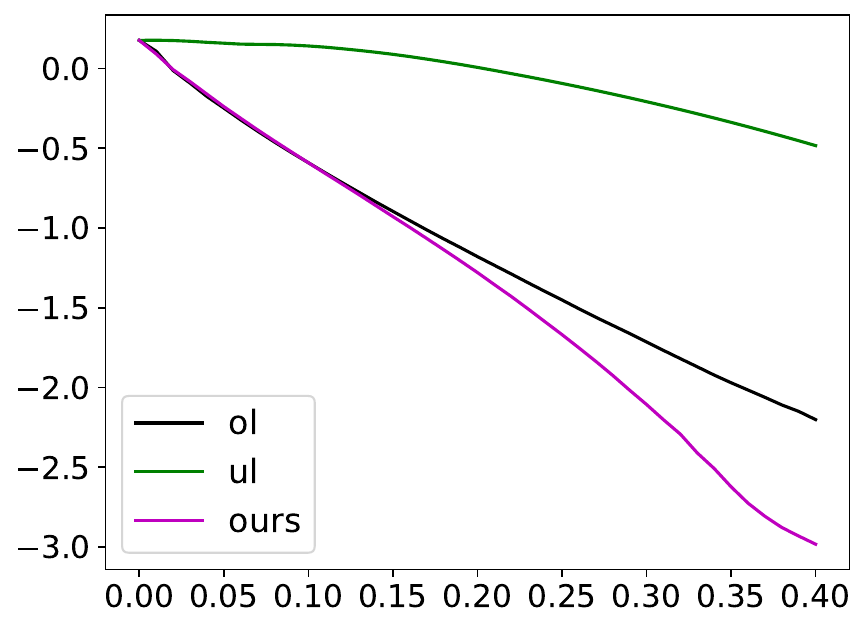}
            \caption{KL decay}\label{fig:bl_kl}
        \end{subfigure}
        \hfill
        \begin{subfigure}[b]{0.24\linewidth}  
            \centering 
            \includegraphics[scale=0.21]{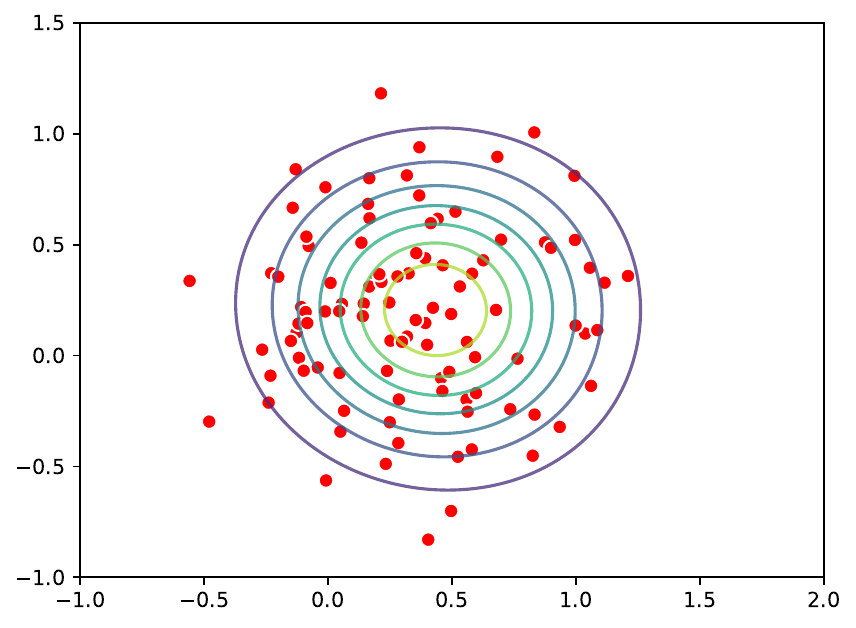}
             \caption{ol}\label{fig:bl_ol}
        \end{subfigure}
        \begin{subfigure}[b]{0.24\linewidth}   
            \centering 
            \includegraphics[scale=0.21]{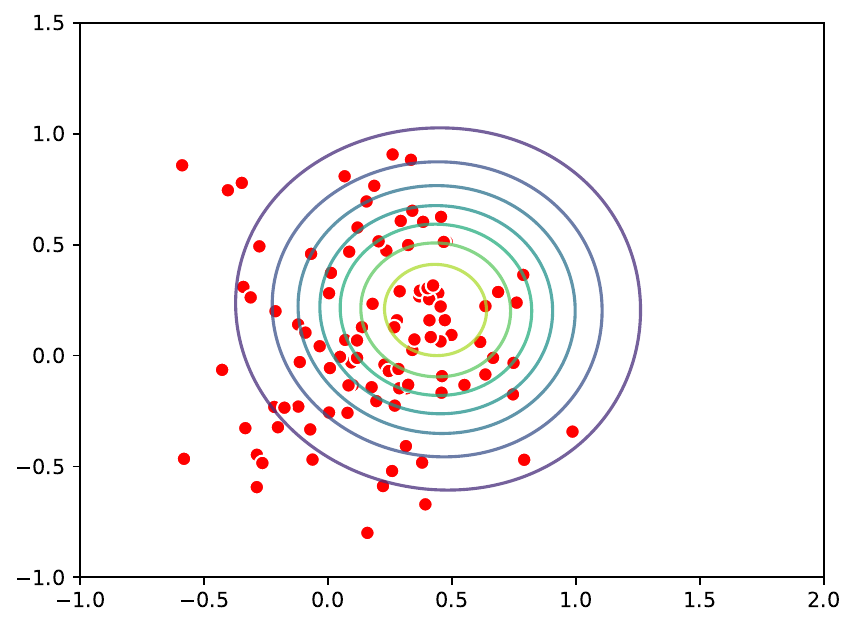}
            \caption{ul}\label{fig:bl_ul}
        \end{subfigure}
        \hfill
        \begin{subfigure}[b]{0.24\linewidth}   
            \centering 
            \includegraphics[scale=0.21]{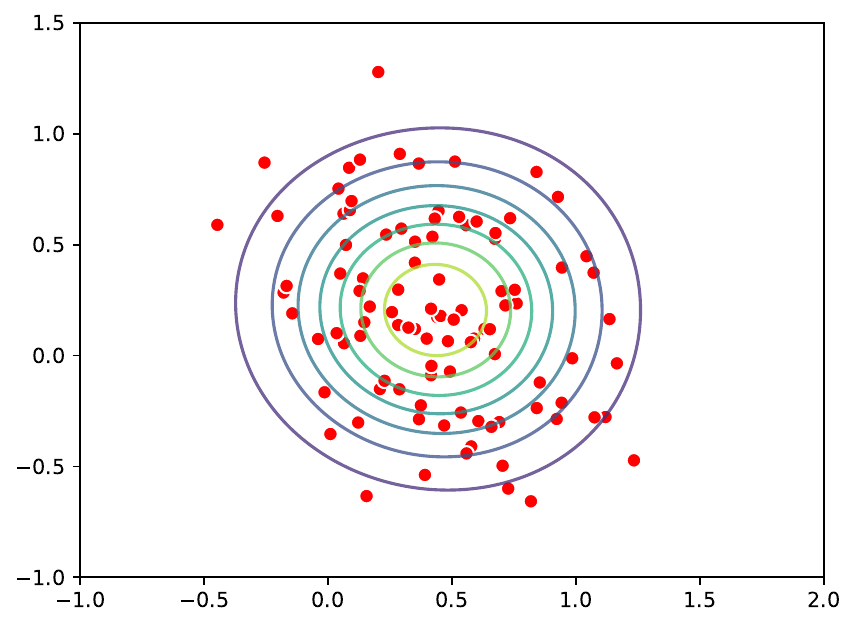}
            \caption{gaul}\label{fig:bl_pdd}
        \end{subfigure}
         \caption{Convergence and scatter plots for Bayesian logistic regression.}
        \label{fig:bl}
    \end{figure*}

Suppose we impose a prior distribution on the parameter $\theta \sim \mathcal{N}(0,\Sigma_X)$, where $\Sigma_X = \frac{1}{n}X^T X$ is the sample covariance of $X$. Then the posterior distribution for $\theta$ can be calculated by 
$$
p(\theta|X,Y) \propto \exp \Big[Y^T X \theta - \sum_{i=1}^n \log\big(1+\exp(\theta^T x_i)\big) - \frac{\alpha}{2} \theta^T \Sigma_X \theta \Big]\,,
$$
where $\alpha>0$ is a regularization parameter. The potential function is 
$$
f(\theta)= -Y^T X \theta + \sum_{i=1}^n \log\big(1+\exp(\theta^T x_i)\big) + \frac{\alpha}{2} \theta^T \Sigma_X \theta\,.
$$
Its gradient is 
$$
\nabla f(\theta) = -X^T Y + \sum_{i=1}^n \frac{x_i}{1+\exp(-\theta^T x_i)} + \alpha \Sigma_X \theta \,. 
$$
As shown in \cite{dwivedi2019log}, the Hessian of $f$ is upper bounded by $L = (0.25n+\alpha)\lambda_{\mathrm{max}}(\Sigma_X)$ and lower bounded by $m =\alpha \lambda_{\mathrm{min}}$. To generate $X$ and $Y$, we set $x_{i,j} $ to be independent Rademacher random variables for each $i$ and $j$. And each $y_i$ is generated according to \cref{eq:generate_y} with $\theta= \theta^* = (1,1)^T$. We set $\alpha=0.5$, $d=2$,  $n=50$, $\gamma_{ul} =2m^{1/2}$ and $\gamma_{pdd}=2m^{1/2}+m$.  To sample the posterior distribution, 
we use time step size $h=10^{-3}$ and run for 400 iterations. The initial distribution of particles is $\mathcal{N}(0,L^{-1}\bI)$. As for evaluation metric, we directly evaluate the KL divergence between the sampled posterior and the true posterior. We use $10^6$ particles and $n^2=2500$ bins to evaluate the KL divergence as before. This is different from the choice by \cite{dwivedi2019log} and \cite{tan2023noise}, where \cite{dwivedi2019log} compared the samples with $\theta^*$. \cite{tan2023noise} compared samples with the true minimizer of $f(\theta)$, i.e. the maximum a posteriori (MAP) estimate in the Bayesian optimization literature. We believe that directly measuring the KL divergence gives a better understanding of how `close' our samples are to the true posterior distribution. The results are presented in \Cref{fig:bl}. 

\subsection{Bayesian neural network}\label{sec:BNN}

\begin{figure*}
        \centering
        \begin{subfigure}[b]{0.475\textwidth}
            \centering
            \includegraphics[width=\textwidth]{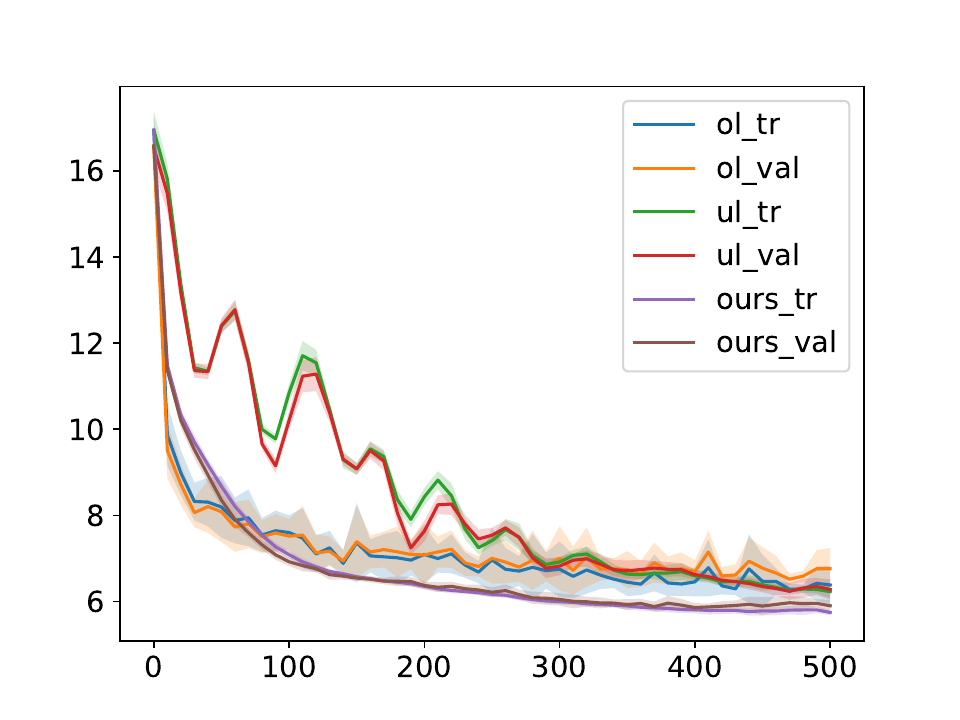}
            \caption{Concrete dataset}\label{fig:concrete}
        \end{subfigure}
        \hfill
        \begin{subfigure}[b]{0.475\textwidth}  
            \centering 
            \includegraphics[width=\textwidth]{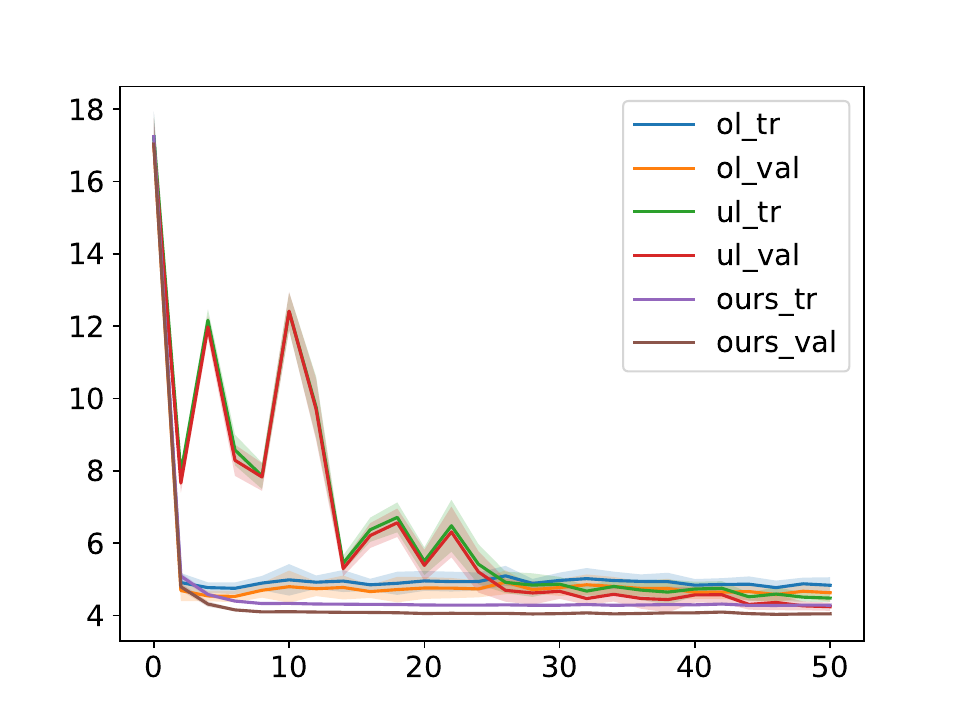}
             \caption{CCPP dataset}\label{fig:ccpp}
        \end{subfigure}
         \caption{Convergence comparison. x-axis represents number of epochs. y-axis represents rMSE averaged over 10 experiments. }
        \label{fig:bnn}
    \end{figure*}

In this section, we compare GAUL with overdamped (`ol') and underdamped Langevin (`ul') dynamics in training Bayesian neural network. We test a one-hidden-layer fully connected neural network with 50 hidden neurons and ReLU activation function on the UCI concrete dataset. We use $h=10^{-3}$, $a=0.1$, $\gamma=0.5$. For each method, we sample $M=20$ particles (each particle corresponds to a neural network) and take the average output as the final output. \Cref{fig:concrete} and \Cref{tab:err} show the rMSE averaged over 10 experiments. We see that `ul' can achieve smaller training and validation error than `ol'. However, `ul' also exhibits a slow start and an oscillatory behavior at the beginning of training as is commonly seen in acceleration methods in optimization. GAUL can get rid of the oscillation and achieve a even smaller training and validation error as is demonstrated in \Cref{tab:err}. We have also tested out the three methods using the Combined Cycle Power Plant (CCPP) dataset. We choose the same parameter as the concrete experiment. The results are presented in \Cref{fig:ccpp} and \Cref{tab:err}. 
\begin{table}
\begin{center}
\begin{tabular}{|c|c|c|c|}
            \hline
             & ol & ul & gaul \\
            \hline\hline
            concrete tr err  & $6.39\pm0.44$  & $6.23\pm0.15$  & $\boldsymbol{5.74\pm0.06}$  \\
            \hline
            concrete val err  & $6.76\pm0.49$  & $6.28\pm0.24$  & $\boldsymbol{5.90\pm0.14}$ \\ 
           \hline\hline 
           ccpp tr err & $4.84\pm0.22 $ & $4.48\pm0.11$ & $\boldsymbol{4.28 \pm0.03}$ \\
           \hline 
           ccpp val err & $4.63 \pm 0.25$ & $4.25\pm0.11$ & $\boldsymbol{4.04\pm0.04}$ \\
           \hline 
        \end{tabular}
\end{center}
        \caption{Training and validation rmse. }\label{tab:err}
\end{table}

\section{Conclusions}
In this work, we introduced gradient-adjusted underdamped Langevin dynamics (GAUL) inspired by primal-dual damping dynamics and Hessian-driven damping dynamics. We demonstrated that GAUL admitted the correct stationary target distribution $\pi \propto \exp(-f)$ under appropriate conditions and achieves exponential convergence for quadratic functions, outperforming both the overdamped and underdamped Langevin dynamics in terms of convergence speed. Our numerical experiments further illustrate the practical advantages of GAUL, showcasing faster convergence and more efficient sampling compared to classical methods, such as overdamped and underdamped Langevin dynamics. 

{We also note a connection between the primal-dual damping dynamics and GAUL. A key challenge in the primal-dual damping algorithm is the design of preconditioner matrices, which can accelerate the algorithm’s convergence compared to the gradient descent method. In the context of solving a linear problem where  f  is a quadratic function and the diffusion constant is zero, \cite{zuo2023primal} demonstrates that the convergence rate depends on the square root of the smallest eigenvalue. In this paper, we extend the study from a sampling perspective, where $f$ is also a quadratic function but the diffusion is non-zero. Towards a Gaussian target distribution, GAUL converges to a biased target distribution with the mixing time depending on $\sqrt{\kappa}$. This is in contrast with overdamped and underdamped Langevin sampling algorithms.}

Several possible future directions are worth exploring. First, can we show that GAUL converges faster than overdamped and underdamped Langevin dynamics for more general potentials, which is beyond the current study of Gaussian distributions? One common assumption is that the potential $f$ is strongly log-concave \cite{cao2020complexity,chewi2022query,chewi2021optimal,dalalyan2017further,durmus2019analysis,dwivedi2019log,he2020ergodicity,lee2020logsmooth,li2022hessian}. Recently, \cite{cao2023explicit} proved that for a class of distributions that satisfy a Poincar\'e-type inequality, underdamped Langevin dynamics converges in $L_2$ with rate $\exp(-\sqrt{m}t)$ where $m$ is the Poincar\'e constant. Then it is interesting to study for the same class of distributions, whether GAUL could converge at an even faster rate. Another direction is to study the convergence of GAUL under different metrics. From a more practical perspective, designing new time discretization schemes \cite{shen2019randomized,cheng2018underdamped,mou2019high,tan2023noise,li2022hessian} for implementing GAUL is also an important direction. We proved that using the Euler-Maruyama discretization, GAUL will converge to a biased target distribution, which is not surprising since ULA is also biased. Therefore, another promising direction could be to combine GAUL with MCMC methods \cite{besag1994comments,dwivedi2019log}, such as Metropolis-Hastings algorithms, to design a hybrid method with accept/reject options so that the algorithm converges to the correct target distribution in the discrete-time update. Finally, choosing the preconditioner $\mC$ to accelerate convergence is an important topic. The difficulty of picking $\mC$ arises from the positive semidefinite constraint on $\mathrm{sym}(\bQ)$ in \cref{eq:general_Q}, which we should explore in future work.
\bibliography{main}

\newpage

\begin{appendices}
    
\section{Euler-Maruyama Discretization}\label{sec:euler_maruyama}
The Euler-Maruyama scheme of \cref{eq:SDE} with step size $h$ and $\mC = \bI$ reads 
\begin{subequations}\label{eq:Euler_dis}
\begin{align}
    \vx_{t+1} &= \vx_t - a \nabla f(\vx_t) h + \vp_t h + \sqrt{2ah}\vz^{(1)}\,,  \\
    \vp_{t+1} &= \vp_t - \nabla f(\vx_t) h - \gamma \vp_t h + \sqrt{2 \gamma  h} \vz^{(2)} \,.
\end{align}
\end{subequations}
$\vz^{(i)}$ is a standard Gaussian random variable for $i=1,2$. 

\section{A matrix lemma}
Let $a\geq0$, $s>0$, $\gamma>0$, and consider the $3\times3$ matrix 
\begin{equation}\label{eq:D}
    \mathbf{D} = \begin{pmatrix}
        -2as&-2\gamma s^{-1} & -s^{-1}\\
        0&0&1 \\
        2s^2 &2(-1-a\gamma)s^{-1}-2\gamma^2 &-3\gamma - as 
    \end{pmatrix}\,.
\end{equation}
A direct calculation shows that the eigenvalues are given by
\begin{subequations}\label{eq:lambda_pm}
    \begin{align} 
    \lambda_0(a,\gamma,s) &= -a s - \gamma \,, \\
    \lambda_-(a,\gamma,s) &= -as - \gamma  - \sqrt{\gamma^2-2a\gamma s + s(-4+a^2 s)}\,, \\
    \lambda_+(a,\gamma,s) &= -as - \gamma  + \sqrt{\gamma^2-2a\gamma s + s(-4+a^2 s)}\,.
\end{align}
\end{subequations}
We have the following lemmas regarding the eigenvalues given by \cref{eq:lambda_pm}. 
\begin{lemma}\label{lemma:matrix_1}
Let $\mathbf{D}$ be as \cref{eq:D}. If $a=0$, then 
$$
\argmin_{\gamma>0} \Re\big(\lambda_+(0,\gamma,s)\big) = 2\sqrt{s}\,.
$$
\end{lemma}
\begin{proof}
We have that $\lambda_+(0,\gamma,s) =\frac{1}{2}\left(-\gamma +\sqrt{\gamma^2-4s}\right)$. If $\gamma \leq 2\sqrt{s}$, then $\Re\big(\lambda_+(0,\gamma,s)\big)\geq -\sqrt{s}$. When $\gamma \geq 2\sqrt{s}$, we have that $\Re\big(\lambda_+(0,\gamma,s)\big) =\lambda_+(0,\gamma,s) $. And $\frac{\partial }{\partial \gamma}\lambda_+(0,\gamma,s) \geq 0 $. Therefore, the minimum of $\Re\big(\lambda_+(0,\gamma,s)\big)$ takes place at $\gamma = 2\sqrt{s}$.    
\end{proof}

\begin{lemma}\label{lemma:matrix_2}
Let $\mathbf{D}$ be as \cref{eq:D}. Let $\gamma>0$ be fixed. Then 
\begin{equation}
    \argmin_{a>0} \Re\big(\lambda_+(a,\gamma,s)\big) = \frac{\gamma}{s} + \frac{2}{\sqrt{s}}\,. 
\end{equation}
\end{lemma}
\begin{proof}
Let us define $\Delta(a) = \gamma^2 -2a\gamma s + s(a^2 s -4)$. It can be seen that $\Delta$ is a quadratic function of $a$. The two roots of $\Delta$ are given by 
$$
a_{\pm} = \frac{\gamma}{s} \pm \frac{2}{\sqrt{s}}\,. 
$$
When $a \in [a_{-},a_+]$, $\Delta(a)\leq 0$ and $$\Re\big(\lambda_+(a,\gamma,s)\big)=\frac{1}{2}(-\gamma-as) \geq \frac{1}{2}(-\gamma-a_+s) =\Re\big(\lambda_+(a_+,\gamma,s)\big)= -\gamma-\sqrt{s}\,. $$
When $a < a_-$, we can calculate that 
$$
\frac{\partial}{\partial a}\lambda_+(a,\gamma,s) = -s + \frac{-\gamma s + a s^2}{\sqrt{\Delta}} < 0 \,.
$$
This implies that $\lambda_+(a_- - \varepsilon,\gamma,s)>\lambda_+(a_-,\gamma,s)$ for any $\varepsilon>0$. Similarly, when $a > a_+$, we have that $\frac{\partial}{\partial a}\lambda_+(a,\gamma,s) >0$. Thus, $\lambda_+(a_+ + \varepsilon,\gamma,s)>\lambda_+(a_-,\gamma,s)$ for any $\varepsilon>0$. Combining the above results, we conclude our proof.  
\end{proof}

\begin{lemma}\label{lemma:matrix_gamma_single}
Let $\mathbf{D}$ be as \cref{eq:D}. Let $a>0$ be fixed. Then 
\begin{equation}
    \argmin_{\gamma>0} \Re\big(\lambda_+(a,\gamma,s)\big) =as + 2\sqrt{s}\,. 
\end{equation}
\end{lemma}
\begin{proof}
The proof will be similar to that of \Cref{lemma:matrix_2}. This time we define $\Delta(\gamma) =\gamma^2 -2a\gamma s + s(a^2 s -4)$. It can be seen that $\Delta(\gamma)$ is a quadratic function of $\gamma$. The two roots of $\Delta(\gamma)$ are given by 
$$
\gamma_{\pm} = as \pm 2\sqrt{s}\,. 
$$
When $\gamma \in [\gamma_-,\gamma_+], \Delta(\gamma)<0$ and 
$$\Re\big(\lambda_+(a,\gamma,s)\big)=\frac{1}{2}(-\gamma-as) \geq \frac{1}{2}(-\gamma_+-as) =\Re\big(\lambda_+(a,\gamma_+,s)\big)= -as-\sqrt{s}\,. $$
When $\gamma <\gamma_-$, we have 
\begin{align*}
    \frac{\partial }{\partial \gamma}\lambda_+(a,\gamma,s)&= -1 + \frac{\gamma-as}{\sqrt{(\gamma-as)^2-4s}} \\
    &\leq -1<0 \,,
\end{align*}
since $\gamma-as < 0$. When $\gamma > \gamma_+$, we have 
\begin{align*}
    \frac{\partial }{\partial \gamma}\lambda_+(a,\gamma,s)&= -1 + \frac{\gamma-as}{\sqrt{(\gamma-as)^2-4s}} \\
    &\geq -1 + 1 = 0\,.
\end{align*}
Combining the above arguments, we conclude that the optimal $\gamma$ is $\gamma_+$. 
\end{proof}

We now turn to a more general setting. Let $a\geq 0$, $\gamma>0$ and define 
\begin{equation}\label{eq:D_high_dim}
    \mathcal{D} = \begin{pmatrix}
        -2a\mS & -2\gamma \mS^{-1} & -\mS^{-1}  \\
        0&0&\bI \\
        2 \mS^2 &2(-1-a\gamma)\mS^{-1}-2\gamma^2 \bI &-3\gamma \bI -a\mS 
    \end{pmatrix}\,,
\end{equation}
where now $\mS$ is a diagonal matrix whose diagonal is given by $s_1\geq s_2 \geq \ldots \geq s_d >0$. And $\bI$ is the identity matrix. Just like \Cref{lemma:matrix_1}, \Cref{lemma:matrix_2}, and \Cref{lemma:matrix_gamma} we want to characterize the eigenvalues of $\mathcal{D}$. In particular, we would like to characterize the largest real part of the eigenvalue of $\mathcal{D}$ in terms of $a$ and $\gamma$. 
\begin{proposition}\label{prop:eig_system}
The eigenvalues for $\mathcal{D}$ are given by
\begin{subequations}\label{eq:lambda_pm_hd}
    \begin{align} 
    \lambda^{(i)}_0(a,\gamma,\mS) &= -a s_i - \gamma \,, \\
    \lambda^{(i)}_-(a,\gamma,\mS) &= -as_i - \gamma  - \sqrt{\gamma^2-2a\gamma s_i + s_i(-4+a^2 s_i)}\,, \\
    \lambda^{(i)}_+(a,\gamma,\mS) &= -as_i - \gamma  + \sqrt{\gamma^2-2a\gamma s_i + s_i(-4+a^2 s_i)}\,,
\end{align}
\end{subequations}
for $i=1,\ldots,d$. The corresponding eigenvectors are sparse and take the following form. (Here we only write out the non-zero part of the eigenvectors) 
\begin{subequations}
    \begin{align}
        v^{(i)}_{0,i} &= \frac{-1}{s_i(\gamma + a s_i)}\,, \\
        v^{(i)}_{0,i+d}&= \frac{-1}{\gamma + a s_i}\,,\\
        v^{(i)}_{0,i+2d}&= 1 \,,
    \end{align}
\end{subequations}
\begin{subequations}
    \begin{align}
        v^{(i)}_{-,i} &= \frac{2\gamma-\sqrt{\gamma^2-2a\gamma s_i+s_i(a^2 s_i-4)} -\frac{2(\gamma^2+s_i+a\gamma s_i)}{\gamma + as_i+\sqrt{\gamma^2-2a\gamma s_i+s_i(a^2 s_i-4)}}}{2s_i^2}\,, \\
        v^{(i)}_{-,i+d}&= \frac{-1}{\gamma + as_i+\sqrt{\gamma^2-2a\gamma s_i+s_i(a^2 s_i-4)}}\,,\\
        v^{(i)}_{-,i+2d}&= 1 \,,
    \end{align}
\end{subequations}
\begin{subequations}
    \begin{align}
        v^{(i)}_{+,i} &= \frac{2\gamma+\sqrt{\gamma^2-2a\gamma s_i+s_i(a^2 s_i-4)} -\frac{2(\gamma^2+s_i+a\gamma s_i)}{\gamma + as_i-\sqrt{\gamma^2-2a\gamma s_i+s_i(a^2 s_i-4)}}}{2s_i^2} \,,\\
        v^{(i)}_{+,i+d}&= \frac{-1}{\gamma + as_i-\sqrt{\gamma^2-2a\gamma s_i+s_i(a^2 s_i-4)}}\,,\\
        v^{(i)}_{+,i+2d}&= 1\,.
    \end{align}
\end{subequations}
In the above, $v^{(i)}_{*,j}$ represents the $j$-th component of the eigenvector corresponding to the eigenvalue $\lambda^{(i)}_*$, where $*\in \{0,+,-\}$. 

Moreover, when $\gamma$ is chosen according to \Cref{lemma:matrix_gamma}, we have a defective eigenvalue $\lambda^{(d)}_0 = \lambda^{(d)}_{\pm} = -as_d - \gamma$, which is accompanied with two generalized eigenvectors $\eta$, $\xi$ that satisfy $(\mathcal{D}-\lambda^{(d)}_0)\eta = v^{(d)}_0$, $(\mathcal{D}-\lambda^{(d)}_0)\xi = v^{(d)}_0$. In details, the nonzero components of $v^{(d)}_0$, $\eta$ and $\xi$ are given by 
\begin{subequations}
    \begin{align}
         v^{(d)}_{0,d} &= \frac{-1}{s_d(\gamma + a s_d)}\,, \\
        v^{(d)}_{0,2d}&= \frac{-1}{\gamma + a s_d}\,,\\
        v^{(d)}_{0,3d}&= 1 \,,
    \end{align}
\end{subequations}
\begin{subequations}
    \begin{align}
        \eta_{d} &= \frac{\gamma-as}{2s_d^2}\,, \\
        \eta_{3d} &= 1\,,
    \end{align}
\end{subequations}
\begin{subequations}
    \begin{align}
        \xi_d &= \frac{\gamma^2-(1+a\gamma)s_d}{s_d^2}\,,\\
        \xi_{2d} &= 1\,.
    \end{align}
\end{subequations}
\end{proposition}

\begin{proof}
One can directly verify that the above computation gives an eigensystem for $\mathcal{ D}$. 
\end{proof}
From the sparsity structure of $v^{(j)}_{\pm}$ and $v^{(j)}_{0}$, we immediately have the following corollary. 
\begin{corollary}\label{cor:orthogonal}
$v^{(j)}_{*}$ is orthogonal to $v^{(k)}_{\star}$ for $*,\star\in\{0,+,-\}$ if $j\neq k$. 
\end{corollary}

\begin{lemma}
    Let $\mathcal{D}$ be as \cref{eq:D_high_dim}. If $a=0$, then 
    \begin{equation}
        \argmin_{\gamma >0} \max_{j} \Re(\lambda^{(j)}_+(0,\gamma,\mS)) = 2\sqrt{s_d}\,.
    \end{equation}
\end{lemma}
\begin{proof}
Plugging $a=0$ into \cref{eq:lambda_pm_hd} we have 
$$
\lambda^{(j)}_+(0,\gamma,\mS) = \frac{1}{2}\left(-\gamma +\sqrt{\gamma^2  -4 s_j}\right)\,.
$$
We first note that since $s_d\leq s_{d-1} \leq \ldots \leq s_1$, if $\gamma \leq 2 \sqrt{s_d}$ then $\Re(\lambda^{(j)}_+(0,\gamma,\mS))=-\gamma/2$ for all $1\leq j\leq d$. In particular, this implies that 
$$
\argmin_{0<\gamma \leq 2 \sqrt{s_d}} \max_{j} \Re(\lambda^{(j)}_+(0,\gamma,\mS)) = 2\sqrt{s_d}\,.
$$
We then need to show that if $\gamma > 2\sqrt{s_d}$, $\max_{j} \Re(\lambda^{(j)}_+(0,\gamma,\mS)) > -\sqrt{s_d}$. This will be very similar to the argument in the proof of \Cref{lemma:matrix_1}. Now consider $\gamma>2\sqrt{s_d}$. We showed in the proof of \Cref{lemma:matrix_1} that $\Re\big(\lambda^{(n)}_+(0,\gamma,\mS)\big) =\lambda^{(n)}_+(0,\gamma,\mS) $. And $\frac{\partial }{\partial \gamma}\lambda^{(n)}_+(0,\gamma,\mS) \geq 0 $. Hence, we have 
$$
\max_{j} \Re(\lambda^{(j)}_+(0,\gamma,\mS)) >  \Re(\lambda^{(n)}_+(0,\gamma,\mS))=\lambda^{(n)}_+(0,\gamma,\mS)\geq \lambda^{(n)}_+(0,2\sqrt{s_d},\mS) =-\sqrt{s_d}\,.
$$
This concludes our proof.
\end{proof}

\begin{lemma}\label{lemma:matrix_gamma}
Let $\mathcal{D}$ be as \cref{eq:D_high_dim}. Let $a>0$. Then      \begin{equation}
    \argmin_{\gamma>0} \max_j \Re(\lambda^{(j)}_+(a,\gamma, \mS)) = as_d + 2\sqrt{s_d}\,. 
\end{equation}
\end{lemma}
\begin{proof}
    Let us define $\Delta(\gamma,s) = \gamma^2 -2a\gamma s + s(a^2 s -4)$. A straightforward calculation shows that the two roots of $\Delta(\gamma,s_j)$ (when viewing $\Delta$ as a function of $\gamma$) are given by 
$$
\gamma^{(j)}_{\pm} = a s_j \pm 2\sqrt{s_j}\,.
$$
We have shown in \Cref{lemma:matrix_gamma_single} that 
$$
 \argmin_{\gamma>0}  \Re(\lambda^{(d)}_+(a,\gamma, \mS)) = a s_d + 2\sqrt{s_d}\,.
$$
Denote by $\gamma^*(a) = a s_d + 2\sqrt{s_d}$. Let us consider $\tilde s > s_d$. If $\Delta(\gamma^*(a),\tilde s)\leq 0$, then we have 
\begin{align}\label{eq:ineq0}
    \Re\left(-\gamma^*(a)-a  \tilde s +\sqrt{\gamma^*(a)^2 -2a\gamma^*(a) \tilde s + \tilde s(a^2 \tilde s -4)}\right)&= -\gamma^*(a)-a \tilde s \nonumber  \\
    &\leq -\gamma^*(a)-a  s_d  \nonumber \\
    &= \Re(\lambda^{(d)}_+(a,\gamma^*(a), \mS))  \,,
\end{align}
where the last line follows from $\Delta(\gamma^*(a),s_d) = 0$ by definition of $\gamma^*(a)$. If $\Delta(\gamma^*(a),\tilde s)>0$, we compute 
\begin{align}\label{eq:derivative}
    &\frac{\partial}{\partial s} \left(-\gamma^*(a)-a  s +\sqrt{\gamma^*(a)^2 -2a\gamma^*(a)  s +  s(a^2  s -4)}\right)|_{s=\tilde s}\nonumber \\
    &= -a + \frac{-a \gamma^*(a) + a^2 \tilde 
    s - 2}{\sqrt{\gamma^*(a)^2 -2a\gamma^*(a) \tilde   s +  \tilde  s(a^2 \tilde  s -4)}} > 0 \,.
\end{align}
We now verify that the above derivative is indeed positive. First observe that given $\tilde s > s_d $, the two roots for $\Delta(\gamma,\tilde s)$ are 
$$
\tilde \gamma_{\pm} = a \tilde s \pm 2 \sqrt{\tilde s}\,. 
$$
Clearly, $\tilde \gamma_+ > \gamma^*(a)$. Hence, $\Delta(\gamma^*(a),\tilde s) >0$ implies that $\gamma^*(a) < \tilde \gamma_-$, or equivalently $\tilde s > s_d +(2\sqrt{s_d}+2\sqrt{\tilde s})/a$. This further implies $\sqrt{\tilde s}> 2/a$. Therefore, 
\begin{align}
    -a \gamma^*(a) + a^2 \tilde s - 2 &> a^2 (s_d +(2\sqrt{s_d}+2\sqrt{\tilde s})/a) -a \gamma^*(a) - 2 \nonumber \\
    &=2a \sqrt{\tilde s} - 2 \nonumber \\
    &> 2a \frac{2}{a} - 2 > 0 \nonumber \,. 
\end{align}
Knowing that the numerator in the second term of \cref{eq:derivative} is positive, we know that \cref{eq:derivative} is positive if and only if 
$$
(-a \gamma^*(a) + a^2 \tilde s - 2)^2 > a^2 (\gamma^*(a)^2 -2a\gamma^*(a) \tilde   s +  \tilde  s(a^2 \tilde  s -4))\,,
$$
which can be verified by expanding the square on the left hand side and comparing with the right hand side directly. 

Since the derivative in \cref{eq:derivative} is positive, let us examine the limit 
\begin{align}\label{eq:ineq}
    &\lim_{s\to \infty} -\gamma^*(a)-a  s +\sqrt{\gamma^*(a)^2 -2a\gamma^*(a)  s +  s(a^2  s -4)}\nonumber \\
    = &\lim_{s\to \infty} -\gamma^*(a)-a  s  + s \sqrt{\gamma^*(a)^2 s^{-2} -2a\gamma^*(a)  s^{-1} +  a^2 -4 s^{-1}}\nonumber\\
    =& \lim_{s\to \infty} -\gamma^*(a)-a  s  +as - (\gamma^*(a) + \frac{2}{a}) + \mathcal{O}(s^{-1})\nonumber\\
    =& -2\gamma^*(a) - \frac{2}{a}\nonumber\\
    =& -2(a s_d + 2\sqrt{s_d})-\frac{2}{a} < \Re(\lambda^{(d)}_+(a,\gamma^*(a), \mS))\,. 
\end{align}
Combining \cref{eq:ineq0}, \cref{eq:derivative} and \cref{eq:ineq}, we obtain that for $1\leq j\leq d$ 
$$
\Re(\lambda^{(j)}_+(a,\gamma^*(a), \mS))  \leq \lambda^{(d)}_+(a,\gamma^*(a), \mS) = \Re(\lambda^{(d)}_+(a,\gamma^*(a), \mS))\,,
$$
which implies 
$$
\min_{\gamma>0} \max_j \Re(\lambda^{(j)}_+(a,\gamma, \mS)) \leq \max_j \Re(\lambda^{(j)}_+(a,\gamma^*(a), \mS)) = \Re(\lambda^{(d)}_+(a,\gamma^*(a), \mS))\,. 
$$
Finally, by \Cref{lemma:matrix_gamma_single} again, we have
$$
\min_{\gamma>0} \max_j \Re(\lambda^{(j)}_+(a,\gamma, \mS)) \geq \min_{\gamma>0} \Re(\lambda^{(d)}_+(a,\gamma, \mS)) = \Re(\lambda^{(d)}_+(a,\gamma^*(a), \mS))\,. 
$$
We now conclude that 
$$
\argmin_{\gamma>0} \max_j \Re(\lambda^{(j)}_+(a,\gamma, \mS)) = \gamma^*(a)\, .
$$
\end{proof}

\begin{lemma}\label{lemma:const_C1_bound}
    The constant $C_1$ in \Cref{eq:frobenius_bound } depends at most polynomially on $d$, $s_1$, $1/s_d$, i.e. $C_1 = \mathrm{poly}(d,s_1,s_d^{-1})\leq \mathrm{poly}(d,\kappa)$. 
\end{lemma}
\begin{proof}
First, we show that $C_1$ depends linearly on the dimension $d$. Let us recall the following fact from linear ODE: if $\dot x = A x$ for some constant matrix $A\in \mathbb R^{d\times d}$, with eigenvalues $\lambda_1, \ldots, \lambda_d$ and eigenvectors $v_1,\ldots,v_d$, then the solution is of the form $x(t) = \sum_{i} a_i e^{\lambda_i t} v_i$. In case there are repeated eigenvalues (e.g. $\lambda_i$) and generalized eigenvectors, the corresponding term in the sum will be replaced with some $\sum_j b_j t^{k-j} e^{\lambda_i t} v_i$ where the sum is over $j=1,\ldots,k$ and $k$ is the dimension of the generalized eigenspace associated with $\lambda_i$. Let $\mathcal{D}$ and $\mathbf{T}$ be as defined in \cref{eq:Sigma_ODE}. By our choice of $\gamma$, we know that eigenvalues of $\mathcal{D}$ are nonzero. Therefore, $\mathcal D$ is invertible. Denote by 
$$
\mathbf{Y}(t) = \begin{pmatrix}
        \Sigma_{11}(t)\\
        \Sigma_{22}(t)\\
        \dot \Sigma_{22}(t)
    \end{pmatrix} + \mathcal{D}^{-1} \mathbf{T}\,.
$$
Then \cref{eq:Sigma_ODE} reads 
\begin{equation}\label{eq:Sigma_no-const}
\frac{\dd }{\dd t} \mathbf{Y} = \mathcal{D}\mathbf{Y}\,.
\end{equation}
We follow the notation in \Cref{prop:eig_system} and use $(\lambda^{(i)}_{*}, v^{(i)}_{*})$ to represent an eigenvalue eigenvector pair of $\mathcal{D}$, for $i=1,\ldots,d$, and $*\in \{0,+,-\}$. Note that for our choice of $\gamma = a s_d + 2\sqrt{s_d}$, we have $\lambda^{(d)}_0 = \lambda^{(d)}_{\pm} $. Correspondingly, there will be generalized eigenvectors. Following the notation in \Cref{prop:eig_system}, we use $v^{(d)}_0$ to represent the eigenvector associated with $\lambda^{(d)}_0$; and we use $\eta$ and $\xi$ to represent the generalized eigenvectors associated with $\lambda^{(d)}_0$. We have already shown in \Cref{prop:eig_system} that both $\eta$ and $\xi$ are generalized eigenvector of rank 2. Therefore, the solution to \cref{eq:Sigma_no-const} takes the form 
\begin{align}
    \mathbf{Y}(t) &= \left(\sum_{i=1}^{d-1} \sum_{*\in\{0,+,-\}}  \alpha^{(i)}_* e^{\lambda^{(i)}_* t} v^{(i)}_*\right) + \alpha^{(d)}_0 e^{\lambda^{(d)}_0 t} v^{(d)}_0 + \alpha^{(d)}_- e^{\lambda^{(d)}_0 t} (t  v^{(d)}_0 + \eta) \nonumber \\
    &\qquad + 
    \alpha^{(d)}_+ e^{\lambda^{(d)}_0 t} (t v^{(d)}_0   + \xi) \,,
\end{align}
where the constants $\alpha^{(i)}_*$ are to be determined by $\mathbf{Y}(0)$. By \Cref{lemma:matrix_gamma} and our choice of $\gamma$, we have that 
$$
\max_{i\leq d } \max_{*\in \{0,+,-\}} \Re(\lambda^{(i)}_*)=\lambda^{(d)}_0 = -2as_d - 2\sqrt{s_d}\,. 
$$
Without loss of generality, consider $t\geq 1$. We have
\allowdisplaybreaks
\begin{align}\label{eq:Yt_bound}
    \| \mathbf{ Y}(t) \|^2 &= \Bigg\|\left(\sum_{i=1}^{d-1} \sum_{*\in\{0,+,-\}}  \alpha^{(i)}_* e^{\lambda^{(i)}_* t} v^{(i)}_*\right) + \alpha^{(d)}_0 e^{\lambda^{(d)}_0 t} v^{(d)}_0 + \alpha^{(d)}_- e^{\lambda^{(d)}_0 t} (t  v^{(d)}_0 + \eta) \nonumber \\
    &\qquad + 
    \alpha^{(d)}_+ e^{\lambda^{(d)}_0 t} (t v^{(d)}_0   +  \xi) \Bigg\|^2 \nonumber \\
    &= \sum_{i=1}^{d-1}\left\| \sum_{*\in\{0,+,-\}}  \alpha^{(i)}_* e^{\lambda^{(i)}_* t} v^{(i)}_*\right\|^2+ \Big\|\alpha^{(d)}_0 e^{\lambda^{(d)}_0 t} v^{(d)}_0 + \alpha^{(d)}_- e^{\lambda^{(d)}_0 t} (t  v^{(d)}_0 + \eta) \nonumber \\
    &\qquad + 
    \alpha^{(d)}_+ e^{\lambda^{(d)}_0 t} (t v^{(d)}_0  + \xi)\Big\|^2 \nonumber \\
    &\leq \sum_{i=1}^{d-1} \sum_{*\in\{0,+,-\}} 3 \left\| \alpha^{(i)}_* e^{\lambda^{(i)}_* t} v^{(i)}_*\right\|^2 + 3 \left\|\alpha^{(d)}_0 e^{\lambda^{(d)}_0 t} v^{(d)}_0 \right\|^2 + 3 \left\|\alpha^{(d)}_- e^{\lambda^{(d)}_0 t} (t  v^{(d)}_0 + \eta)\right\|^2 \nonumber \\
    &\qquad + 3 \left\| \alpha^{(d)}_+ e^{\lambda^{(d)}_0 t} (t v^{(d)}_0  + \xi) \right\|^2 \nonumber \\
    &\leq 3  t^2 e^{2\lambda^{(d)}_0 t} \Bigg[\left( \sum_{i=1}^{d-1} \sum_{*\in\{0,+,-\}} \left\| \alpha^{(i)}_* v^{(i)}_*\right\|^2\right)  + \left\|  v^{(d)}_0 \right\|^2 \big((\alpha^{(d)}_0)^2 + 2(\alpha^{(d)}_-)^2 + 2(\alpha^{(d)}_+)^2\big) \nonumber \\
    &\qquad +  2\left\| \eta \right\|^2 (\alpha^{(d)}_- )^2  +  2 \left\|\xi\right\|^2 (\alpha^{(d)}_+ )^2  \Bigg] \nonumber \\
    &\leq 6  t^2 e^{2\lambda^{(d)}_0 t} \Bigg[\left( \sum_{i=1}^{d-1} \sum_{*\in\{0,+,-\}} \left\| \alpha^{(i)}_* v^{(i)}_*\right\|^2\right)  + \left\|  v^{(d)}_0 \right\|^2 \big((\alpha^{(d)}_0)^2 + (\alpha^{(d)}_-)^2 + (\alpha^{(d)}_+)^2\big) \nonumber \\
    &\qquad +  \left\| \eta \right\|^2  (\alpha^{(d)}_- )^2  +  \left\|\xi\right\|^2 (\alpha^{(d)}_+ )^2  \Bigg] \nonumber \\
    &\leq 6  t^2 e^{2\lambda^{(d)}_0 t} \Bigg[\left( \sum_{i=1}^{d-1} \sum_{*\in\{0,+,-\}} \left\| \alpha^{(i)}_* v^{(i)}_*\right\|^2\right)+ \Big( \left\|\alpha^{(d)}_+  v^{(d)}_0 \right\|^2 + \left\|\alpha^{(d)}_- \eta \right\|^2 + \left\|\alpha^{(d)}_+ \xi\right\|^2 \Big)\nonumber \\
    &\qquad \left(1 + \frac{ \Big\|v^{(d)}_0\Big\|^2}{\|\xi\|^2} + \frac{\Big\|v^{(d)}_0\Big\|^2}{\|\eta\|^2}\right) \Bigg]\nonumber \\
    &\leq 6  t^2 e^{2\lambda^{(d)}_0 t} \left(1 + \frac{ \Big\|v^{(d)}_0\Big\|^2}{\|\xi\|^2} + \frac{\Big\|v^{(d)}_0\Big\|^2}{\|\eta\|^2}\right) \Bigg[ \sum_{i=1}^{d-1} \sum_{*\in\{0,+,-\}} \left\| \alpha^{(i)}_* v^{(i)}_*\right\|^2 \nonumber \\
    &\qquad + \left\|\alpha^{(d)}_+  v^{(d)}_0 \right\|^2 + \left\|\alpha^{(d)}_- \eta \right\|^2 + \left\|\alpha^{(d)}_+ \xi\right\|^2 \Bigg]\,.
\end{align}
$$
\| \mathbf{ Y}(0) \|^2 = \sum_{i=1}^{d-1}\left\| \sum_{*\in\{0,+,-\}}  \alpha^{(i)}_* v^{(i)}_*\right\|^2+ \Big\|\alpha^{(d)}_0  v^{(d)}_0 + \alpha^{(d)}_-  \eta+ 
    \alpha^{(d)}_+   \xi\Big\|^2\,.
$$
Denote by $\mathbf{Y}(0)^{(i)}$ the projection of $\mathbf{Y}(0)$ onto the subspace $\Phi_i = \mathrm{Span}( \{v^{(i)}_0,v^{(i)}_+,v^{(i)}_-\})$. And accordingly, $\Phi_d = \mathrm{Span}(\{v^{(d)}_0,\eta,\xi\})$. By \Cref{cor:orthogonal}, we know that $\Phi_i$ is orthogonal to $\Phi_j$ for $i\neq j$. Therefore, $|\alpha^{(i)}_*|$ depends on the inverse of the Gram matrix of $\{v^{(i)}_0,v^{(i)}_+,v^{(i)}_-\}$ as well as $\|\mathbf{Y}(0)^{(i)}\|$. This inverse Gram matrix can be computed analytically since it is a 3 by 3 matrix for each $1\leq i \leq d$. However, the exact computation does not add more insights to the proof and we will not include the computation. Since each eigenvector and generalized eigenvector depends on $\{s_1,\ldots\,s_d,s_1^{-1},\ldots\,s_d^{-1}\}$ polynomially, we know that the inverse of the Gram matrix also also depends on $\{s_1,\ldots\,s_d,s_1^{-1},\ldots\,s_d^{-1}\}$ polynomially. From \cref{eq:Yt_bound}, we conclude that 
$$
\| \mathbf{ Y}(t) \|^2 = \mathcal{O}\big(t^2 e^{2\lambda^{(d)}_0 t} d^2 \cdot \mathrm{poly}(s_1,s_d^{-1})\big) =\mathcal{O}\big(t^2 e^{2\lambda^{(d)}_0 t} d^2 \cdot \mathrm{poly}(\kappa)\big) \,.
$$
\end{proof}

\begin{lemma}\label{lemma:X_zero}
    Suppose $X\in \mathbb S^n$ satisfies $X = AXA^T$ for some $A\in \mathbb R^n$. If all eigenvalues of $A$ has absolute value less than 1, then $X$ is the zero matrix. 
\end{lemma}
\begin{proof}
    Let us first assume that $A^T$ is diagonalizable: $A^T = QDQ^{-1}$, where $D$ is a diagonal matrix of eigenvalues $d_1, \ldots, d_n$, and the columns of $Q$ contains the eigenvectors $q_1, \ldots ,q_n$. Then it follows that 
    $$
    |q_i^T X q_j| = |d_i d_j| |q_i^T X q_j|\,. 
    $$
    This implies $|q_i^T X q_j|=0$ for all $1\leq i,j \leq n$, since $|d_i d_j|<1$ by assumption. Now suppose that $A$ has some generalized eigenvalues. Without loss of generality, assume that $d_{n-1}$ is a generalized eigenvalue such that $A^T q_{n-1} = d_{n-1}q_{n-1}$ and $A^T q_n = d_{n-1} q_n + q_{n-1}$. Let $q_i$ be an eigenvector. Then we still have $q_i^T X q_{n-1}=0$ as before. And 
    $$
    |q_i^T X q_n| = |d_i d_{n-1}q_i^T X q_{n} + d_i q_i^T X q_{n-1}| = |d_i d_{n-1}q_i^T X q_{n} | = |d_i d_{n-1}|\, |q_i^T X q_{n} |\,.
    $$
    Again this implies $|q_i^T X q_{n} |=0$. The case where $d_{n-1}$ has algebraic multiplicity greater than 2 or $q_i$ is a generalized eigenvector can be proved in a similar fashion. Therefore, we have shown that if $A$ has Jordan decomposition $A = PJP^{-1}$, then $q_i^T X q_j = 0$ where $q_i$ and $q_j$ are the $i$-th and $j$-th column of $P$. Equivalently, we have $P^T X P = 0$. This proves that $X=0$.  
\end{proof}

\begin{corollary}\label{cor:uniqueness}
Suppose $X, Y\in \mathbb S^n$ satisfy $X = AXA^T +B$, $Y=AYA^T+B$ for some $B \in \mathbb S^n$. If all eigenvalues of $A$ have absolute value less than 1, then $X=Y$. 
\end{corollary}
\begin{proof}
 The proof follows by \Cref{lemma:X_zero} and that $X-Y = A(X-Y)A^T$. 
\end{proof}
Taking inspiration from system of linear ODE, we have the following lemma regarding the solution to the iteration $X_{k+1}=AX_k A^T$. 
\begin{lemma}\label{lemma:X_k_Jordan}
    Let $A \in \mathbb R^{n\times n}$ be given by $A = \bI - h \tilde G$ for some $\tilde G\in \mathbb R^{n\times n}$, $h >0$. Suppose $\tilde G$ has Jordan decomposition $\tilde G = PJP^{-1}$. And consider the iteration $X_{k+1} = A X_k A^T$. If $q_i$ is an eigenvector of $\tilde G$ with associated eigenvalue $d_i$ and $X_0=q_i q_i^T$, then $X_k = (1-h d_i)^{2k} X_0$. Moreover, if $q_i$ is a generalized eigenvector of $\tilde G$ of algebraic multiplicity 2, i.e. $\tilde Gq_i = d_j q_i + q_j$ for some eigenvector $q_j$ and eigenvalue $d_j$, and $X_0 = q_i q_i^T$, then $X_k = \big((1-hd_j)^k q_i - kh(1-hd_j)^{k-1} q_j\big)\big((1-hd_j)^k q_i - kh(1-hd_j)^{k-1} q_j\big)^T $
\end{lemma}

\begin{lemma}\label{lemma:discrete_eig}
    The eigenvalues of $\mG$ in \cref{eq:A} are given by the following 
\begin{align}\label{eq:discrete_eigvalue_approx}
    \tilde \lambda^{(i)}_{\pm} = h\,\frac{(a s_i + \gamma) \pm 
    \sqrt{(a s_i-\gamma)^2 - 4s_i}}{2}\,.
\end{align}
\end{lemma}
\begin{proof}
    The proof follows by a direct computation. 
\end{proof}
\begin{lemma}\label{lemma:a}
Consider $\gamma= \gamma^* = as_d + 2\sqrt{s_d}$. Let $s > s_d$. Then $a 
\leq \frac{2}{\sqrt{s}-\sqrt{s_d}}$ if and only if $(a s-\gamma^*)^2 - 4s \leq 0$. 
\end{lemma}
\begin{proof}
Multiplying by $s-s_d$, we obtain 
$$
a \leq \frac{2}{\sqrt{s}-\sqrt{s_d}} \Longleftrightarrow a(s-s_d) \leq 2\sqrt{s} + 2\sqrt{s_d} \Longleftrightarrow as - 2\sqrt{s} \leq  \gamma^*\,. 
$$
And it is straightforward to verify that $2\sqrt{s}> -as + \gamma^*$ always holds. Squaring on both hand sides completes the proof. 
\end{proof}
\begin{lemma}\label{lemma:discrete_eig_bound}
    Consider $\tilde \lambda^{(i)}_{\pm}$ given by \cref{eq:discrete_eigvalue_approx}. Suppose $a\geq \frac{2}{\sqrt{s_{1}}-\sqrt{s_d}}$.  If the step size $h$ satisfies $0<h\leq 1/(a s_1 + \gamma)$, and $\gamma = \gamma^* = as_d + 2\sqrt{s_d}$, then 
    $$
    \max_i |1-\tilde\lambda^{(i)}_{\pm}| \leq 1 - \frac{h}{2}(as_d + \sqrt{s_d})\,.
    $$
\end{lemma}
\begin{proof}
Observe that the eigenvalues given in \cref{eq:discrete_eigvalue_approx} is almost the same as the eigenvalues given in \cref{eq:lambda_pm_hd} except for an extra factor of $h/2$. This allows us to use previous lemma regarding the eigenvalues from \cref{eq:lambda_pm_hd}. We consider two cases. Define
$$
j = \inf\left \{n: a \leq \frac{2}{\sqrt{s_{n}}-\sqrt{s_d}}\right\}\,.
$$
\textbf{Case 1:} Consider $i\leq j-1$ (if $j=1$, we directly consider Case 2). Then $a \geq \frac{2}{\sqrt{s_{i}}-\sqrt{s_d}}$. By \Cref{lemma:a} and our assumption on $a$, we have $(a s_i-\gamma^*)^2 - 4s_i \geq 0$. Then, one can verify that $0<h\leq \frac{1}{as_1 + \gamma^*}$ is a sufficient condition for $1- \tilde\lambda^{(i)}_{\pm} >0$. Indeed, we compute 
\begin{align}
    \tilde\lambda^{(i)}_{\pm} &\leq \frac{1}{as_1 + \gamma^*} \frac{(a s_i + \gamma^*) + 
    \sqrt{(a s_i-\gamma^*)^2 - 4s_i}}{2} \nonumber \\
    &\leq \frac{1}{as_1 + \gamma^*} \frac{(a s_i + \gamma^*) + 
    \sqrt{(a s_i+\gamma^*)^2}}{2} \nonumber \\
    &= \frac{as_i + \gamma^* }{as_1 + \gamma^*} \nonumber \\
    &\leq 1\,. 
\end{align}
Moreover, we clearly have $\tilde\lambda^{(i)}_{\pm} >0$. Therefore, $|1-\tilde\lambda^{(i)}_{\pm}| \leq 1$. On the other hand, by \cref{eq:derivative} and \cref{eq:ineq}, we have that 
\begin{align*}
     \tilde\lambda^{(i)}_{\pm} & \geq  \lim_{s \to \infty}\frac{h}{2} \Big((a s + \gamma^*) + 
    \sqrt{(a s-\gamma^*)^2 - 4s}\Big) \\
    &= \frac{h}{2} \Big(2\gamma^* + \frac{2}{a}\Big) \\
    &\geq  h\gamma^* \,.
\end{align*}
Therefore, 
$$
 \max_{i\leq j-1} |1-\tilde\lambda^{(i)}_{\pm}| \leq 1-h(as_d + 2\sqrt{s_d}) \,. 
$$

\noindent\textbf{Case 2:} Consider $i\geq j$. Note that for a complex number $z = z_1 + i z_2$ and $h>0$, we have that 
$$
|1-hz|^2 = (1-hz_1)^2 + h^2 z_2^2 \leq 1-hz_1\leq (1-hz_1/2)^2\,,
$$
where the first inequality holds if and only if $h\leq z_1/(z_1^2 + z_2^2)$. Therefore, we have 
\begin{equation*}
    |1-\tilde\lambda^{(i)}_{\pm}|^2 \leq  \Big(1-\frac{\Re(\tilde\lambda^{(i)}_{\pm})}{2}\Big)^2\,, 
\end{equation*}
if 
$$
h \leq \frac{2(as_i+\gamma^*)}{(as_i+\gamma^*)^2 + 4s_i-(as_i-\gamma^*)^2} = \frac{as_i + \gamma^*}{2as_i \gamma^* + 2s_i } \,. 
$$
We now verify that $h\leq 1/(a s_1 + \gamma^*)$ is a sufficient condition. We have 
$$
\frac{1}{a s_1 + \gamma^*} \leq \frac{as_i + \gamma^*}{2as_i \gamma^* + 2s_i } \Longleftrightarrow a^2 s_1 s_i + \gamma^* a s_1 + (\gamma^*)^2 \geq as_i \gamma^* + 2s_i \,.
$$
By \Cref{lemma:a}, we have that 
\begin{align*}
    a^2 s_1^2 + (\gamma^*)^2 - 2as_1 \gamma^* &\geq 4s_1 \\
    a^2 s_1 +\frac{(\gamma^*)^2}{s_1} - 2a \gamma^* &\geq 4 \\
    a^2 s_1 &\geq 4 + 2a \gamma^* - \frac{(\gamma^*)^2 }{s_1}\,.
\end{align*}
Then 
\begin{align*}
    a^2 s_1 s_i + \gamma^* a s_1 + (\gamma^*)^2 & \geq s_i \Big(4 + 2a \gamma^* - \frac{(\gamma^*)^2 }{s_1}\Big) + \gamma^* a s_1 + (\gamma^*)^2 \\
    &= 4s_i + 2a\gamma^* s_i - \frac{s_i (\gamma^*)^2}{s_1} + \gamma^* a s_1 + (\gamma^*)^2 \\
    &\geq 4s_i + 2a\gamma^* s_i + \gamma^* a s_1 \\
    &> as_i \gamma^* + 2s_i\,. 
\end{align*}
This shows that $h\leq 1/(a s_1 + \gamma^*)$ is sufficient. By \Cref{lemma:matrix_gamma} and our choice of $h$, we obtain that 
$$
\max_{i\geq j} |1-\tilde\lambda^{(i)}_{\pm}| < \max_i 1-\frac{\Re(\tilde\lambda^{(i)}_{\pm})}{2} \leq 1 - \frac{h}{2}(as_d + \sqrt{s_d})\,.
$$
Combining the two cases, we complete the proof.
\end{proof}

\begin{lemma}\label{lemma:h_size_ul}
    Consider $\tilde \lambda^{(i)}_{\pm}$ given by \cref{eq:discrete_eigvalue_approx}. Suppose $a=0$ and $\gamma =\gamma^*= 2\sqrt{s_d}$. Then  
    $$
    \max_i |1-\tilde\lambda^{(i)}_{\pm}| \leq 1\,,
    $$
    if and only if $h\leq 2\sqrt{s_d}/s_1$. 
\end{lemma}
\begin{proof}
    We directly compute 
    \begin{align*}
        |1-\tilde\lambda^{(i)}_{\pm}|^2 \leq 1 &\Longleftrightarrow |1-h\sqrt{s_d}\mp h\sqrt{s_d-s_i}|^2 \leq 1 \\
        &\Longleftrightarrow 1-2h\sqrt{s_d} + h^2 s_i \leq 1\\
        &\Longleftrightarrow h\leq 2\sqrt{s_d}/s_1\,. 
     \end{align*}
\end{proof}

\begin{theorem}\label{thm:conv_rate}
Consider the iteration given in \Cref{cor:Yn}. Suppose $a\geq \frac{2}{\sqrt{s_{1}}-\sqrt{s_d}}$. We choose $\gamma = \gamma^* = as_d + 2\sqrt{s_d}$ and $0<h\leq 1/(a s_1 + \gamma^*)$. Then for $k\geq 1/h$ we have $\|Y_k\|_\mathrm{F} \leq \widetilde C h^2k^2 (1 - \frac{h}{2}(as_d + \sqrt{s_d})^{2k-2} $, where the constant $\widetilde C =  d^2 \cdot \mathcal{O}(\mathrm{poly}(\kappa))$.   
\end{theorem}
\begin{proof}
Let us denote by $A = PJP^{-1}$ the Jordan decomposition of $A$. Then we know from \cref{eq:discrete_eigvalue_approx} that $A$ has precisely $2d-1$ eigenvectors and one generalized eigenvector of algebraic multiplicity 2. Let $q^{(i)}_{\pm},\ldots,q^{(d-1)}_{\pm}$ be the eigenvectors with associated eigenvalues $\lambda^{(i)}_{\pm}=1-\tilde \lambda^{(i)}_{\pm}$, where $\tilde \lambda^{(i)}_{\pm}$ are from \cref{eq:discrete_eigvalue_approx}. With $\gamma=\gamma^*$, one has that $\tilde \lambda^{(d)}_{+}=\tilde \lambda^{(d)}_{-}$ is a generalized eigenvalue. Abusing notation, let us use $q^{(d)}_+$ to represent the eigenvector and $q^{(d)}_-$ to represent the generalized eigenvector of $\lambda^{(d)}_{-}=\lambda^{(d)}_{+}$. This means 
$$
A q^{(d)}_+ = \lambda^{(d)}_+ q^{(d)}_+ \,,\qquad A q^{(d)}_- = \lambda^{(d)}_- q^{(d)}_- +  q^{(d)}_+\,. 
$$
We can express $Y_0$ by a basis representation
$$
Y_0 = \sum_{\star,\diamond\in\{\pm\}} \sum_{i,j\leq d} \alpha^{i,j}_{\star,\diamond}\, q^{(i)}_{\star} (q^{(j)}_{\diamond})^T\,.
$$
Then using \Cref{lemma:X_k_Jordan}, we have that for $k\geq 1/h$, 
\begin{align}
    \|Y_k\|_\mathrm{F} &\leq 4d^2 h^2 k^2 \max_{i}|\lambda^{(i)}_{\pm}|^{2k-2} \max_{i,j,\star,\diamond}|\alpha^{i,j}_{\star,\diamond}|\|q^{(i)}_{\star} q^{(j)}_{\diamond}\|_\mathrm{F}   \nonumber \\
    &\leq 4d^2 h^2 k^2 \left(1 - \frac{h}{2}(as_d + \sqrt{s_d})\right)^{2k-2}\max_{i,j,\star,\diamond}|\alpha^{i,j}_{\star,\diamond}|\|q^{(i)}_{\star} q^{(j)}_{\diamond}\|_\mathrm{F}\,.
\end{align}
The second inequality is due to \Cref{lemma:discrete_eig_bound}. The maximum in the above is over $1\leq i,j\leq d$ and $\star,\diamond \in\{\pm\}$. It remains to show that $\max_{i,j,\star,\diamond}|\alpha^{i,j}_{\star,\diamond}|\|q^{(i)}_{\star} q^{(j)}_{\diamond}\|_\mathrm{F} = \mathcal{O}(\mathrm{poly}(\kappa))$. Note that $\mA$ in \Cref{cor:Yn} can be written as $\mA = \bI - h \tilde G$ where $\tilde G$ does not depend on $h$ when taking the first order approximation as in \Cref{lemma:discrete_eig}.  The rest of the argument is very similar to the proof of \Cref{lemma:const_C1_bound} which we will not present due to brevity. We conclude that 
\begin{align*}
    \|Y_k\|_\mathrm{F} &\leq d^2 h^2 k^2 \left(1 - \frac{h}{2}(as_d + \sqrt{s_d})\right)^{2k-2} \mathcal{O}(\mathrm{poly}(\kappa)) \\
    &= \widetilde C h^2 k^2 \left(1 - \frac{h}{2}(as_d + \sqrt{s_d})\right)^{2k-2} \,. 
\end{align*}
\end{proof}

\begin{lemma}\label{lemma:fixed_point}
A solution to the fixed point equation $\mY^*=\mA \mY^* \mA^T + \mL \mL^T $ where $\mA$ and $\mL$ are given in \Cref{prop:discrete_update_cheng}, is given by 
$$
\mY^* = \begin{pmatrix}
    Y^*_{11} & Y^*_{12}\\
    Y^*_{12} & Y^*_{22}
\end{pmatrix}\,,
$$
where $Y^*_{ij} \in \mathbb R^d$ are diagonal matrices. And the diagonal elements of $Y^*_{ij}$ are given by 
\begin{align}
    Y^*_{11,i} &= \frac{1}{s_i} \Big( 1-\frac{hs_i (4+(h+a(h \gamma-2))(hs_i - \gamma + as_i(h\gamma -1)))}{(hs_i - \gamma + as_i(h\gamma-1))(4+h(hs_i -2\gamma + as_i(h\gamma -2)))}\Big)\,,\label{eq:x_element}  \\
    Y^*_{12,i} &= \frac{2h(as_i - \gamma)}{(hs_i - \gamma + as_i(h\gamma-1))(4+h(hs_i -2\gamma + as_i(h\gamma -2)))}\,, \\
    Y^*_{22,i} &= \frac{-4\gamma -2as_i (2+h(hs_i-3\gamma + as_i(h \gamma -1)))}{(hs_i - \gamma + as_i(h\gamma-1))(4+h(hs_i -2\gamma + as_i(h\gamma -2)))}\,.
\end{align}
\end{lemma}

\section{Postponed proofs}\label{sec:post_proof}
\begin{proof}[proof of \Cref{prop:decomposition}]
We directly plug \cref{eq:FPK_Gamma} into \cref{eq:FPK_decomp} and verify that we recover \cref{eq:FPK_general}.  
\begin{align*}
    &\nabla \cdot \Big(\rho \,\mathrm{sym}(\bQ) \nabla \log \frac{\rho}{\Pi} \Big) + \nabla \cdot \big(\rho (\mathrm{sym}(\bQ) \nabla \log(\Pi) + \bQ \nabla H)\big) \\
    &= \mathrm{sym}(\bQ) : \nabla^2 \rho + \nabla \rho\, \mathrm{sym}(\bQ) \nabla H + \rho \, \mathrm{sym}(\bQ): \nabla^2 H + \nabla \rho \,\mathrm{sym}(\bQ) \nabla \log(\Pi) \\
    &\qquad + \rho \, \mathrm{sym}(\bQ): \nabla^2 \log(\Pi) + \nabla \cdot \big(\rho \bQ \nabla H\big) \\
    &= \mathrm{sym}(\bQ) : \nabla^2 \rho  + \nabla \cdot \big(\rho \bQ \nabla H\big) \\ 
    &=  \nabla \cdot (\bQ \nabla H \rho) + \sum_{i,j=1}^{2d} \frac{\partial^2}{\partial X_i \partial X_j} (Q_{ij} \rho)\,, 
\end{align*}
where we denote by $\bA:\mathbf B = \sum_{i,j=1}^{2d} A_{ij} B_{ij}$ for $\bA, \mathbf B \in \mathbb R^{d\times d}$. We have used $\nabla \log(\Pi) = -\nabla H$ and $\nabla^2 \log(\Pi) = -\nabla^2 H$ to get the second equality. 
\end{proof}

\begin{proof}[proof of \Cref{prop:pi}]
    We just need to verify that when $\rho(\vX,t) = \Pi(\vX)$, we have  $\frac{\partial \rho}{\partial t} = 0$. It is clear that when $\rho(\vX,t) = \Pi(\vX)$, the first term on the right hand side of \cref{eq:FPK_decomp} is 0, since $\nabla \log(\frac{\rho}{\Pi})=0$. For the second term, let us use \cref{eq:FPK_Gamma} to get 
    \begin{align}
        \nabla \cdot (\Pi \Gamma) &= \nabla \cdot \big(\Pi \bQ \nabla H - \Pi \, \mathrm{sym}(\bQ) \nabla \log(\Pi) \big) \nonumber \\
        &= \nabla \Pi \bQ \nabla H + \Pi \bQ : \nabla^2 H + \nabla \Pi \,\mathrm{sym}(\bQ) \nabla \log(\Pi) + \Pi \,\mathrm{sym}(\bQ) : \nabla^2 \log(\Pi) \nonumber \\
        &= -\Pi \nabla H \bQ \nabla H + \Pi \bQ : \nabla^2 H + \Pi \nabla H  \mathrm{sym}(\bQ) \nabla H +  \Pi \,\mathrm{sym}(\bQ) : \nabla^2 \log(\Pi) \nonumber \\
        &= \Pi \bQ : \nabla^2 H +  \Pi \,\mathrm{sym}(\bQ) : \nabla^2 \log(\Pi) \nonumber \\
        &= \Pi \bQ : \nabla^2 H - \Pi \,\mathrm{sym}(\bQ) : \nabla^2 H \nonumber \\
        &= 0 \nonumber  \,,
    \end{align}
      
    We have used $\nabla \Pi = -\Pi \nabla H$ to get the third equality. And we used $\nabla^2 \log(\Pi) = -\nabla^2 H$ to get the fifth equality. This proves that when $\rho = \Pi$, we indeed have 
    $$
     \frac{\partial \rho}{\partial t}\Big\vert_{\rho=\Pi} = \nabla \cdot \Big(\Pi \,\mathrm{sym}(\bQ) \nabla \log \frac{\Pi}{\Pi} \Big) + \nabla \cdot (\Pi \Gamma) = 0 + 0 = 0 \,. 
    $$
\end{proof}
\begin{proof}[proof of \Cref{prop:Sigma_ODE}]
With our choice of $H$, \cref{eq:SDE} is a multidimensional OU process. And since $\bX_0$ follows a Gaussian distribution, it shows that $\bX_t$ will also be a Gaussian distribution. It is well known that the solution to \cref{eq:SDE} with $H$ given by \cref{eq:quadratic H} is 
$$
\vX_t = e^{-t \bQ \widetilde \Sigma^{-1}} \vX_0 + \int_0^t  e^{-(t-\tau) \bQ \widetilde \Sigma^{-1}} \sqrt{2\,\mathrm{sym}(\bQ)} \,d \mB_{\tau}\,.
$$
The mean of $\bX_t$ is given by 
$$
\mathbb E \bX_t = e^{-t \bQ \widetilde \Sigma^{-1} }  \mathbb E \bX_0 = 0.
$$
We can compute the covariance $\Sigma(t)$ of $\bX_t$. Since $\bX_t$ has zero mean, we obtain the following using Ito's isometry 
\begin{align}\label{eq:cov_soln}
    \Sigma(t) = \mathbb E \bX_t \bX_t^T &= 2\int_0^t e^{-(t-\tau) \bQ \widetilde \Sigma^{-1}}  \mathrm{sym}(\bQ) \Big(e^{-(t-\tau) \bQ \widetilde \Sigma^{-1}} \Big)^T d \tau  + \mathbb E \bX_0 \bX_0^T\,. 
\end{align}
From the above expression, $\Sigma(t)$ is clearly well-defined, symmetric, positive definite for all $t>0$. We proceed by differentiating $\Sigma(t)$
\begin{align*}
    \dot \Sigma(t) &= 2 \frac{d}{dt} \int_0^t e^{-(t-\tau) \bQ \widetilde \Sigma^{-1}}  \mathrm{sym}(\bQ) \Big(e^{-(t-\tau) \bQ \widetilde \Sigma^{-1}} \Big)^T d \tau \\ 
    &= 2 \,\mathrm{sym}(\bQ) + \int_0^t \frac{d}{dt} e^{-(t-\tau) \bQ \widetilde \Sigma^{-1}}  \mathrm{sym}(\bQ) \Big(e^{-(t-\tau) \bQ \widetilde \Sigma^{-1}} \Big)^T d \tau \\
    &= 2\, \mathrm{sym}(\bQ)  - \bQ \widetilde \Sigma^{-1} \Sigma(t) - \Sigma(t) \widetilde \Sigma^{-1} \bQ^T \\
    &= 2\,\mathrm{sym}(\bQ(\bI-\widetilde\Sigma^{-1}\Sigma )) \,. 
\end{align*}
This finishes the proof. 
\end{proof}
\end{appendices}

\end{document}